\documentclass{article}

\usepackage{arxiv}

\usepackage[utf8]{inputenc} 
\usepackage[T1]{fontenc}    
\usepackage{hyperref}       
\usepackage{url}            
\usepackage{booktabs}       
\usepackage{amsfonts}       
\usepackage{nicefrac}       
\usepackage{microtype}      
\usepackage{lipsum}

\usepackage{morefloats}
\usepackage{color}
\usepackage{multirow}
\usepackage{latexsym}
\usepackage{amsmath,amssymb,amsthm,graphicx,mathrsfs}
\usepackage{dsfont}
\usepackage{enumitem}
\usepackage{extarrows}
\usepackage{hyperref}

\usepackage{morefloats}
\usepackage{dsfont}

\usepackage[title]{appendix}

\usepackage{rotating}
\newtheoremstyle{thm}
{9pt}
{9pt}
{\itshape}
{}
{\bfseries}
{.}
{ }
{}
\theoremstyle{thm}

\newtheorem{theorem}{Theorem}[section]
\newtheorem{lemma}[theorem]{Lemma}
\newtheorem{corollary}[theorem]{Corollary}

\newtheoremstyle{def}
{9pt}
{9pt}
{}
{}
{\bfseries}
{.}
{ }
{}
\theoremstyle{def}

\def\R{\mathbb{R}}

\def\N{\mathbb{N}}
\def\P{\mathbb{P}}

\def\PP{\mathbb{P}}
\def\BE{\mathbb{E}}

\def\dint{{\rm d}}
\DeclareMathOperator{\Var}{Var}

\DeclareMathOperator{\Vol}{Vol}

\renewcommand{\footnoterule}{%
	\kern -3.5pt
	\hrule width \textwidth height 1pt
	\kern 3.5pt
}

\makeatletter
\def\blfootnote{\xdef\@thefnmark{}\@footnotetext}
\makeatother

\title{Testing multivariate uniformity based on random geometric graphs}

\author{Bruno Ebner\\
 Institute of Stochastics, \\
Karlsruhe Institute of Technology (KIT), \\
Englerstr. 2, D-76133 Karlsruhe, Germany \\
\texttt{Bruno.Ebner@kit.edu}\\
\And
Franz Nestmann\\
Institute of Stochastics, \\
Karlsruhe Institute of Technology (KIT), \\
Englerstr. 2, D-76133 Karlsruhe, Germany \\
\texttt{Franz.Nestmann2@kit.edu}\\
\And
Matthias Schulte\\
Heriot-Watt University,\\
School of Mathematical and Computer Sciences, \\
Edinburgh EH14 4AS, UK\\
\texttt{M.Schulte@hw.ac.uk}
}

\begin{document}

\date{\today}
\maketitle

\blfootnote{ {\em MSC 2010 subject
classifications.} Primary 62G10 Secondary 62G20, 60D05}
\blfootnote{
{\em Key words and phrases} Multivariate goodness-of-fit test; uniform distribution; random geometric graph; Gilbert graph; $U$-statistics, contiguous alternatives}

\begin{abstract}
We present new families of goodness-of-fit tests of uniformity on a full-dimensional set $W\subset\R^d$ based on statistics related to edge lengths of random geometric graphs. Asymptotic normality of these statistics is proven under the null hypothesis as well as under fixed alternatives. The derived tests are consistent and their behaviour for some contiguous alternatives can be controlled. A simulation study suggests that the procedures can compete with or are better than established goodness-of-fit tests. We show with a real data example that the new tests can detect non-uniformity of a small sample data set, where most of the competitors fail.
\end{abstract}

\section{Introduction}

The analysis of point patterns in a given study area is of particular interest in a wide variety of fields, such as astronomy (e.g.\ occurrence of high energetic events in a sky map), biology (e.g.\ locations of sightings of threatened species) or geology (e.g.\ locations of raw materials). The concept of uniformity of the observations stands for the absence of structure in the data. Thus, testing uniformity of random vectors is a natural starting point for serious statistical inference involving any cluster analysis or multimodality assumption. To be specific, let $n\in\N$ and
\begin{equation*}
{\mathscr X}_n := \{X_1,\ldots,X_n\}
\end{equation*}
be the data set, where $X_1,\ldots,X_n$ are independent identically distributed (i.i.d.) random vectors taking values in a given measurable set $W\subset\R^d$, $d\ge1$, of positive finite volume, called the observation window. Without loss of generality we assume that $\Vol(W)=1$. We want to test the {\em null hypothesis}
\begin{equation}\label{H0}
H_0:\,X\sim {\cal U}(W)
\end{equation}
with $X$ being an independent copy of $X_1$ and ${\cal U}(W)$ denoting the uniform distribution on $W$ against general alternatives. This situation also arises in the investigation of pseudo random number generators, see e.g.\ \cite[Section 3.3]{knuth}. Testing if i.i.d.\ random vectors in $\R^d$ follow a given absolutely continuous distribution is, by the Rosenblatt transformation, see \cite{22}, theoretically equivalent to testing uniformity on the $d$-dimensional unit cube $[0,1]^d$, although this transformation is hard to compute in many cases. The problem of testing uniformity has been investigated in classical papers in the univariate case, see \cite{12} for a survey and \cite{Baringhaus2018} for a recent article, and, hitherto far less studied, in the multivariate setting, see \cite{BartoszynskiEtAl,03,05,24,02,04,SR:2005,20,Yang2017}, for which an empirical study was conducted in \cite{08}. The cited methods include classical goodness-of-fit testing approaches as the Kolmogorov-Smirnov test, see \cite{02}, nearest neighbour concepts, see \cite{24} and the references therein, the distances of data points to the boundary of the observation window, see \cite{05}, or the volume of the largest ball that can be placed in the observation window and does not cover any data point, see \cite{03}. The related problem of testing for complete spatial randomness of a point pattern (i.e., the points are a realisation of a homogeneous Poisson point process) is also of ongoing interest, see e.g.\ monographs like \cite{spatstat,25} or the recent publications \cite{26,27}.

\begin{figure}[h!]
\centering
\includegraphics[scale=0.7]{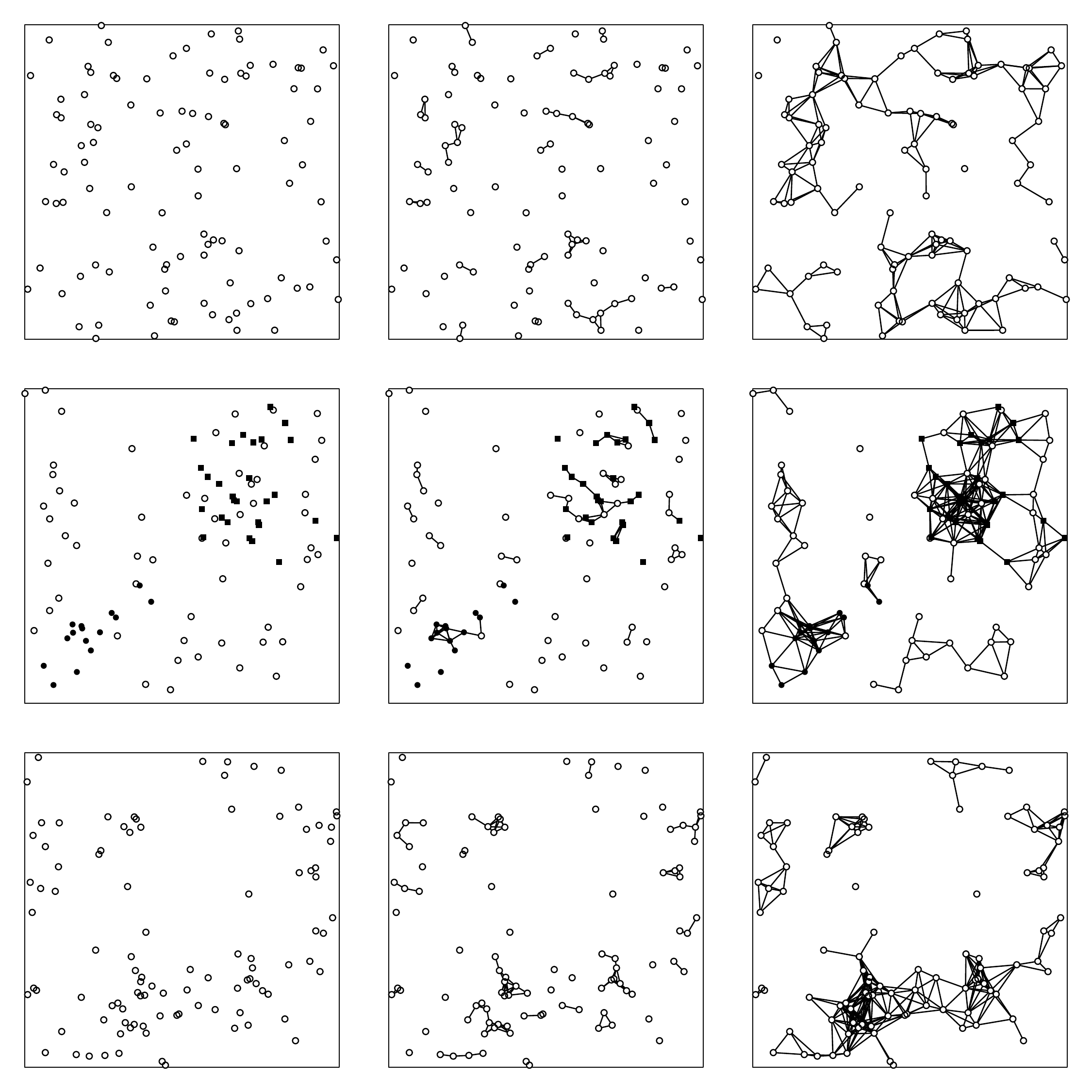}
\caption{Realisations of uniform data in $W=[0,1]^2$ (first row), the CON alternative (second row) and the CLU alternative (third row). Point data (first column), $RGG({\mathscr X}_n, 0.03)$ (second column) and $RGG({\mathscr X}_n, 0.06)$ (third column), $n=100$.}\label{fig:9plots}
\end{figure}

We approach the testing problem \eqref{H0} by examining the local properties of the data by means of random graphs. Using random graphs for testing uniformity is a known but not widely used concept, see \cite{Godehardt1997,21,08}. Our new approach is to consider statistics of the random geometric graph $RGG({\mathscr X}_n, r_n)$, $r_n>0$: It has the realisations of the random vectors in ${\mathscr X}_n$ as vertices, and any two distinct vertices $x,y\in{\mathscr X}_n$ are connected by an edge if $\|x-y\|\le r_n$, where $\|\cdot\|$ stands for the Euclidean norm. This random graph model was introduced by Gilbert for an underlying Poisson point process in \cite{Gilbert1961} and is thus also called Gilbert graph. For further details see \cite{09} and the references cited therein. Figure \ref{fig:9plots} provides a visualisation of different point data and selected random geometric graphs. For definitions of the CLU and CON alternatives we refer to Section \ref{Simu}.

Our test statistics are related to the edge lengths of $RGG({\mathscr X}_n, r_n)$ and are defined by
$$
L_n(\beta):=\frac12\sum_{(x,y)\in{\mathscr X}^2_{n,\neq}}\mathbf{1}\{\|x-y\|\le r_n\}\|x-y\|^\beta,\quad \beta\in\R.
$$
Here $\sum_{(x,y)\in{\mathscr X}^2_{n,\neq}}$ stands for the sum over all pairs of distinct points of ${\mathscr X}_n$ (such sums are called $U$-statistics), and $\mathbf{1}\{\cdot\}$ is the indicator function. Notice that $L_n(0)$ counts the number of edges and $L_n(1)$ is the total edge length of $RGG({\mathscr X}_n, r_n)$. These statistics differ from nearest neighbour methods, see e.g.\ \cite{BiauDevroye,24} and the references therein, as such that they rely on all interpoint distances not exceeding $r_n$, whereas nearest neighbour methods take only distances between points and their $k$-nearest neighbours into account. In order to analyse point processes in spatial statistics, one often studies Ripley's $K$-function (see e.g.\ \cite{25,27} and the references therein). For $\beta=0$, $L_n(\beta)$ is - up to a rescaling - an estimator of Ripley's $K$-function at $r_n$. While one usually considers Ripley's $K$-function for a range of arguments, we choose here only one value $r_n$ that depends on the sample size $n$. An extensive theory of properties and the asymptotic behaviour of $L_n(\beta)$ in the complete spatial randomness setting can be found in \cite{00}.

Based on the asymptotically standardised statistics $L_n(\beta)$, we propose the test statistics
\begin{equation*}
T_{e,n}(\beta):=\left(\frac{L_n(\beta)- \frac{1}{2}n(n-1)\int_{W^2} {\bf 1}\{\|x-y\|\leq r_n\} \, \|x-y\|^\beta\, \dint (x,y)}{\sqrt{\frac{d\kappa_d}{2(2\beta+d)} }n r_n^{\beta+d/2}}\right)^2
\end{equation*}
and
\begin{equation*}
T_{a,n}(\beta):=\left(\frac{L_n(\beta)-\frac{d\kappa_d}{2(\beta+d)} \, n(n-1) r_n^{\beta+d}}{\sqrt{\frac{d\kappa_d}{2(2\beta+d)} }n r_n^{\beta+d/2}}\right)^2,
\end{equation*}
where $\beta>-d/2$ and rejection of $H_0$ will be for large values of $T_{j,n}(\beta)$, $j\in\{e,a\}$. The indices $e$ and $a$ are abbreviations for 'exact' and 'asymptotic', and they point out that $T_{e,n}(\beta)$ involves $\BE L_n(\beta)$, which can be difficult to compute depending on the shape of the observation window $W$, while $T_{a,n}(\beta)$ uses a simple asymptotic approximation of $\BE L_n(\beta)$, see Theorem \ref{thm:ELn}.

In order to derive distributional limit theorems for $L_n(\beta)$, $T_{e,n}(\beta)$ and $T_{a,n}(\beta)$, we apply a central limit theorem from \cite{19} for triangular schemes of $U$-statistics. For $\beta=0$ the statistic $L_n(\beta)$ was considered as application in \cite{19}. Here, we generalise these findings to $\beta\in (-d/2,\infty)$, which is technical for $\beta\in(-d/2,0)$, and present them in more detail. Moreover, the focus of the present paper is on statistical tests based on $L_n(\beta)$ and their properties, which even for $\beta=0$ goes clearly beyond what was studied in \cite{19}. In \cite{Yang2017} some $U$-statistics based on interpoint distances are proposed as test statistics for uniformity on the unit cube (beside two other statistics based on data depth and normal quantiles). In contrast to $L_n(\beta)$, these $U$-statistics take all interpoint distances into account and not only the small ones, whence their kernels do not depend on $n$ (i.e., the summand associated with two given points from the sample is the same for all $n\in\N$). The tests for multivariate uniformity studied in \cite{BartoszynskiEtAl,SR:2005} are also based on $U$-statistics with fixed kernels, which are more involved to compute than the distances between the sample points. For $U$-statistics with fixed kernels as considered in \cite{BartoszynskiEtAl,SR:2005,Yang2017}, the asymptotic behaviour is much easier to analyse than for $L_n(\beta)$, where the kernels depend on the parameters $n$ and $r_n$ and their interplay.

Simulations on the $d$-dimensional unit cube indicate that the power of $T_{j,n}(\beta)$, $j\in\{e,a\}$, against alternatives depends on the parameters $\beta$ and $r_n$. We show for several parameters that they are serious competitors to time-honoured tests and demonstrate the applicability of the new procedures by analysing the real dataset \texttt{finpines}. Clearly, we leave open questions for further research, as e.g. to find an optimal (automatic) selection of the parameters.

This paper is organised as follows. In Section \ref{sec:Asymptotics} we derive the theory for $L_n(\beta)$ in a general setting, including formulae for the mean and the variance as well as central limit theorems. The two families of test statistics $T_{j,n}(\beta)$, $j\in\{e,a\}$, are discussed in Section \ref{TfU}, and their limiting behaviour is given under $H_0$ and under fixed alternatives. The behaviour for some contiguous alternatives is studied in Section \ref{sec:ContiguousAlternatives}. Section \ref{Simu} provides a simulation study and a comparison to existing methods. We finish the paper by applying the new tests to a real data set in Section \ref{sec:RealData}, and with comments on open problems and research perspectives in Section \ref{sec:Conclusion}.

\section{Properties of $L_n(\beta)$}\label{sec:Asymptotics}
Let ${\mathscr X}_n:=\{X_1,\ldots,X_n\}$, where $n\ge2$ and $X_1,\ldots,X_n$ are i.i.d.\ random vectors distributed according to a density $f$, whose support is contained in a measurable set $W\subset\R^d$ of positive finite volume. In the following, we assume without loss of generality that $\Vol(W)=1$, i.e., $W$ has volume one. For some of our results we need the additional assumption that
\begin{equation}\label{eqn:AdditionalAssumptionW}
\limsup_{r\to 0} \frac{\Vol(\{x\in W: d(x,\partial W)\leq r\})}{r}<\infty.
\end{equation}
Here, we use the notation $d(x,A):=\inf_{y\in A} \|x-y\|$ for $x\in\R^d$ and $A\subset\R^d$. The assumption \eqref{eqn:AdditionalAssumptionW} requires that the volume of the set of points in $W$ that are in the $r$-neighbourhood of the boundary of $W$ is at most of order $r$ and seems to be no significant restriction. For many sets $W$, for example all compact and convex $W$, the limit superior in \eqref{eqn:AdditionalAssumptionW} equals the surface area of $W$. The expression in \eqref{eqn:AdditionalAssumptionW} is related to the so-called (outer) Minkowski content. For a definition as well as some results on its finiteness we refer to \cite{AmbrosioColesantiVilla2008}.

Let $(r_n)$ be a sequence of positive real numbers such that $r_n\to0$ as $n\to\infty$. In the following $B^d(x,r)$ stands for the $d$-dimensional closed ball with centre $x\in\R^d$ and radius $r>0$, and $\kappa_d:=\pi^{d/2}/ \Gamma(d/2 + 1)$ is the volume of the $d$-dimensional unit ball $B^d(0,1)$.
For $p>0$ we denote by $L^p(W)$ the space of all measurable functions on $W$ for which the Lebesgue integral of the $p$-th power of the absolute value is finite.
For the special case $\beta=0$ the formulae of the following theorem can also be found in \cite[Equations (4.2) and (4.3)]{19}.

\begin{theorem}\label{thm:ELn}
For $\beta>-d$ and all densities $f\in L^2(W)$,
\begin{equation}\label{eqn:ELnexact}
\BE L_n(\beta) = \frac{n(n-1)}{2} \int_{W^2} {\bf 1}\{\|x-y\|\leq r_n\} \, \|x-y\|^{\beta} \, f(x) \, f(y) \, \dint (x,y)
\end{equation}
and
\begin{equation}\label{eqn:ELnasymptotic}
\lim_{n\to\infty} \frac{\BE L_n(\beta)}{n^2 r_n^{\beta+d}}=\frac{d\kappa_d}{2(\beta+d)} \int_W f(x)^2 \, \dint x.
\end{equation}
\end{theorem}
Theorem \ref{thm:ELn}, which we prove in Appendix \ref{app:proof:2.1}, states exact formulae for the mean and easy to compute asymptotic approximations under fairly general assumptions. The behaviour of $\BE L_n(\beta)$ under $H_0$ in the next corollary is a direct consequence. We write $g\equiv h$ to indicate that two functions $g,h: W\to\R$ are identical almost everywhere.

\begin{corollary}\label{cor:Mean}
If $\beta>-d$ and $f\equiv{\bf 1}_W$, then
$$
\BE L_n(\beta) = \frac{n(n-1)}{2} \int_{W^2} {\bf 1}\{\|x-y\|\leq r_n\} \, \|x-y\|^{\beta} \, \dint (x,y)
$$
and
$$
\lim_{n\to\infty} \frac{\BE L_n(\beta)}{n^2 r_n^{\beta+d}}=\frac{d\kappa_d}{2(\beta+d)}.
$$
\end{corollary}
Recall that the degree of a vertex in a graph is the number of edges emanating from it. The average degree $\bar{D}_n$ of the vertices in $\operatorname{RGG}(\mathscr{X}_n,r_n)$ is given by $\bar{D}_n = 2 L_n(0)/n$. Thus, it follows from Theorem \ref{thm:ELn} that $\BE \bar{D}_n$ is of the same order as $n r_n^d$ as $n\to\infty$. For the special choice of uniformity $f\equiv{\bf 1}_W$ Corollary \ref{cor:Mean} implies
\begin{equation}\label{eqn:ExpectedAverageDegree}
\lim_{n\to\infty} \frac{\BE \bar{D}_n}{\kappa_d n r_n^d}=1.
\end{equation}

In the next theorem we present exact and asymptotic formulae for the variance of $L_n(\beta)$, which generalise the findings from \cite[Section 4]{19} for $\beta=0$. The proof of the theorem is provided in Appendix \ref{app:proof:2.3}.

\begin{theorem}\label{thm:Variance}
Let $f\in L^3(W)$ and $\beta>-d/2$.
\begin{itemize}
\item [(a)] Then,
\begin{equation}\label{eqn:VarianceLnexact}
\begin{split}
\Var L_n(\beta) 
& = \frac{n(n-1)}{2} \int_{W^2} {\bf 1}\{\|x-y\|\leq r_n\} \, \|x-y\|^{2\beta} \, f(x) \, f(y) \, \dint (x,y)\\
& \quad + n(n-1)(n-2) \int_{W} \bigg( \int_W {\bf 1}\{\|x-y\|\leq r_n\} \, \|x-y\|^\beta \, f(y) \, \dint y\bigg)^2 \, f(x) \, \dint x\\
& \quad - n(n-1)(n-3/2) \bigg( \int_{W^2} {\bf 1}\{\|x-y\|\leq r_n\} \, \|x-y\|^{\beta} \, f(x) \, f(y) \, \dint (x,y)\bigg)^2.
\end{split}
\end{equation}
\item [(b)] For $f\not\equiv \mathbf{1}_W$,
\begin{equation}\label{eqn:VarianceLnasymptotic}
\lim_{n\to\infty}\frac{\Var L_n(\beta)}{\sigma^{(1)}_{\beta,f} n^2 r_n^{2\beta+d}+\sigma^{(2)}_{\beta,f} n^3 r_n^{2\beta+2d}}=1,
\end{equation}
where
$$
\sigma^{(1)}_{\beta,f} := \frac{d\kappa_d}{2(2\beta+d)} \int_W f(x)^2 \, \dint x
$$
and
$$
\sigma^{(2)}_{\beta,f} := \frac{d^2\kappa_d^2}{(\beta+d)^2} \bigg( \int_W f(x)^3 \, \dint x - \bigg(\int_W f(x)^2 \, \dint x \bigg)^2 \bigg).
$$
\item [(c)] If $f\equiv \mathbf{1}_W$, $W$ satisfies \eqref{eqn:AdditionalAssumptionW} and $nr_n^{d+1}\to0$ as $n\to\infty$, then
\begin{equation} \label{eqn:VarianceUniform}
\lim_{n\to\infty}\frac{\Var L_n(\beta)}{n^2 r_n^{2\beta+d}}=\frac{d\kappa_d}{2(2\beta+d)}.
\end{equation}
\end{itemize}
\end{theorem}
Notice that the orders of the two terms in the denominator in \eqref{eqn:VarianceLnasymptotic} differ by $n r_n^d$, which is the order of the expected average degree. For $\sigma^{(1)}_{\beta,f},\sigma^{(2)}_{\beta,f}>0$ this means that the first (second) term dominates if $\BE \bar{D}_n\to 0$ ($\BE \bar{D}_n\to \infty$) as $n\to\infty$, while both terms contribute to the limit if $\BE\bar{D}_n\to c\in(0,\infty)$ as $n\to\infty$. For all densities $f\in L^3(W)$ we have $\sigma^{(1)}_{\beta,f}>0$. The Cauchy-Schwarz inequality implies
$$
 \bigg(\int_W f(x)^2 \, \dint x \bigg)^2 \leq \int_W f(x)^3 \, \dint x \ \int_W f(x) \, \dint x = \int_W f(x)^3 \, \dint x
$$
with equality if and only if $f\equiv {\bf 1}_W$. So $\sigma^{(2)}_{\beta,f}\geq0$ with equality if and only if $f\equiv {\bf 1}_W$.

The formula \eqref{eqn:VarianceUniform} coincides with \eqref{eqn:VarianceLnasymptotic} for $f\equiv \mathbf{1}_W$. Nevertheless we have to impose for \eqref{eqn:VarianceUniform} additional conditions on the boundary of $W$ and on the sequence $(r_n)$. They ensure that the sum of the second and the third term in \eqref{eqn:VarianceLnexact} does not have an asymptotic order that is less than $n^3 r_n^{2\beta+2d}$ but still larger than $n^2r_n^{2\beta+d}$. The following example shows that this can happen due to boundary effects (see also \cite[Section 4]{19}). For $W=[0,1]$, $f\equiv \mathbf{1}_W$, $\beta=0$ and $r_n<1/2$, we have
\begin{align*}
\int_0^1 \bigg( \int_0^1 \mathbf{1}\{|x-y|\leq r_n\} \, \dint y \bigg)^2 \, \dint x & = 2 \int_0^{r_n} (r_n+x)^2 \, \dint x + (1-2r_n) 4r_n^2 \\
& = \frac{2}{3} (8r_n^3 - r_n^3) + (1-2r_n) 4r_n^2 = 4r_n^2-\frac{10}{3} r_n^3
\end{align*}
and
\begin{align*}
\int_{[0,1]^2} {\bf 1}\{|x-y|\leq r_n\} \, \dint(x,y) & = 2\int_0^{r_n} r_n+x \, \dint x + (1-2r_n) 2 r_n\\
& = 4r_n^2 - r_n^2 + (1-2r_n) 2 r_n = 2r_n - r_n^2.
\end{align*}
Thus, the sum of the second and the third term in \eqref{eqn:VarianceLnexact} equals
$$
n(n-1)(n-2) \left(4r_n^2-\frac{10}{3} r_n^3 - 4r_n^2 + 4r_n^3 - r_n^4\right) - \frac{1}{2}n(n-1) (2r_n - r_n^2)^2.
$$
If $nr_n^2\to\infty$ as $n\to\infty$, this is of a higher order than the first term in \eqref{eqn:VarianceLnexact}.

Theorem 3.3 in \cite{00} states asymptotic variances for the same statistics $L_n(\beta)$ with an underlying homogeneous Poisson point process of intensity $n$ (i.e., $f\equiv\mathbf{1}_W$ and the number of points is Poisson-distributed with mean $n$). In contrast to \eqref{eqn:VarianceUniform}, these formulae show the same phase transition depending on the behaviour of $n r_n^d$ as we have in \eqref{eqn:VarianceLnasymptotic} for $f\not\equiv {\bf 1}_W$.

In the following we use the abbreviation
$$
\sigma_{\beta,f,n} := \sqrt{\sigma^{(1)}_{\beta,f} n^2 r_n^{2\beta+d}+\sigma^{(2)}_{\beta,f} n^3 r_n^{2\beta+2d}},
$$
with $\sigma^{(1)}_{\beta,f},\sigma^{(2)}_{\beta,f}$ as in Theorem \ref{thm:Variance} for $\beta>-d/2$ and $n\in\N$. Moreover, we write $\overset{\mathcal{D}}{\longrightarrow}$ for convergence in distribution and $N_m(\mu,\Sigma)$ for an $m$-dimensional Gaussian random vector with mean vector $\mu\in\R^m$ and positive semidefinite covariance matrix $\Sigma\in\R^{m\times m}$. In the univariate case the index $m$ is omitted.

\begin{theorem}\label{thm:CLT}
Let $f\in L^3(W)$, $\beta>-d/2$
and assume that $n^2 r_n^d\to\infty$ as $n\to\infty$. If $f\not\equiv \mathbf{1}_W$ or if $f\equiv\mathbf{1}_W$, $W$ satisfies \eqref{eqn:AdditionalAssumptionW} and $nr_n^{d+1}\to0$ as $n\to\infty$, then
$$
\frac{L_n(\beta)-\BE L_n(\beta)}{\sigma_{\beta,f,n}} \overset{\mathcal{D}}{\longrightarrow} N(0,1) \quad \text{as} \quad n\to\infty.
$$
\end{theorem}
The proof of Theorem \ref{thm:CLT} is provided in Appendix \ref{app:proof:2.4}. For $\beta=0$ a central limit theorem as Theorem \ref{thm:CLT} is established in \cite[Section 4]{19}; see also \cite{Weber1983} and the references therein. In \cite[Section 3.5]{09} central limit theorems for subgraph counts of random geometric graphs are derived, which include the number of edges $L_n(0)$ as special case.
Notice that $n^2r_n^d\to\infty$ as $n\to\infty$ means that the expected number of edges goes to infinity as $n\to\infty$ (see Theorem \ref{thm:ELn}), which is a reasonable assumption for a central limit theorem involving edge lengths. The additional assumptions for $f\equiv \mathbf{1}_W$ are the same as in Theorem \ref{thm:Variance}(c) and are used to ensure that the rescaled variances converge to one.

The following corollary concerning the behaviour under the null hypothesis is proven in Appendix \ref{app:proof:2.5}.

\begin{corollary}\label{cor:CLTuniform}
Let $\beta>-d/2$, $f\equiv{\bf 1}_W$ and assume that $W$ satisfies \eqref{eqn:AdditionalAssumptionW}.
\begin{itemize}[noitemsep,topsep=0pt,parsep=0pt,partopsep=0pt]
\item [(a)] If $n^2 r_n^d\to\infty$ and $nr_n^{d+1}\to 0$ as $n\to\infty$, then
$$
\hspace{-0.5cm} \frac{L_n(\beta)-\frac{n(n-1)}{2} \int_{W^2} {\bf 1}\{\|x-y\|\leq r_n\} \, \|x-y\|^\beta\, \dint (x,y)}{\sqrt{\frac{d\kappa_d}{2(2\beta+d)} }n r_n^{\beta+d/2}} \overset{\mathcal{D}}{\longrightarrow} N(0,1) \quad \text{as} \quad n\to\infty.
$$
\item [(b)] If $n^2 r_n^d\to\infty$ and $n^2 r_n^{d+2}\to0$ as $n\to\infty$, then
$$
\frac{L_n(\beta)-\frac{d\kappa_d}{2(\beta+d)} \, n(n-1) r_n^{\beta+d}}{\sqrt{\frac{d\kappa_d}{2(2\beta+d)} }n r_n^{\beta+d/2}} \overset{\mathcal{D}}{\longrightarrow} N(0,1) \quad \text{as} \quad n\to\infty.
$$
\end{itemize}
\end{corollary}
It can be seen from Corollary \ref{cor:Mean} that in part (a) of the previous corollary $L_n(\beta)$ is centred with its expectation, while in (b) the asymptotic expectation is used. In the latter situation, the assumptions on $(r_n)$ are stricter. For the statistics $L_n(\beta)$ with respect to an underlying homogeneous Poisson point process (i.e.\ the case of complete spatial randomness) central limit theorems are shown in \cite[Section 5.1]{00}.

\section{Testing for uniformity}\label{TfU}
Motivated by Corollary \ref{cor:CLTuniform} we propose testing goodness-of-fit of $H_0$ in (\ref{H0}) against general alternatives based on the families of statistics
\begin{equation}\label{Te}
T_{e,n}(\beta)=\left(\frac{L_n(\beta)- \frac{1}{2}n(n-1)\int_{W^2} {\bf 1}\{\|x-y\|\leq r_n\} \, \|x-y\|^\beta\, \dint (x,y)}{\sqrt{\frac{d\kappa_d}{2(2\beta+d)} }n r_n^{\beta+d/2}}\right)^2
\end{equation}
and
\begin{equation}\label{Ta}
T_{a,n}(\beta)=\left(\frac{L_n(\beta)-\frac{d\kappa_d}{2(\beta+d)} \, n(n-1) r_n^{\beta+d}}{\sqrt{\frac{d\kappa_d}{2(2\beta+d)} }n r_n^{\beta+d/2}}\right)^2,
\end{equation}
depending on $\beta>-\frac d2$ and $r_n\in(0,\infty)$. The choice of the sequence $(r_n)$ is discussed in Section \ref{Simu}, where we introduce a parameter $k$, see \eqref{eq:rn.1} and \eqref{eq:rn.2}. Rejection of $H_0$ will be for large values of $T_{j,n}(\beta),\,j\in\{a,e\}$. Empirical critical values for $W=[0,1]^d$ can be found in Tables \ref{tab:EC2} to \ref{tab:EC3} for dimensions $d=2,3$ and sample sizes $n\in\{50,100,200,500\}$. Notice that under $H_0$ and some mild assumptions on $(r_n)$ and $W$ the continuous mapping theorem and Corollary \ref{cor:CLTuniform} yield
\begin{equation*}\label{TeaH0}
T_{j,n}(\beta) \overset{\mathcal{D}}{\longrightarrow} \chi^2_1 \quad \text{as} \quad n\to\infty,\quad j\in\{a,e\}, \beta>-d/2.
\end{equation*}
Here $\chi^2_1$ denotes a random variable having a chi-squared distribution with one degree of freedom. In the following theorem we consider the asymptotic behaviour of $T_{e,n}(\beta)$ and $T_{a,n}(\beta)$ under fixed alternatives. We write $\overset{\PP}{\longrightarrow}$ for convergence in probability and prove the next theorem in Appendix \ref{app:proof:3.1}.

\begin{theorem}\label{thm:ConsistencyExpectation}
Let $\beta>-d/2$ and $f\not\equiv{\bf 1}_W$. If $n^2 r_n^d\to\infty$ as $n\to\infty$, then
$$
T_{e,n}(\beta) \overset{\PP}{\longrightarrow} \infty \quad \text{and} \quad T_{a,n}(\beta) \overset{\PP}{\longrightarrow}\infty \quad \text{as} \quad n\to\infty.
$$
\end{theorem}
Theorem \ref{thm:ConsistencyExpectation} yields consistency of $T_{e,n}(\beta)$ and $T_{a,n}(\beta)$ against each fixed alternative $f\not\equiv{\bf 1}_W$.

\section{Behaviour under contiguous alternatives}\label{sec:ContiguousAlternatives}

Let $g\in L^3(W)$ be such that $g\not\equiv 0$ and $\int_W g(x) \, \dint x=0$ and let $(a_n)$ be a positive sequence such that $a_n\to 0$ as $n\to\infty$. In the following we always tacitly assume that $1+a_n g(x)\geq 0$ for all $x\in W$ and $n\in\N$. This guarantees that $\mathbf{1}_W+a_n g$ is a density. In the sequel we denote by $\widetilde{T}_{e,n}(\beta)$ and $\widetilde{T}_{a,n}(\beta)$ our test statistics in \eqref{Te} and \eqref{Ta} computed on $n$ i.i.d.\ points $\widetilde{X}_1,\ldots,\widetilde{X}_n$ distributed according to the density $\mathbf{1}_W+a_ng$ (i.e., we have a triangular scheme).

\begin{theorem} \label{thm:CloseAlternatives}
Let $\beta>-d/2$ and assume that $W$ satisfies \eqref{eqn:AdditionalAssumptionW}, that $n^2 r_n^d\to\infty$, $n r_n^{d+1}\to0$ and $\min\{n r_n^{d/2+1}a_n,r_n/a_n \} \to 0$ as $n\to\infty$ and that, for $r>0$,
\begin{equation} \label{eqn:AssumptionWh}
\int_{W} {\bf 1}\{d(x,\partial W)\leq r\} |g(x)| \, \dint x  \leq C_{W,g} r
\end{equation}
with some constant $C_{W,g}\in (0,\infty)$. Then the following assertions hold:
\begin{itemize}
\item [(a)] If $n r_n^{d/2}a_n^2\to \gamma\in[0,\infty)$ as $n\to\infty$, then
$$
\widetilde{T}_{e,n}(\beta) \overset{\mathcal{D}}{\longrightarrow} \left( Z + \frac{\sqrt{d\kappa_d(2\beta+d)}}{\sqrt{2}(\beta+d)} \int_W g(x)^2 \, \dint x \ \gamma\right)^2 \quad \text{as} \quad n\to\infty
$$
with $Z \sim N(0,1)$.

\item [(b)] If $n r_n^{d/2}a_n^2\to \infty$ as $n\to\infty$, then
$$
\widetilde{T}_{e,n}(\beta) \overset{\mathbb{P}}{\longrightarrow}\infty \quad \text{as} \quad n\to\infty.
$$
\item[(c)] If, additionally, $n^2r_n^{d+2}\to 0$ as $n\to\infty$, the statements of (a) and (b) also hold for $\widetilde{T}_{a,n}(\beta)$.
\end{itemize}
\end{theorem}

The condition \eqref{eqn:AssumptionWh} requires that the fluctuations of $g$ in an $r$-neighbourhood of the boundary of $W$ are at most of order $r$. Because we assume \eqref{eqn:AdditionalAssumptionW}, this is always the case if $g$ is bounded.
The limiting random variable in Theorem \ref{thm:CloseAlternatives}(a) follows a non-central chi-squared distribution with one degree of freedom. For $n r_n^{d/2}a_n^2\to 0$ as $n\to\infty$ Theorem \ref{thm:CloseAlternatives} implies that $\widetilde{T}_{e,n}(\beta)$ and $\widetilde{T}_{a,n}(\beta)$ behave exactly as $T_{e,n}(\beta)$ and $T_{a,n}(\beta)$ under $H_0$. As the following result shows one can slightly modify Theorem \ref{thm:CloseAlternatives} if $g$ vanishes close to the boundary of $W$. By $\operatorname{supp} g$, we denote the support of $g$, i.e., the set of all $x\in W$ such that $g(x)\neq 0$. For $A,B\subset\R^d$ let $d(A,B):=\inf_{x\in A,y\in B}\|x-y\|$.

\begin{theorem} \label{thm:CloseAlternativesBoundary}
Let $\beta>-d/2$ and assume that $d(\operatorname{supp}g, \partial W)>0$, that $W$ satisfies \eqref{eqn:AdditionalAssumptionW} and that $n^2 r_n^d\to\infty$ and $nr_n^{d+1}\to0$ as $n\to\infty$. Then, (a), (b) and (c) of Theorem \ref{thm:CloseAlternatives} hold.
\end{theorem}

Theorem \ref{thm:CloseAlternatives} and Theorem \ref{thm:CloseAlternativesBoundary} are proven in Appendix \ref{app:proof:4.1u4.2}. Following these theorems, we conclude that under the stated assumptions the tests based on $\widetilde{T}_{a,n}(\beta)$ and $\widetilde{T}_{e,n}(\beta)$ are able to detect alternatives which converge to the uniform distribution at rate $a_n$. Moreover, the theorems could be the foundation of establishing local optimality of the tests by applying the third Le Cam lemma, see Section 5.2 of \cite{28} for a short review of the needed methodology.

\section{Simulation}\label{Simu}
In this section we compare the finite-sample power performance of the test statistics $T_{e,n}(\beta)$ and $T_{a,n}(\beta)$, $\beta>-d/2$, $n\in\N$, with that of some competitors. Since the $d$-dimensional hypercube $[0,1]^d$ is the mostly used observation window, we restrict our simulation study to this case with $d \in \{2,3\}$. Particular interest will be given to the influence on the finite-sample power of $\beta$ and $r_n$ in dependence of the chosen alternatives. In each scenario, we consider the sample sizes $n\in\{50,100,200,500\}$ and set the nominal level of significance to $0.05$. Since the test statistics depend on the parameter $\beta$ and the choice of $r_n$ and the empirical finite sample quantile is in some cases far away from the quantile $\chi^2_{1,0.95}\approx3.8415$ of the limiting distribution, we simulated critical values for $T_{e,n}(\beta)$ and $T_{a,n}(\beta)$ with $100~000$ replications, see Tables \ref{tab:EC2} to \ref{tab:EC3}. Each stated empirical power of the tests in Tables \ref{tab:Comp.d.2} to \ref{tab:Ta.d.3} is based on $10~000$ replications and the asterisk $*$ denotes a rejection rate of $100\%$.

Since there is a vast variety of ways to choose the parameters $\beta$ and $r_n$, we chose the parameter configurations to fit the limiting regimes of Corollary \ref{cor:CLTuniform} as well as the following additional property: From \eqref{eqn:ExpectedAverageDegree} we know that the expectation of the average degree $\bar{D}_n$ behaves as $\kappa_d n r_n^d$ for $n\to\infty$ under $H_0$. This observation motivates the following choices of the radius $r_n$ for $T_{e,n}(\beta)$, namely
\begin{equation}\label{eq:rn.1}
r_n = \left(\frac{k}{n \kappa_d}\right)^{\frac{1}{d}}, \quad k \in \{1, \dotsc, 10, 15, 20, 25\},
\end{equation}
which satisfies $n^2 r_n^d \to \infty$ and $n r_n^{d+1} \to 0$ as $n\to\infty$ and ensures $\BE \bar{D}_n\to k$ as $n\to\infty$ under $H_0$. For the test statistic $T_{a,n}(\beta)$ the additional condition $n^2 r_n^{d+2}\to0$ as $n\to\infty$ has to be fulfilled, so we choose
\begin{equation}\label{eq:rn.2}
r_n = \left( \frac{k}{n^{\frac{3}{2}}\kappa_d} \right)^{\frac{1}{d}}, \quad k \in \{1, \dotsc, 10, 15, 20, 25\},
\end{equation}
to guarantee this additional assumption for $d\in\{2,3\}$. In this case we have $\BE \bar{D}_n\to 0$ as $n\to\infty$, which for $d=2$ is always the case if $n^2 r_n^{d+2}\to0$ as $n\to\infty$.

The expected value $\BE L_n(\beta)$ depends on the observation window $W$ as well as on the dimension $d\ge2$. The following lemma provides exact formulae of $\BE L_n(\beta)$ for each of the cases simulated and is proved in Appendix \ref{app:proof:5.1}.

\begin{lemma}\label{lem:expect}
Assume $\beta>-d$ and $f\equiv{\bf 1}_W$.
\begin{enumerate}
\item [(a)] If $d=2$, $W=[0,1]^2$ and $r_n\leq 1$, then
\begin{align*}
\BE L_n(\beta) = \frac{n(n-1)}{2}\left( \frac{2\pi}{\beta+2} r_n^{\beta+2}-\frac{8}{\beta+3} r_n^{\beta+3}+\frac{2}{\beta+4} r_n^{\beta+4} \right).
\end{align*}
\item [(b)] If $d=3$, $W=[0,1]^3$ and $r_n\leq 1$, then
\begin{align*}
\BE L_n(\beta) = \frac{n(n-1)}{2}\left( \frac{4\pi}{\beta+3} r_n^{\beta+3}-\frac{6\pi}{\beta+4} r_n^{\beta+4}+\frac{8}{\beta+5} r_n^{\beta+5}-\frac{1}{\beta+6} r_n^{\beta+6}\right).
\end{align*}
\end{enumerate}
\end{lemma}

As competitors to the new test statistics we consider the distance to boundary test ($DB$-test), see \cite{05}, the maximal spacing test ($MS$-test), see \cite{03,henze_2018}, the nearest neighbour type test ($NN$-test) of \cite{24} as well as the Bickel-Rosenblatt test ($BR$-test) presented in \cite{20}. We follow the descriptions of the $DB$- and $MS$-tests given in \cite{24}.

For the $NN$-test we consider the family of statistics
\begin{equation*} \label{EHY}
NN_{n,J}^{(\beta)}:= \sum_{x \in \mathscr{X}_n} \xi_{n,J}^{(\beta)}(x, \mathscr{X}_n)
\end{equation*}
in dependence of $\beta\in (0, \infty)$, where $J$ is the number of nearest neighbours, with $x^{(k)}$ being the $k$-nearest neighbour of $x\in\mathscr{X}_n$ and
\begin{equation*}
\xi_{n,J}^{(\beta)}(x, \mathscr{X}_n):= \sum_{k = 1}^J (\kappa_d \|n^{1/d}(x - x^{(k)})\|^d)^{\beta}.
\end{equation*}
To avoid boundary problems in the computation of the $NN$-test, we used the same toroid metric in the simulation as in \cite{24}. Since rejection rates depend crucially on the
power $\beta$ and the number of neighbours $J$ taken into account, we chose different values for $\beta$ and $J$ for the two alternatives where the choice was motivated by Table 2 in \cite{24}. Notice that this test is consistent, but one has to be careful to choose the correct rejection region, which depends on the choice of $\beta$.


As a further competitor we consider the fixed bandwidth  $BR$-test on the unit cube, studied in \cite{20}. The corresponding test statistic is
$${BR}_n^2(h) = -I_n^{2,1}(h) + I_n^{2,2}(h) + V_h(0) + n (V_h \star \overline{U} \star U)(0),$$
with
$$I_n^{2,1}(h) = 2 \sum_{i=1}^n (V_h \star U)(X_i) \quad \text{and} \quad I_n^{2,2}(h) = \frac{2}{n} \sum_{1\le i< j \le n} V_h(X_i - X_j),$$
where $h>0$ is a fixed bandwidth. For the sake of completeness we restate the following abbreviations, see \cite{20}. The convolution product operator is denoted by $\star$,
$U=\mathbf{1}_{[0,1]^d}$ is the density of the uniform distribution over the unit hypercube $[0,1]^d$ and for any function $g$ we define $\overline{g}(x) := g(-x)$ and $g_h(x) := g\left(\frac{x}{h}\right)/h^d$ with $h>0$. Furthermore, we set $V := \overline{K}\star K$, where $K$ is a product kernel on $\R^d$,
that is, $K(u) = \prod_{i=1}^d k(u_i)$, $u=(u_1,\dotsc,u_d)\in\R^d$ with a kernel $k$ on $\R$ (so $k$ is bounded and integrable).
Using the arguments and techniques in \cite{20}, direct calculations for $d=2$ and $k(x)=\frac{1}{\sqrt{2\pi}}\exp\left(-\frac{x^2}{2}\right)$, $x\in\R$, being the standard Gaussian density function, give for $h>0$,
\begin{align*}
I_n^{2,1}(h) &= 2 \sum_{i=1}^n \left( \Phi\left( \frac{X_{i,1}-1}{\sqrt{2}h} \right) - \Phi\left( \frac{X_{i,1}}{\sqrt{2}h} \right) \right)
\left( \Phi\left( \frac{X_{i,2}-1}{\sqrt{2}h} \right) - \Phi\left( \frac{X_{i,2}}{\sqrt{2}h} \right) \right),\\
I_n^{2,2}(h) &= \frac{1}{2\pi n h^4}\sum_{1\le i < j \le n}\exp\left(-\frac{(X_{i,1}-X_{j,1})^2}{4h^2}\right)\exp\left(-\frac{(X_{i,2}-X_{j,2})^2}{4h^2}\right)\\
\end{align*}
and
\begin{equation*}
V_h(0) = \frac{1}{4\pi h^4}, \quad n(V_h \star \overline{U} \star U)(0) = \frac{4n}{\pi h^2} \left[ \sqrt{\pi}\left( \Phi\left( \frac{1}{\sqrt{2}h} \right) - \frac{1}{2} \right) + h \left( \exp\left( -\frac{1}{4h^2} \right) -1 \right) \right]^2,
\end{equation*}
where $\Phi$ is the standard Gaussian distribution function and $X_{i,j}$ denotes the $j$-th component of the random vector $X_i$, with $i\in\{1, \dotsc, n\}$ and $j\in\{1,2\}$.

The $BR$-test rejects the null hypothesis for large values of ${BR}_n^2(h)$. Notice that the asymptotic distribution of ${BR}_n^2(h)$ is known, see \cite{20}, but not in a closed form. Hence we simulated critical values of ${BR}_n^2(h)$ for $h \in \{0.1, 0.25, 0.5, 1\}$, which can be found in Table \ref{tab:crit.val.BR}.
\begin{table}[t]
\centering
\begin{tabular}{c|cccc}
$n \backslash h$ & $0.1$ & $0.25$ & $0.5$ & $1$ \\\hline
50 & 9113.028 & 827.3781 & 72.70593 & 0.01799183\\
100 & 17048.245 & 1616.5611 & 144.16370 & 0.01799641\\
200 & 32839.801 & 3186.1990 & 286.73272 & 0.01795072\\
500 & 80073.212 & 7876.7591 & 713.74621 & 0.01787482
\end{tabular}
\caption{Critical values of the $BR$-statistic ${BR}_n^2(h)$}\label{tab:crit.val.BR}
\end{table}

Following the studies in \cite{06,24}, we simulated a contamination and a clustering model as alternatives to the uniform distribution.
In addition, we considered an alternative consisting of a single point source within uniformly distributed points.
The contamination alternative (CON) is given by the mixture
\begin{align*}
(1-q_1-q_2) {\cal U} ([0,1]^d) + q_1 N_d (c_1, \sigma_1^2 I_d) + q_2 N_d (c_2,\sigma_2^2 I_d),
\end{align*}
under the condition that all simulated points are located in $[0,1]^d$. Here, $I_d \in \R^{d\times d}$ denotes the identity matrix of order $d$. The chosen parameters are given in Table \ref{tab:pc.con}, where $\Phi^{-1}(p)$, $p \in (0,1)$, denotes the $p$-quantile of a standard Gaussian distribution. See Figure \ref{fig:9plots}, second row, for a realisation of this model, where the normally distributed contamination points are filled points and filled squares, respectively.

\begin{table}[b]
\renewcommand{\arraystretch}{.7}
\centering
\begin{tabular}{c||c|c|c|c|c|c}
$d$ & $q_1$ & $q_2$ & $c_1$ & $c_2$ & $\sigma_1$ & $\sigma_2$ \\
\hline
\hline
2 & 0.135 & 0.24 & $(0.25, 0.25)$ & $(0.7, 0.7)$ & $0.15 \cdot \Phi^{-1}(\sqrt{0.9})$ & $0.2 \cdot \Phi^{-1}(\sqrt{0.9})$ \\
\hline
3 & 0.135 & 0.24 & $(0.25, 0.25, 0.25)$ & $(0.7, 0.7, 0.7)$ & $0.15 \cdot \Phi^{-1}(\sqrt[3]{0.9})$ & $0.2 \cdot \Phi^{-1}(\sqrt[3]{0.9})$ \\
\end{tabular}
\caption{Parameter configuration of the CON-alternatives}\label{tab:pc.con}
\end{table}
The clustering alternative (CLU) is motivated by a fixed number of data points version of a Mat\'ern cluster process, see Section 12.3 in \cite{spatstat}, and is designed to destroy the independence. One first chooses a radius $r_{\text{clu}}$ and simulates $\frac{n}{5}$ random points with the uniform distribution ${\cal U}([-r_{\text{clu}},1+r_{\text{clu}}]^d)$, that act as centres of clusters. These points will not be part of the final sample. In a second step, one generates $5$ points around each centre in a ball with radius $r_{\text{clu}}$. These points are generated independently of each other and follow uniform distributions on the mentioned balls. If a point falls outside $[0,1]^d$, it is replaced by a point that follows a ${\cal U}([0,1]^d)$ distribution. In the following we set $r_{\text{clu}}=0.1$ and a realisation of this model can be found in Figure \ref{fig:9plots}, third row. The clustering alternative is not included in the framework of our theoretical results since the points are, by construction, not independent. Nevertheless it is interesting to see how the test statistics behave for such alternatives, which were also considered in the simulation study in \cite{05}.

For the single point source alternative (SPS), we simulate a large number of uniformly distributed points and disturb them with a few points from a single source. In detail, on average $95\%$ of the points are uniformly distributed on $[0,1]^d$. The remaining $5\%$ of the points are derived from a $N_d(c,\sigma^2 I_d)$ distribution under the condition that all simulated points are located in $[0,1]^d$. Here, the parameters are given by $c=(0.5,\ldots,0.5)\in\R^d$ and $\sigma=0.01$. This alternative is designed to emphasise the dependency of the statistics on the parameter $\beta$.

We now present the simulation results for $d=2$. Table \ref{tab:Comp.d.2} exhibits the empirical percentage of rejection of the competing procedures under discussion. An asterisk stands for power of $100\%$, and in each row the best performing procedures have been highlighted using boldface ciphers. Clearly, $BR_n^2(0.1)$ and $NN_{n,15}^{(0.5)}$ dominate the other procedures for the CON-alternative, but as noted in \cite{24} the performance of $NN_{n,J}^{(0.5)}$ might even increase for bigger values of $J$. Comparison with $T_{e,n}(\beta)$ for $\beta=-0.5$ (see Table \ref{tab:Te.d.2}) shows that the presented new methods are for sample sizes of $n=100,200,500$ as good as and for $n=50$ nearly as good as the best competitor $BR_n^2(0.1)$. As one can witness throughout the Tables \ref{tab:Te.d.2} and \ref{tab:Ta.d.2}, $T_{e,n}(\beta)$ dominates $T_{a,n}(\beta)$ for small sample sizes, while the power is similar to the best competitors. In case of the CLU alternative $T_{e,n}(\beta)$ gives the overall highest performance for $\beta=-0.5$ over small sample sizes of $n=50,100,200$, while the only procedure that is better for $n=500$ is again $NN_{n,15}^{(0.5)}$. Notice that the asymptotic version $T_{a,n}(\beta)$ might even achieve higher performance if one considers bigger radii, since it attains the highest rates for the biggest values of $k$. A closer look at these tables reveals the dependency of the new tests on the choice of $\beta$ and $k$. Interestingly, the highest performance is given for both alternatives and $T_{j,n}(\beta)$, $j\in\{a,e\}$, for the choice of $\beta=-0.5$. The best choice of $k$ obviously depends on the sample size.
The dependency of the test statistics on the parameter $\beta$ becomes even clearer in the Tables \ref{tab:SPS.Te} and \ref{tab:SPS.Ta}, which contain the empirical rejection rates under the SPS alternative. Here, the best choice is obviously $\beta=-0.5$. One explanation for this behaviour could be that for $\beta=-0.5$ very small distances between the data points are taken more into account. Under the SPS alternative, some of the data points actually are very close to each other. Thus in case $\beta=-0.5$ the presented test statistics seem to be particularly suitable to detect a single point source between uniformly distributed points.

Observe that the simulation results for $d=3$ in Tables \ref{tab:Te.d.3} and \ref{tab:Ta.d.3} show higher rejection rates for $T_{j,n}(\beta)$ than in the bivariate setting. Since the other methods were too time consuming to implement or to simulate we restrict the comparison to the $DB$-test. As can be seen in Table \ref{tab:Te.d.3} the new tests dominate the $DB$-method for $\beta=-0.5$ and nearly for every value of $k$.

\section{Real data example: Finnish Pines}\label{sec:RealData}
\begin{figure}[t]
\centering
\includegraphics[scale=0.65]{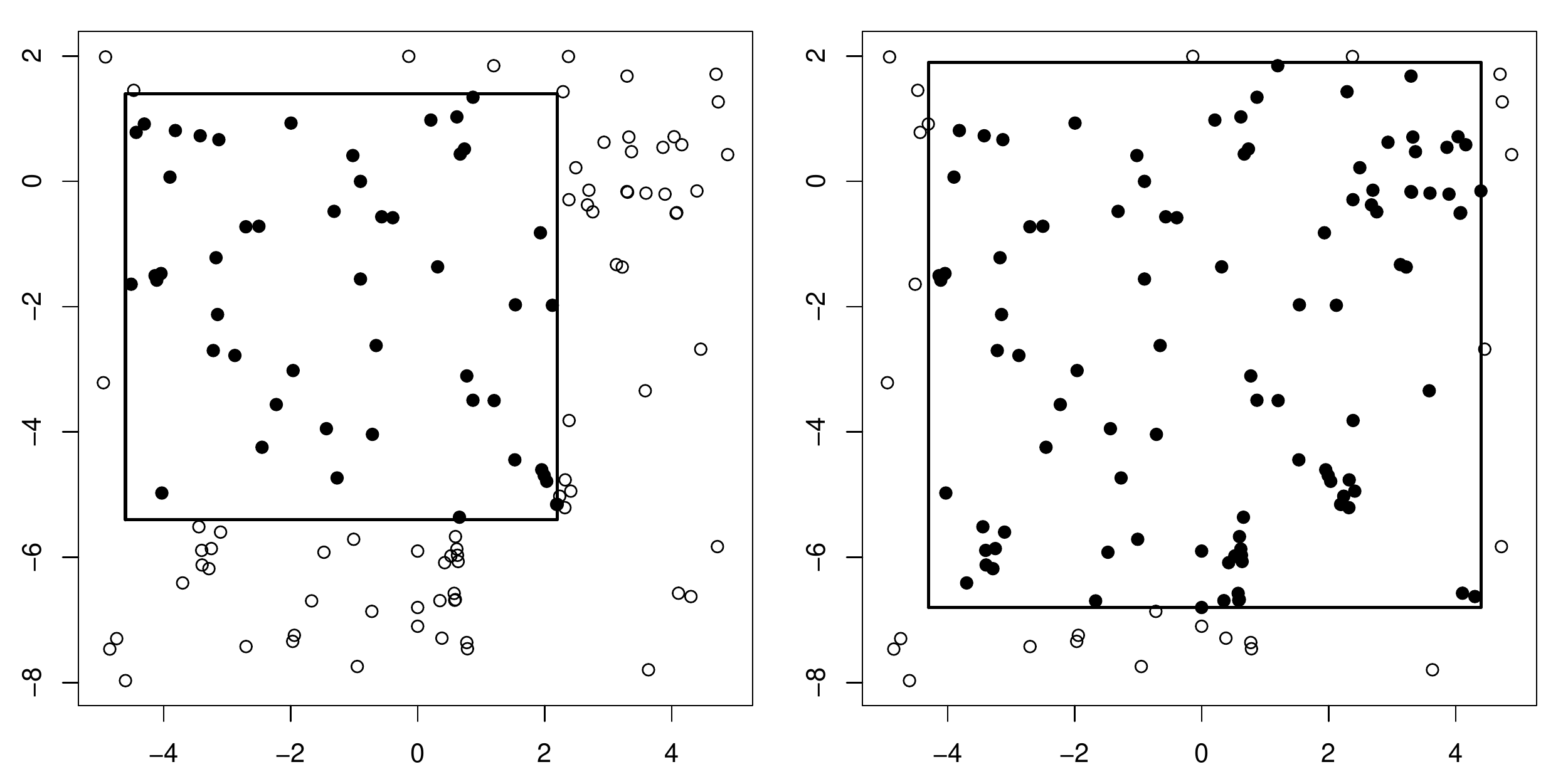}
\caption{Selection of $n=50$ (left) and $n=100$ (right) data points in the data set \texttt{finpines} }\label{dat:finpines}
\end{figure}
We apply our methods to the data set \texttt{finpines} included in the \texttt{R} package \texttt{spatstat}, see \cite{BT:2005}, which presents the locations of 126 pine saplings in a Finnish forest, the locations are given in metres (to six significant digits). In order to compute small sample sizes, i.e. $n=50$ and $n=100$, we restricted the data set to two specific observation windows $W$, see Figure \ref{dat:finpines}.
We test the hypothesis $H_0$ in (\ref{H0}), i.e. if the points are uniformly distributed in $W$, and apply the new methods as well as all presented tests from Section \ref{Simu}. Results are reported as empirical $p$-values (based on 10~000 replications) for all procedures and are found for the tests $T_{e,n}$ and $T_{a,n}$ in Table \ref{tab:finpines2} as well as for the competing tests in Table \ref{tab:finpines1}. Interestingly, in the first example ($n=50$) the distance to boundary test, the maximal spacings test, the nearest neighbour tests and the Bickel-Rosenblatt tests uniformly fail to reject the hypothesis of uniformity on a $5\%$ level, whereas $T_a$ rejects the hypothesis for $k\le10$ and $\beta=-0.5$ as well as for $k=1$ and $\beta\in\{0,1\}$. In the second example ($n=100$), again most of the competitors fail to reject $H_0$ on a $5\%$ level, with exception of the nearest neighbour tests, which show an empirical $p$-value of 0. Impressively, $T_{e,n}$ as well as $T_{a,n}$ reject $H_0$ for $\beta=-0.5$ for every $k\in\{1,\ldots,10,15,20,25\}$, showing for the negative exponent the overall best power. The nearest neighbour tests also show a similar behaviour, which is not very surprising due to the related concepts of both procedures.
\begin{table}[t]
\renewcommand{\arraystretch}{.99}
\centering
\small
\begin{tabular}{c|c|c|ccccccccccccc}
&$n$ & $\beta\backslash k$ & 1 & 2 & 3 & 4 & 5 & 6 & 7 & 8 & 9 & 10 & 15 & 20 & 25 \\
\hline
\multirow{8}{*}{$T_{e,n}$}&\multirow{4}{*}{$50$}
&$-0.5$&0.9&2.5&12.3&34.0&59.5&89.3&60.4&72.5&78.5&69.3&51.8&50.5&38.5\\
&&$0$&13.9&33.2&95.4&61.6&35.6&8.9&4.5&8.9&14.0&11.8&11.6&16.4&12.5\\
&&$1$&46.7&78.8&52.5&20.1&10.1&1.7&0.9&3.6&8.0&7.3&8.9&18.4&13.0\\
&&$5$&67.7&95.1&46.6&15.3&9.2&0.8&0.7&16.3&39.0&29.0&14.5&55.0&20.6\\
\cline{2-16}
&\multirow{4}{*}{$100$}
&$-0.5$&0.0&0.0&0.0&0.0&0.0&0.0&0.0&0.0&0.0&0.0&0.1&0.4&1.2\\
&&$0$&0.0&0.0&0.0&0.0&0.0&0.0&0.1&0.8&4.5&9.9&28.1&68.5&87.2\\
&&$1$&0.0&0.0&0.1&0.1&1.4&8.6&22.0&91.3&49.5&35.9&51.4&30.7&31.9\\
&&$5$&5.5&1.1&59.3&9.7&70.7&73.4&47.4&4.5&0.9&1.4&59.8&25.3&35.3\\
\hline
\hline
\multirow{8}{*}{$T_{a,n}$}&\multirow{4}{*}{$50$}
&$-0.5$&0.1&0.5&0.9&0.4&0.5&2.6&1.2&5.5&1.4&2.8&22.5&70.3&86.8\\
&&$0$&0.3&11.3&26.1&5.2&10.9&52.1&27.7&72.8&28.0&45.5&100&50.5&53.2\\
&&$1$&1.1&92.9&91.8&18.9&38.5&79.2&69.5&63.4&57.0&84.2&56.6&23.5&22.5\\
&&$5$&6.9&17.6&64.0&1.0&31.7&37.6&82.7&24.2&43.0&95.1&33.0&14.0&23.2\\
\cline{2-16}
&\multirow{4}{*}{$100$}
&$-0.5$&0.0&0.0&0.0&0.0&0.0&0.0&0.0&0.0&0.0&0.0&0.0&0.0&0.0\\
&&$0$&0.0&0.0&0.0&0.0&0.0&0.0&0.0&0.0&0.0&0.0&0.0&0.0&0.0\\
&&$1$&0.0&0.0&0.0&0.0&0.0&0.0&0.0&0.0&0.0&0.0&0.0&0.0&0.1\\
&&$5$&0.0&0.1&3.9&0.0&0.1&0.5&0.0&0.0&0.9&13.0&5.2&9.7&40.6\\
\hline
\end{tabular}
\caption{Empirical $p$-values for $T_{e}$ and $T_{a}$ for the subsets of size $n$ of the \texttt{finpines} data set}\label{tab:finpines2}
\end{table}

\begin{table}[t]
\centering
\begin{tabular}{c|cccc|cccc}
$n$ & $BR_n^2(0.1)$ & $BR_n^2(0.25)$ & $BR_n^2(0.5)$ & $BR_n^2(1)$ & $NN_{n,1}^{(0.5)}$ & $NN_{n,15}^{(0.5)}$ & DB & MS \\
\hline
50  & 88.8 & 95.4 & 97.1 & 66.0  & 5.7 & 45.7 & 29.9 & 67.2\\
\hline
100 & 8.6 & 93.2 & 95.1 & 34.6 & 0.1 & 0.0 & 31.2 & 65.3 \\
\hline
\end{tabular}
\caption{Empirical $p$-values for the competing tests for the subsets of size $n$ of the \texttt{finpines} data set}\label{tab:finpines1}
\end{table}

\section{Conclusions and open problems}\label{sec:Conclusion}
We have theoretically investigated statistics related to the edge length of the random geometric graph of a point pattern in an observation window under fairly general assumptions. From these findings, we introduced two new families of consistent goodness-of-fit tests of uniformity based on random geometric graphs. As the simulation section shows, the presented methods are serious competitors to existing methods, even dominating them for right choices of the parameters $\beta$ and $r_n$ (or $k$). Clearly, a natural question is to find (data dependent) best choices of them. Another obvious extension of the presented methods would be to find tests of uniformity on (lower dimensional) manifolds, including special cases of directional statistics as the circle or the sphere (for existing methods see Chapter 6 of either \cite{28} or \cite{29}). Section \ref{sec:ContiguousAlternatives} invites to further investigate in view of concepts of locally optimal tests. Since the approach is fairly general, an extension would be testing the fit of $X_1,\ldots,X_n$ to some parametric family $\{f(\cdot,\vartheta):\vartheta\in\Theta\}$ of densities for a specific parameter space $\Theta$ (eventually the procedures would use a suitable estimator $\widehat{\vartheta}_n$ of $\vartheta$). In view of the special interest in the case of unknown support of the data, see \cite{03,06}, we indicate that the definition of $T_{a,n}(\beta)$ is not dependent on the shape of the underlying observation window and therefore is applicable in this setting (as long as the observation window has volume one).

\section*{Acknowledgements}
The authors thank Norbert Henze and Steffen Winter for fruitful discussions.

\bibliography{literatur}

\begin{thebibliography}{10}

\bibitem{AmbrosioColesantiVilla2008}
L.~Ambrosio, A.~Colesanti, and E.~Villa.
\newblock Outer {M}inkowski content for some classes of closed sets.
\newblock {\em Math. Ann.}, 342(4):727--748, 2008.

\bibitem{spatstat}
A.~Baddeley, E.~Rubak, and R.~Turner.
\newblock {\em Spatial Point Patterns: Methodology and Applications with {R}}.
\newblock Chapman and Hall/CRC Press, 2015.

\bibitem{BT:2005}
A.~Baddeley and R.~Turner.
\newblock {spatstat}: An {R} package for analyzing spatial point patterns.
\newblock {\em Journal of Statistical Software}, 12(6):1--42, 2005.

\bibitem{Baringhaus2018}
L.~Baringhaus, D.~Gaigall, and J.~P. Thiele.
\newblock Statistical inference for ${L}^2$-distances to uniformity.
\newblock {\em Comput. Stat.}, 33(4):1863--1896, 2018.

\bibitem{BartoszynskiEtAl}
R.~Bartoszy\'{n}ski, D.~K. Pearl, and J.~Lawrence.
\newblock A multidimensional goodness-of-fit test based on interpoint
  distances.
\newblock {\em J. Amer. Statist. Assoc.}, 92(438):577--586, 1997.

\bibitem{03}
J.~R. Berrendero, A.~Cuevas, and B.~Pateiro-L\'{o}pez.
\newblock A multivariate uniformity test for the case of unknown support.
\newblock {\em Stat. Comput.}, 22(1):259--271, 2012.

\bibitem{06}
J.~R. Berrendero, A.~Cuevas, and B.~Pateiro-L\'{o}pez.
\newblock Testing uniformity for the case of a planar unknown support.
\newblock {\em Canad. J. Statist.}, 40(2):378--395, 2012.

\bibitem{05}
J.~R. Berrendero, A.~Cuevas, and F.~V\'{a}zquez-Grande.
\newblock Testing multivariate uniformity: The distance-to-boundary method.
\newblock {\em Canad. J. Statist.}, 34(4):693--707, 2006.

\bibitem{BiauDevroye}
G.~Biau and L.~Devroye.
\newblock {\em Lectures on the Nearest Neighbor Method}.
\newblock Springer, 2015.

\bibitem{Billingsley1979}
P.~Billingsley.
\newblock {\em Probability and Measure}.
\newblock Wiley, 1979.

\bibitem{25}
N.~Cressie.
\newblock {\em Statistics for Spatial Data}.
\newblock Wiley, 1993.

\bibitem{26}
B.~Ebner, N.~Henze, M.~A. Klatt, and K.~Mecke.
\newblock Goodness-of-fit tests for complete spatial randomness based on
  {M}inkowski functionals of binary images.
\newblock {\em Electron. J. Stat.}, 12(2):2873--2904, 2018.

\bibitem{24}
B.~Ebner, N.~Henze, and J.~E. Yukich.
\newblock Multivariate goodness-of-fit on flat and curved spaces via nearest
  neighbor distances.
\newblock {\em J. Multivar. Anal.}, 165:231--242, 2018.

\bibitem{Gilbert1961}
E.~N. Gilbert.
\newblock Random plane networks.
\newblock {\em J. Soc. Ind. Appl. Math.}, 9:533--543, 1961.

\bibitem{Godehardt1997}
E.~Godehardt, J.~Jaworski, and D.~Godehardt.
\newblock The application of random coincidence graphs for testing the
  homogeneity of data.
\newblock In {\em Classification, Data Analysis, and Data Highways (Potsdam,
  1997)}, pages 35--45. Springer, 1998.

\bibitem{27}
L.~Heinrich.
\newblock {G}aussian limits of empirical multiparameter ${K}$-functions of
  homogeneous {P}oisson processes and tests for complete spatial randomness.
\newblock {\em Lith. Math. J.}, 55(1):72--90, 2015.

\bibitem{henze_2018}
N.~Henze.
\newblock On the consistency of the spacings test for multivariate uniformity,
  including on manifolds.
\newblock {\em J. Appl. Probab.}, 55(2):659--665, 2018.

\bibitem{19}
S.~R. Jammalamadaka and S.~Janson.
\newblock Limit theorems for a triangular scheme of {U}-statistics with
  applications to inter-point distances.
\newblock {\em Ann. Probab.}, 14:1347--1358, 1986.

\bibitem{02}
A.~Justel, D.~Pena, and R.~Zamar.
\newblock A multivariate {K}olmogorov-{S}mirnov test of goodness of fit.
\newblock {\em Statist. Probab. Lett.}, 35:251--259, 1997.

\bibitem{knuth}
D.~E. Knuth.
\newblock {\em The art of computer programming}, volume 2: Seminumerical
  algorithms.
\newblock Addison-Wesley, 3rd. edition, 2011.

\bibitem{21}
E.~O. Krauczi.
\newblock Joint cluster counts from uniform distribution.
\newblock {\em Probab. Math. Statist.}, 33(1):93--106, 2013.

\bibitem{28}
C.~Ley and T.~Verdebout.
\newblock {\em Modern Directional Statistics}.
\newblock Chapman and Hall/CRC Press, New York, 2017.

\bibitem{04}
J.~Liang, K.~Fang, F.~Hickernell, and R.~Li.
\newblock Testing multivariate uniformity and it's applications.
\newblock {\em Math. Comp.}, 70(233):337--355, 2001.

\bibitem{29}
K.~V. Mardia and P.~E. Jupp.
\newblock {\em Directional Statistics}.
\newblock Wiley, 2000.

\bibitem{12}
Y.~Marhuenda, D.~Morales, and M.~C. Pardo.
\newblock A comparison of uniformity tests.
\newblock {\em Statistics}, 39(4):315--328, 2005.

\bibitem{09}
M.~Penrose.
\newblock {\em Random geometric graphs}.
\newblock Oxford University Press, 2003.

\bibitem{08}
A.~Petrie and T.~R. Willemain.
\newblock An empirical study of tests for uniformity in multidimensional data.
\newblock {\em Comput. Statist. Data Anal.}, 64:253--268, 2013.

\bibitem{00}
M.~Reitzner, M.~Schulte, and C.~Th{\"a}le.
\newblock Limit theory for the {G}ilbert graph.
\newblock {\em Adv. Appl. Math.}, 88:26--61, 2017.

\bibitem{22}
M.~Rosenblatt.
\newblock Remarks on a multivariate transformation.
\newblock {\em Ann. Math. Stat.}, 23(3):470--472, 1952.

\bibitem{Rudin1974}
W.~Rudin.
\newblock {\em Real and complex analysis}.
\newblock McGraw-Hill, 2nd edition, 1974.

\bibitem{SR:2005}
G.~J. {Sz\'ekely} and M.~L. {Rizzo}.
\newblock A new test for multivariate normality.
\newblock {\em J. Multivariate Anal.}, 93(1):58--80, 2005.

\bibitem{20}
C.~Tenreiro.
\newblock On the finite sample behavior of fixed bandwidth
  {B}ickel-{R}osenblatt test for univariate and multivariate uniformity.
\newblock {\em Comm. Statist. Simulation Comput.}, 36:827--846, 2007.

\bibitem{Weber1983}
N.~C. Weber.
\newblock Central limit theorems for a class of symmetric statistics.
\newblock {\em Math. Proc. Cambridge Philos. Soc.}, 94(2):307--313, 1983.

\bibitem{Yang2017}
M.~Yang and R.~Modarres.
\newblock Multivariate tests of uniformity.
\newblock {\em Stat. Papers}, 58(3):627--639, 2017.

\end{thebibliography}

\begin{appendices}

\section{Proofs}\label{app:proofs}

\subsection{Proof of Thoerem \ref{thm:ELn}}\label{app:proof:2.1}
\noindent {\it Proof:} Equation \eqref{eqn:ELnexact} follows from
$$
\BE L_n(\beta) = \frac{n(n-1)}{2} \BE {\bf 1}\{\|X-Y\|\leq r_n\} \|X-Y\|^\beta,
$$
where $X$ and $Y$ are independent random vectors distributed according to the density $f$. Notice that
\begin{align*}
\BE {\bf 1}\{\|X-Y\|\leq r_n\} \|X-Y\|^\beta & = \int_{W^2} {\bf 1}\{\|x-y\|\leq r_n\} \, \|x-y\|^{\beta} \, f(x) f(y) \, \dint (x,y) \\
& \leq \int_{W^2} {\bf 1}\{\|x-y\|\leq r_n\} \, \|x-y\|^{\beta} \, f(x)^2 \, \dint (x,y)\\
& \leq \frac{d\kappa_d}{\beta+d} r_n^{\beta+d}\int_W f(x)^2 \, \dint x,
\end{align*}
where we used the inequality of arithmetic and geometric means and spherical coordinates. This yields
\begin{equation}\label{eqn:limsupE}
 \limsup_{n\to\infty} \frac{\BE L_n(\beta)}{n^2 r_n^{\beta+d}} \leq \frac{d\kappa_d}{2(\beta+d)} \int_W f(x)^2 \, \dint x.
\end{equation}
For $C>0$ we use the shorthand notation $f_C(x):=\min\{f(x),C\}$ for $x\in W$ and $f_C(x):=0$ for $x\notin W$. It follows from Lemma \ref{lem:Differentiation} that, for any $C>0$,
$$
\lim_{n\to\infty} \frac{1}{r_n^{\beta+d}}\int_{B^d(x,r_n)} \|x-y\|^\beta f_C(y) \, \dint y = \frac{d\kappa_d}{\beta+d} f_C(x)
$$
for almost all $x\in W$. Now the dominated convergence theorem yields
$$
\lim_{n\to\infty} \frac{1}{r_n^{\beta+d}}\int_{W^2} {\bf 1}\{\|x-y\|\leq r_n\} \, \|x-y\|^{\beta} \, f_C(x) f_C(y) \, \dint (x,y)
 = \frac{d\kappa_d}{\beta+d}\int_{W} f_C(x)^2\, \dint x.
$$
Together with
\begin{align*}
&  \int_{W^2} {\bf 1}\{\|x-y\|\leq r_n\} \, \|x-y\|^{\beta} \, f(x) f(y) \, \dint (x,y) \\
&  \geq \int_{W^2} {\bf 1}\{\|x-y\|\leq r_n\} \, \|x-y\|^{\beta} \, f_C(x) f_C(y) \, \dint (x,y)
\end{align*}
we obtain
$$
\liminf_{n\to\infty} \frac{\BE L_n(\beta)}{n^2 r_n^{\beta+d}} \geq \frac{d\kappa_d}{2(\beta+d)} \int_W f_C(x)^2 \, \dint x.
$$
Now letting $C\to\infty$ and the monotone convergence theorem yield
$$
\liminf_{n\to\infty} \frac{\BE L_n(\beta)}{n^2 r_n^{\beta+d}} \geq \frac{d\kappa_d}{2(\beta+d)} \int_W f(x)^2 \, \dint x.
$$
Combining this with \eqref{eqn:limsupE} proves \eqref{eqn:ELnasymptotic}.\hfill$\square$

\subsection{Proof of Theorem \ref{thm:Variance}}\label{app:proof:2.3}
\noindent {\it Proof:} A straightforward computation shows that
\begin{align*}
\BE L_n(\beta)^2 
 & = \frac{n(n-1)}{2} \BE {\bf 1}\{\|X_1-X_2\|\leq r_n\} \|X_1-X_2\|^{2\beta} \\
& \quad + n(n-1)(n-2) \BE {\bf 1}\{\|X_1-X_2\|,\|X_1-X_3\|\leq r_n\} \|X_1-X_2\|^\beta \|X_1-X_3\|^\beta \\
& \quad + \frac{n(n-1)(n-2)(n-3)}{4} \BE {\bf 1}\{\|X_1-X_2\|,\|X_3-X_4\|\leq r_n\} \|X_1-X_2\|^\beta \|X_3-X_4\|^\beta.
\end{align*}
Here, $X_1,\ldots,X_4$ are independent random vectors with density $f$. Combining this with \eqref{eqn:ELnexact} yields \eqref{eqn:VarianceLnexact}.

Observe that the asymptotic behaviour of the first and the third term in \eqref{eqn:VarianceLnexact} follows immediately from Theorem \ref{thm:ELn}. By the inequality of arithmetic and geometric means and spherical coordinates, we obtain
\begin{equation}\label{eqn:ProofVarUpperBound}
\begin{split}
& \frac{1}{r_n^{2\beta+2d}}\int_{W} \bigg( \int_W {\bf 1}\{\|x-y\|\leq r_n\} \, \|x-y\|^\beta \, f(y) \, \dint y\bigg)^2 \, f(x) \, \dint x \\
& \leq \frac{1}{3 r_n^{2\beta+2d}} \int_{W^3} {\bf 1}\{\|x_1-x_2\|,\|x_1-x_3\|\leq r_n\} \|x_1-x_2\|^\beta \|x_1-x_3\|^\beta \\
& \hskip 2.75cm (f(x_1)^3+f(x_2)^3+f(x_3)^3) \, \dint(x_1,x_2,x_3)\\
& \leq \frac{d^2\kappa_d^2}{(\beta+d)^2} \int_W f(x)^3 \, \dint x.
\end{split}
\end{equation}
On the other hand, Lemma \ref{lem:Differentiation} and the dominated convergence theorem imply
\begin{align*}
 \lim_{n\to\infty} \frac{1}{r_n^{2\beta+2d}}\int_{W} \bigg( \int_W {\bf 1}\{\|x-y\|\leq r_n\} \, \|x-y\|^\beta \, f_C(y) \, \dint y\bigg)^2 \, f_C(x) \, \dint x 
 = \frac{d^2\kappa_d^2}{(\beta+d)^2} \int_W f_C(x)^3 \, \dint x
\end{align*}
for each $C>0$. Recall that $f_C(x)=\min\{f(x),C\}$ for $x\in W$. Now letting $C\to\infty$ and the monotone convergence theorem yield
\begin{align*}
 \liminf_{n\to\infty}\frac{1}{r_n^{2\beta+2d}}\int_{W} \bigg( \int_W {\bf 1}\{\|x-y\|\leq r_n\} \, \|x-y\|^\beta \, f(y) \, \dint y\bigg)^2 \, f(x) \, \dint x 
 \geq \frac{d^2\kappa_d^2}{(\beta+d)^2} \int_W f(x)^3 \, \dint x.
\end{align*}
This, together with \eqref{eqn:ProofVarUpperBound} and the observation that $\sigma^{(1)}_{\beta,f},\sigma^{(2)}_{\beta,f}>0$, completes the proof of \eqref{eqn:VarianceLnasymptotic}.

For the proof of \eqref{eqn:VarianceUniform} we define $W_{-r_n}:=\{x\in W: d(x,\partial W)\geq r_n \}$. Now straightforward computations yield
\begin{align*}
 \frac{d^2\kappa_d^2}{(\beta+d)^2} r_n^{2\beta+2d} \Vol(W_{-r_n}) 
 \leq \int_W \bigg( \int_W  \mathbf{1}\{\|x-y\|\leq r_n\} \|x-y\|^\beta \, \dint y \bigg)^2 \, \dint x \leq \frac{d^2\kappa_d^2}{(\beta+d)^2} r_n^{2\beta+2d} \Vol(W)
\end{align*}
and
\begin{align*}
 \frac{d^2\kappa_d^2}{(\beta+d)^2} r_n^{2\beta+2d} \Vol(W_{-r_n})^2 
 \leq \bigg(\int_{W^2} \mathbf{1}\{\|x-y\|\leq r_n\} \|x-y\|^{\beta}  \, \dint (x,y) \bigg)^2 \leq \frac{d^2\kappa_d^2}{(\beta+d)^2} r_n^{2\beta+2d} \Vol(W)^2.
\end{align*}
It follows from \eqref{eqn:AdditionalAssumptionW} that there exists a constant $C_W\in(0,\infty)$ such that
\begin{equation} \label{eqn:ConsequenceBoundary}
0\leq \Vol(W) - \Vol(W_{-r_n}) \leq  \Vol(\{x\in W: d(x,\partial W)\leq r_n\}) \leq C_W r_n.
\end{equation}
Together with $\Vol(W)=1$ this means that the absolute value of the sum of the second and the third term in \eqref{eqn:VarianceLnexact} can be bounded by
$$
\frac{3d^2\kappa_d^2}{(\beta+d)^2} C_W n^3 r_n^{2\beta+2d+1} + \frac{d^2\kappa_d^2}{2(\beta+d)^2} n^2 r_n^{2\beta+2d}.
$$
Together with the asymptotic order of the first term in \eqref{eqn:VarianceLnexact}, which is as in the proof of \eqref{eqn:VarianceLnasymptotic}, this proves \eqref{eqn:VarianceUniform}.
\hfill$\square$

\subsection{Proof of Theorem \ref{thm:CLT}}\label{app:proof:2.4}
We prepare the proof of Theorem \ref{thm:CLT} by several lemmas, which are formulated for the following more general setting, required later: We assume that the underlying points of $\mathscr{X}_n$ are distributed according to some density $f_n\in L^3(W)$ and that $f_n(x)\to f(x)$ as $n\to\infty$ for almost all $x\in W$.

For $n\in\N$ we define
$$
W_{f_n}:=\{(x,m)\in W\times [0,\infty): m\leq f_n(x)\}
$$
and let $\widehat{X}_1,\ldots,\widehat{X}_n$ be independent and uniformly distributed points in $W_{f_n}$. We denote the collection of these points by $\widehat{\mathscr{X}}_n$. For a point $\hat{x}\in W_{f_n}$ we often use the decomposition $\hat{x}=(x,m)$ with $x\in W$ and $m\in [0,f_n(x)]$. Observe that the first components of $\widehat{X}_1,\ldots,\widehat{X}_n$ are distributed according to the density $f_n$ in $W$. For $\beta\in\R$ we define
\begin{equation}\label{eqn:LnHat}
\widehat{L}_{n}(\beta) := \frac{1}{2} \sum_{((x_1,m_1),(x_2,m_2))\in\widehat{\mathscr{X}}_{n,\neq}^2} {\bf 1}\{\|x_1-x_2\|\leq r_n\} \, \|x_1-x_2\|^\beta.
\end{equation}
If $f_n=f$, $\widehat{L}_{n}(\beta)$ has the same distribution as $L_{n}(\beta)$.
For $M>0$ and $a\geq 0$ we define
\begin{equation}\label{eqn:LnM}
\widehat{L}_{n,M}(\beta) := \frac{1}{2} \sum_{((x_1,m_1),(x_2,m_2))\in\widehat{\mathscr{X}}_{n,\neq}^2} {\bf 1}\{m_1,m_2\leq M\} \, {\bf 1}\{ \|x_1-x_2\|\leq r_n\} \, \|x_1-x_2\|^\beta
\end{equation}
and
\begin{equation*}
\begin{split}
 \widehat{L}_{n,a,M}(\beta) 
 := \frac{1}{2} \sum_{((x_1,m_1),(x_2,m_2))\in\widehat{\mathscr{X}}_{n,\neq}^2} {\bf 1}\{m_1,m_2\leq M\} \, {\bf 1}\{n^{-2/d}a\leq \|x_1-x_2\|\leq r_n\} \, \|x_1-x_2\|^\beta.
\end{split}
\end{equation*}
Moreover, we use the abbreviations $f_{n,M}(x):=\min\{f_n(x),M\}$ and $f_M(x):=\min\{f(x),M\}$ for $x\in W$.

\begin{lemma}\label{lem:VarianceLnaM}
Let $\beta>-d/2$, $M\geq 1$, $a>0$ and assume that $n^2r_n^d\to\infty$ as $n\to\infty$ and that $\lim_{n\to\infty} \Var \widehat{L}_{n,M}(\beta)/\sigma^2_{\beta,f_M,n}=1$. Then,
\begin{equation}\label{eqn:LnaMLnM}
\lim_{n\to\infty} \BE \bigg(\frac{\widehat{L}_{n,M}(\beta) - \BE \widehat{L}_{n,M}(\beta)}{\sigma_{\beta,f_M,n}} - \frac{\widehat{L}_{n,a,M}(\beta) - \BE \widehat{L}_{n,a,M}(\beta)}{\sigma_{\beta,f_M,n}} \bigg)^2=0
\end{equation}
and
\begin{equation}\label{eqn:VarLnaM}
\lim_{n\to\infty} \frac{\Var \widehat{L}_{n,a,M}(\beta)}{\sigma_{\beta,f_M,n}^2} = 1.
\end{equation}
\end{lemma}

\noindent {\it Proof}:
Throughout the proof we assume that $n$ is so large that $n^{-2/d}a<r_n$, which is no restriction since $n^2r_n^d\to\infty$ as $n\to\infty$. By definition, we have that
\begin{align*}
 \widehat{L}_{n,M}(\beta) - \widehat{L}_{n,a,M}(\beta)
 = \frac{1}{2} \sum_{((x_1,m_1),(x_2,m_2))\in\widehat{\mathscr{X}}_{n,\neq}^2} {\bf 1}\{m_1,m_2\leq M\} \, {\bf 1}\{ \|x_1-x_2\| < n^{-2/d}a\} \, \|x_1-x_2\|^\beta.
\end{align*}
Now a similar computation as in the proof of Theorem \ref{thm:Variance}(a) yields that
$$
\Var (\widehat{L}_{n,M}(\beta) - \widehat{L}_{n,a,M}(\beta)) \leq I_1 + I_2
$$
with
\begin{align*}
I_1 & := \frac{n^2}{2} \int_{W^2} \mathbf{1}\{\|x-y\|\leq n^{-2/d}a \} \|x-y\|^{2\beta} f_{n,M}(x) f_{n,M}(y) \, \dint(x,y) \allowdisplaybreaks\\
I_2 & := n^3 \int_W  \bigg( \int_W \mathbf{1}\{\|x-y\|\leq n^{-2/d}a\} \|x-y\|^\beta f_{n,M}(y) \, \dint y \bigg)^2 f_{n,M}(x) \, \dint x.
\end{align*}
Note that $I_1$ and $I_2$ correspond to the first two terms in \eqref{eqn:VarianceLnexact}, whereas the third term in \eqref{eqn:VarianceLnexact} was omitted since it is non-positive. Now short computations show that
$$
\frac{I_1}{\sigma_{\beta,f_M,n}^2} \leq  \frac{d\kappa_dM^2}{2(2\beta+d)}  \frac{n^2 n^{-2-4\beta/d} a^{2\beta+d}}{\sigma^{(1)}_{\beta,f_M}n^2 r_n^{2\beta+d}} = \frac{d\kappa_dM^2a^{2\beta+d}}{2(2\beta+d)\sigma^{(1)}_{\beta,f_M} }  \frac{1}{(n^2r_n^d)^{2\beta/d+1}}
$$
and
$$
\frac{I_2}{\sigma_{\beta,f_M,n}^2} \leq\frac{d^2\kappa_d^2M^3}{(\beta+d)^2}  \frac{n^3 n^{-4-4\beta/d} a^{2\beta+2d}}{\sigma^{(1)}_{\beta,f_M} n^2 r_n^{2\beta+d}}
= \frac{d^2\kappa_d^2M^3a^{2\beta+2d}}{(\beta+d)^2 \sigma^{(1)}_{\beta,f_M} }  \frac{1}{n(n^2r_n^d)^{2\beta/d+1}}.
$$
Since $n^2r_n^d\to\infty$ as $n\to\infty$ and $2\beta/d+1>0$, this provides \eqref{eqn:LnaMLnM}. Now \eqref{eqn:VarLnaM} follows from combining \eqref{eqn:LnaMLnM} and $\lim_{n\to\infty} \Var \widehat{L}_{n,M}(\beta)/\sigma^2_{\beta,f_{M},n}=1$.
\hfill$\square$

\begin{lemma}\label{lem:CLTLnaM}
Let $\beta>-d/2$, $M\geq 1$, $a>0$ and assume that $n^2r_n^d\to\infty$ as $n\to\infty$ and that $\lim_{n\to\infty} \Var \widehat{L}_{n,M}(\beta)/\sigma^2_{\beta,f_M,n}=1$. Then,
\begin{equation}\label{eqn:WeakConvergenceLnaM}
\frac{\widehat{L}_{n,a,M}(\beta)-\BE \widehat{L}_{n,a,M}(\beta)}{\sigma_{\beta,f_M,n}} \overset{\mathcal{D}}{\longrightarrow} N(0,1) \quad \text{as} \quad n\to\infty
\end{equation}
and
\begin{equation}\label{eqn:WeakConvergenceLnM}
\frac{\widehat{L}_{n,M}(\beta)-\BE \widehat{L}_{n,M}(\beta)}{\sigma_{\beta,f_M,n}} \overset{\mathcal{D}}{\longrightarrow} N(0,1) \quad \text{as} \quad n\to\infty.
\end{equation}
\end{lemma}

\noindent {\it Proof:}
From \eqref{eqn:VarLnaM} we know that $\lim_{n\to\infty} \Var \widehat{L}_{n,a,M}(\beta)/\sigma_{\beta,f_M,n}^2=1$. If $\beta\geq 0$, then
\begin{align*}
& \lim_{n\to\infty} \frac{\sup_{(x,m_x),(y,m_y)\in W_{f_n}} \mathbf{1}\{m_x,m_y\leq M\} \mathbf{1}\{n^{-2/d}a\leq \|x-y\|\leq r_n\} \|x-y\|^{\beta} }{n r_n^{\beta+d/2}}\\
& \leq\lim_{n\to\infty} \frac{r_n^\beta}{n r_n^{\beta+d/2}}=\lim_{n\to\infty} \frac{1}{\sqrt{n^2 r_n^d}}=0,
\end{align*}
while for $\beta\in(-d/2,0)$,
\begin{align*}
& \lim_{n\to\infty} \frac{\sup_{(x,m_x),(y,m_y)\in W_{f_n}} \mathbf{1}\{m_x,m_y\leq M\} \mathbf{1}\{n^{-2/d}a\leq \|x-y\|\leq r_n\} \|x-y\|^{\beta} }{n r_n^{\beta+d/2}}\\
& \leq\lim_{n\to\infty} \frac{n^{-2\beta/d}a^{\beta}}{n r_n^{\beta+d/2}}=\lim_{n\to\infty} \frac{a^\beta}{(n^2 r_n^d)^{1/2+\beta/d}}=0.
\end{align*}
Denoting by $(X_1^{(n)},m_{X_1^{(n)}})$ an uniformly distributed point in $W_{f_n}$, we obtain
\begin{align*}
& \lim_{n\to\infty} \frac{n\sup_{(x,m)\in W_{f_n}}\BE \mathbf{1}\{m,m_{X_1^{(n)}}\leq M\} \mathbf{1}\{n^{-2/d}a \leq \|x-X_1^{(n)}\|\leq r_n\} \|x-X_1^{(n)}\|^{\beta}}{n r_n^{\beta+d/2}}\\
& = \lim_{n\to\infty} \frac{1}{r_n^{\beta+d/2}} \sup_{x\in W}\int_W {\bf 1}\{n^{-2/d}a\leq\|x-y\|\leq r_n\} \, \|x-y\|^\beta \, f_{n,M}(y)\, \dint y \\
& \leq	 \lim_{n\to\infty} \frac{d\kappa_d M}{\beta+d} \frac{r_n^{\beta+d}}{r_n^{\beta+d/2}}=\lim_{n\to\infty} \frac{d\kappa_d M}{\beta+d} r_n^{d/2} =0.
\end{align*}
Thus, \eqref{eqn:WeakConvergenceLnaM} follows from Theorem \ref{thm:CLTUstatistics}. Combining the $L^2$-covergence in \eqref{eqn:LnaMLnM} with \eqref{eqn:WeakConvergenceLnaM} yields \eqref{eqn:WeakConvergenceLnM}.\hfill$\square$

In the following we use the abbreviation $\overline{f}_{n,M}(x):=\max\{f_n(x)-M,0\}$ for $x\in W$ and $M\geq 0$.

\begin{lemma}\label{lem:ApproximationM}
For $n\in\N$, $\beta>-d/2$ and $M\geq 1$,
\begin{align*}
& \BE\bigg(\frac{\widehat{L}_{n}(\beta)-\BE \widehat{L}_{n}(\beta)}{\sigma_{\beta,f,n}} - \frac{\widehat{L}_{n,M}(\beta)-\BE \widehat{L}_{n,M}(\beta)}{\sigma_{\beta,f,n}}  \bigg)^2\\
& \leq \frac{d\kappa_d}{2\beta+d} \frac{n^2r_n^{2\beta+d}}{\sigma_{\beta,f,n}^2} \int_{W}  \overline{f}_{n,M}(x)^2+M \overline{f}_{n,M}(x) \, \dint x\\
& \quad + \frac{18d^2\kappa_d^2}{(\beta+d)^2} \frac{n^3r_n^{2\beta+2d}}{\sigma_{\beta,f,n}^2} \int_W M^2 \overline{f}_{n,M}(x) + M \overline{f}_{n,M}(x)^2 + \overline{f}_{n,M}(x)^3 \, \dint x.
\end{align*}
\end{lemma}

\noindent {\it Proof:}
By definition we have
\begin{align*}
\widehat{L}_{n}(\beta)-\widehat{L}_{n,M}(\beta) 
 = \frac{1}{2}\sum_{((x,m_x),(y,m_y))\in\widehat{\mathscr{X}}_{n,\neq}^2} {\bf 1}\{ m_x>M \text{ or } m_y>M\} \, {\bf 1}\{\|x-y\|\leq r_n\} \|x-y\|^\beta.
\end{align*}
From similar arguments as in the proofs of Theorem \ref{thm:Variance}(a) and Lemma \ref{lem:VarianceLnaM}, it follows that
$$
\Var (\widehat{L}_n(\beta)-\widehat{L}_{n,M}(\beta)) \leq I_1+I_2
$$
with
\begin{align*}
I_1 & := \frac{n^2}{2} \int_{W_{f_n}^2} {\bf 1}\{m_1>M \text{ or } m_2>M\} \, {\bf 1}\{\|x_1-x_2\|\leq r_n\} \, \|x_1-x_2\|^{2\beta} 
\, \dint((x_1,m_1),(x_2,m_2)) \allowdisplaybreaks\\
I_2 & := n^3 \int_{W_{f_n}^3} {\bf 1}\{m_1>M \text{ or } m_2>M\} \, {\bf 1}\{m_1>M \text{ or } m_3>M\} \, {\bf 1}\{\|x_1-x_2\|\leq r_n\}\\
& \qquad \qquad \quad \quad \times {\bf 1}\{\|x_1-x_3\|\leq r_n\} \, \|x_1-x_2\|^\beta \, \|x_1-x_3\|^\beta \, \dint((x_1,m_1),(x_2,m_2),(x_3,m_3)).
\end{align*}
For $I_1$ we obtain the bound
\begin{align*}
I_1 & \leq n^2 \int_{W^2} {\bf 1}\{\|x-y\|\leq r_n\} \, \|x-y\|^{2\beta} \, \overline{f}_{n,M}(x) \, (\overline{f}_{n,M}(y)+M) \, \dint(x,y) \\
& \leq n^2 \int_{W^2} {\bf 1}\{\|x-y\|\leq r_n\} \, \|x-y\|^{2\beta} \, (\overline{f}_{n,M}(x)^2+M \overline{f}_{n,M}(x)) \, \dint(x,y) \\
& \leq  \frac{d\kappa_d}{2\beta+d} n^2r_n^{2\beta+d} \int_{W}  \overline{f}_{n,M}(x)^2+M \overline{f}_{n,M}(x) \, \dint x.
\end{align*}
Because of
\begin{align*}
{\bf 1}\{m_1>M \text{ or } m_2>M\} \, {\bf 1}\{m_1>M \text{ or } m_3>M\} 
\leq {\bf 1}\{m_1>M\} +{\bf 1}\{m_2>M,m_3>M\},
\end{align*}
we have
\begin{align*}
I_2 & \leq n^3 \int_{W^3}  {\bf 1}\{\|x_1-x_2\|\leq r_n\}\, {\bf 1}\{\|x_1-x_3\|\leq r_n\} \, \|x_1-x_2\|^\beta \, \|x_1-x_3\|^\beta \\
& \qquad \qquad \quad \quad \times (\overline{f}_{n,M}(x_1) \, f_n(x_2) \, f_n(x_3) + f_n(x_1) \, \overline{f}_{n,M}(x_2) \, \overline{f}_{n,M}(x_3)) \, \dint(x_1,x_2,x_3).
\end{align*}
Using that $f_n(x)\leq \overline{f}_{n,M}(x) + M$ for $x\in\R^d$, we obtain
$$
\overline{f}_{n,M}(x_1) \, f_n(x_2) \, f_n(x_3) + f_n(x_1) \, \overline{f}_{n,M}(x_2) \, \overline{f}_{n,M}(x_3) \leq 6 \max_{k,i\in\{1,2,3\}} M^{3-k} \overline{f}_{n,M}(x_i)^{k}.
$$
This implies
$$
I_2 \leq \frac{18d^2\kappa_d^2}{(\beta+d)^2} n^3r_n^{2\beta+2d} \int_W M^2 \overline{f}_{n,M}(x) + M \overline{f}_{n,M}(x)^2 + \overline{f}_{n,M}(x)^3 \, \dint x,
$$
which completes the proof.\hfill$\square$

We recall that $f_M(x):=\min\{f(x),M\}$ for $x\in W$ and $M\geq 0$.

\begin{lemma}\label{lem:VarianceLnM}
Let $\beta>-d/2$, $M\geq 1$ and $f_n=f$, $n\in\N$. If $f\not\equiv \mathbf{1}_W$ or if $f\equiv\mathbf{1}_W$, $W$ satisfies \eqref{eqn:AdditionalAssumptionW} and $nr_n^{d+1}\to0$ as $n\to\infty$, then
$$
\lim_{n\to\infty}  \frac{\Var \widehat{L}_{n,M}(\beta)}{\sigma_{\beta,f_M,n}^2}=1.
$$
\end{lemma}

\noindent {\it Proof}:
For $M\geq 1$ and $f\equiv \mathbf{1}_W$, $f_M\equiv \mathbf{1}_W$ and the statement is the same as Theorem \ref{thm:Variance}(c) because $\widehat{L}_{n,M}(\beta)$ follows the same distribution as $L_n(\beta)$. For $f\not\equiv\mathbf{1}_W$ one can show as in the proof of Theorem \ref{thm:Variance}(a) that
\begin{align*}
\Var \widehat{L}_{n,M}(\beta) & = \frac{n(n-1)}{2} \int_{W^2} {\bf 1}\{\|x-y\|\leq r_n\} \, \|x-y\|^{2\beta} \, f_M(x) \, f_M(y) \, \dint (x,y)\\
& \quad + n(n-1)(n-2) \int_{W} \bigg( \int_W {\bf 1}\{\|x-y\|\leq r_n\} \, \|x-y\|^\beta \, f_M(y) \, \dint y\bigg)^2 \, f_M(x) \, \dint x\\
& \quad - n(n-1)(n-3/2) \bigg( \int_{W^2} {\bf 1}\{\|x-y\|\leq r_n\} \, \|x-y\|^{\beta} \, f_M(x) \, f_M(y) \, \dint (x,y)\bigg)^2.
\end{align*}
Now the assertion can be proved as Theorem \ref{thm:Variance}(b).
\hfill$\square$

\noindent {\it Proof of Theorem \ref{thm:CLT}:}
We consider the same setting as in the previous lemmas with $f_n=f$ for $n\in\N$ so that $L_n(\beta)$ has the same distribution as $\widehat{L}_n(\beta)$, which we study throughout this proof. For $f\equiv \mathbf{1}_W$ the assertion follows from \eqref{eqn:WeakConvergenceLnM} in Lemma \ref{lem:CLTLnaM} because, for $M\geq1$, $\widehat{L}_n(\beta)$ has the same distribution as $\widehat{L}_{n,M}(\beta)$, $\sigma_{\beta,f_M,n}=\sigma_{\beta,f,n}$ and Lemma \ref{lem:VarianceLnM} guarantees that the variance condition in Lemma \ref{lem:CLTLnaM} is satisfied. So we assume $f\not\equiv\mathbf{1}_W$ in the sequel.

Let $h: \R\to\R$ be a bounded Lipschitz function whose Lipschitz constant is at most one and let $\varepsilon>0$. In the following we show
\begin{equation}\label{eqn:CLTtoprove}
\limsup_{n\to\infty} \left|\BE h\bigg(\frac{\widehat{L}_n(\beta)-\BE \widehat{L}_n(\beta)}{\sigma_{\beta,f,n}}\bigg) - \BE h(N(0,1))\right|\leq \varepsilon,
\end{equation}
which yields the assertion.

For $M\geq 1$ the triangle inequality implies
\begin{equation}\label{eqn:Triangle}
\begin{split}
& \left|\BE h\bigg(\frac{\widehat{L}_n(\beta)-\BE \widehat{L}_n(\beta)}{\sigma_{\beta,f,n}}\bigg) - \BE h(N(0,1))\right|\\
& \leq \left|\BE h\bigg(\frac{\widehat{L}_n(\beta)-\BE \widehat{L}_n(\beta)}{\sigma_{\beta,f,n}}\bigg) - \BE h\bigg(\frac{\widehat{L}_{n,M}(\beta)-\BE \widehat{L}_{n,M}(\beta)}{\sigma_{\beta,f,n}}\bigg)\right|\\
& \quad + \left|\BE h\bigg(\frac{\widehat{L}_{n,M}(\beta)-\BE \widehat{L}_{n,M}(\beta)}{\sigma_{\beta,f,n}}\bigg) - \BE h\bigg(\frac{\widehat{L}_{n,M}(\beta)-\BE \widehat{L}_{n,M}(\beta)}{\sigma_{\beta,f_M,n}}\bigg)\right|\\
& \quad + \left|\BE h\bigg(\frac{\widehat{L}_{n,M}(\beta)-\BE \widehat{L}_{n,M}(\beta)}{\sigma_{\beta,f_M,n}}\bigg) -\BE h(N(0,1))\right|\\
& =: R_{1,n,M}+R_{2,n,M}+R_{3,n,M}.
\end{split}
\end{equation}
It follows from Lemma \ref{lem:CLTLnaM} (notice that the variance condition is satisfied because of Lemma \ref{lem:VarianceLnM}) that $R_{3,n,M}$ vanishes for any $M\geq 1$ as $n\to\infty$. The Lipschitz property of $h$, the Cauchy-Schwarz inequality and Lemma \ref{lem:ApproximationM} imply that
\begin{align*}
R_{1,n,M}^2 & \leq \BE \bigg(\frac{\widehat{L}_{n}(\beta)-\BE \widehat{L}_{n}(\beta)}{\sigma_{\beta,f,n}} - \frac{\widehat{L}_{n,M}(\beta)-\BE \widehat{L}_{n,M}(\beta)}{\sigma_{\beta,f,n}}\bigg)^2\\
& \leq \frac{d\kappa_d}{2\beta+d} \frac{n^2r_n^{2\beta+d}}{\sigma_{\beta,f,n}^2} \int_{W}  \overline{f}_M(x)^2+M \overline{f}_M(x) \, \dint x\\
& \quad + \frac{18d^2\kappa_d^2}{(\beta+d)^2} \frac{n^3r_n^{2\beta+2d}}{\sigma_{\beta,f,n}^2} \int_W M^2 \overline{f}_M(x) + M \overline{f}_M(x)^2 + \overline{f}_M(x)^3 \, \dint x.
\end{align*}
Here the terms depending on $n$ can be bounded by some constants. The dominated convergence theorem with the upper bounds $2f^2$ and $3f^3$ leads to
$$
\lim_{M\to\infty} \int_{W}  \overline{f}_M(x)^2+M \overline{f}_M(x) \, \dint x =0
$$
and
$$
\lim_{M\to\infty} \int_W M^2\overline{f}_M(x) + M \overline{f}_M(x)^2 + \overline{f}_M(x)^3  \, \dint x =0.
$$
Hence, there exists an $M_1\geq 1$ such that $\limsup_{n\to\infty} R_{1,n,M}\leq \varepsilon/2$ for $M>M_1$.

A short computation using the Lipschitz continuity of $h$ and the Cauchy-Schwarz inequality shows that
$$
R_{2,n,M} \leq \bigg| \frac{\sigma_{\beta,f_M,n}}{\sigma_{\beta,f,n}} -1 \bigg| \ \BE \bigg| \frac{\widehat{L}_{n,M}(\beta)-\BE \widehat{L}_{n,M}(\beta)}{\sigma_{\beta,f_M,n}} \bigg| \leq \bigg| \frac{\sigma_{\beta,f_M,n}}{\sigma_{\beta,f,n}} -1 \bigg| \frac{\sqrt{\Var \widehat{L}_{n,M}(\beta)}}{\sigma_{\beta,f_M,n}} .
$$
By the monotone convergence theorem and the assumption $f\not\equiv \mathbf{1}_W$, we have $\sigma^{(1)}_{\beta,f_M}\to \sigma^{(1)}_{\beta,f}>0$ and $\sigma^{(2)}_{\beta,f_M}\to \sigma^{(2)}_{\beta,f}>0$ as $M\to\infty$. Together with the definitions of $\sigma_{\beta,f_M,n}$ and $\sigma_{\beta,f,n}$ this implies that there exists an $M_2\geq 1$ such that
$$
\limsup_{n\to\infty}  \bigg| \frac{\sigma_{\beta,f_M,n}}{\sigma_{\beta,f,n}} -1 \bigg| \leq \frac{\varepsilon}{2}
$$
for $M>M_2$. Since, by Lemma \ref{lem:VarianceLnM}, $\lim_{n\to\infty} \sqrt{\Var \widehat{L}_{n,M}(\beta)}/\sigma_{\beta,f_M,n}=1$, we obtain $\limsup_{n\to\infty} R_{2,n,M} \leq \varepsilon/2$
for $M>M_2$. Thus, choosing $M>\max\{M_1,M_2\}$ in \eqref{eqn:Triangle} and letting $n\to\infty$ yields \eqref{eqn:CLTtoprove} and completes the proof.\hfill$\square$

\subsection{Proof of Corollary \ref{cor:CLTuniform}}\label{app:proof:2.5}
\noindent {\it Proof of Corollary \ref{cor:CLTuniform}:} Part (a) is an immediate consequence of Theorem \ref{thm:CLT}, \eqref{eqn:ELnexact} and the definition of $\sigma_{\beta,f,n}$. For the proof of (b) recall that $W_{-r_n}=\{x\in W: d(x,\partial W)\geq r_n \}$. It follows from \eqref{eqn:ELnexact} that
$$
 \frac{d\kappa_d}{2(\beta+d)} \Vol(W_{-r_n}) n (n-1) r_n^{\beta+d} \leq \BE L_n(\beta) \leq \frac{d\kappa_d}{2(\beta+d)} \Vol(W) n (n-1) r_n^{\beta+d}.
$$
Together with \eqref{eqn:ConsequenceBoundary}, which is valid because we assume \eqref{eqn:AdditionalAssumptionW}, and $\Vol(W)=1$ this yields
$$
|\BE L_n(\beta) - \frac{d\kappa_d}{2(\beta+d)}  n (n-1) r_n^{\beta+d}| \leq \frac{d\kappa_d}{2(\beta+d)} C_W n^2 r_n^{\beta+d+1}
$$
so that
\begin{equation}\label{eqn:DiffExactAsymptotic}
\lim_{n\to\infty}\frac{|\BE L_n(\beta) - \frac{d\kappa_d}{2(\beta+d)}  n (n-1) r_n^{\beta+d}|}{\sqrt{\frac{d\kappa_d}{2(2\beta+d)} }n r_n^{\beta+d/2}}\leq \lim_{n\to\infty} \frac{\sqrt{d\kappa_d(2\beta+d)}}{\sqrt{2}(\beta+d)} C_W n r_n^{d/2+1}=0.
\end{equation}
Hence, the assertion of (b) follows from (a).\hfill$\square$

\subsection{Proof of Theorem \ref{thm:ConsistencyExpectation}}\label{app:proof:3.1}
\noindent {\it Proof:}  Throughout this proof we denote the terms that are squared in \eqref{Te} and \eqref{Ta} by $\overline{L}_{e,n}(\beta)$ and $\overline{L}_{a,n}(\beta)$, respectively. In the following we will show that
\begin{equation}\label{eqn:ConvergenceToInfinity}
\overline{L}_{j,n}(\beta) \overset{\PP}{\longrightarrow} \infty \quad \text{as} \quad n\to\infty
\end{equation}
for $j\in\{a,e\}$, which implies the assertion.

Let $M\geq 1$ and $f_n:=f$ for $n\in\N$. Recall the definitions of $\widehat{L}_n(\beta)$ and $\widehat{L}_{n,M}(\beta)$ from \eqref{eqn:LnHat} and \eqref{eqn:LnM}. Since $\widehat{L}_n(\beta)$ and $L_n(\beta)$ have the same distribution, we can assume without loss of generality that they are identical. All pairs of points that contribute to $\widehat{L}_{n,M}(\beta)$ also contribute to $\widehat{L}_n(\beta)$ so that $\widehat{L}_{n,M}(\beta)\leq \widehat{L}_n(\beta)$. This implies that, for $j\in\{e,a\}$,
$$
\overline{L}_{j,n}(\beta) \geq \frac{\widehat{L}_{n,M}(\beta) - \BE \widehat{L}_{n,M}(\beta)}{\sqrt{\frac{d\kappa_d}{2(2\beta+d)} }n r_n^{\beta+d/2}} + \frac{\BE \widehat{L}_{n,M}(\beta)- m_{n,j}(\beta)}{\sqrt{\frac{d\kappa_d}{2(2\beta+d)} }n r_n^{\beta+d/2}} =: S_{1,n}+S_{2,j,n}\\
$$
with
$$
m_{n,e}(\beta)=\frac{n(n-1)}{2}\int_{W^2} {\bf 1}\{\|x-y\|\leq r_n\} \, \|x-y\|^\beta\, \dint (x,y)
$$
and
$$
m_{n,a}(\beta)=\frac{d\kappa_d}{2(\beta+d)} \, n(n-1)r_n^{\beta+d}.
$$
Using the same arguments as in the proof of Theorem \ref{thm:ELn}, one can show that
\begin{align*}
\lim_{n\to\infty} \frac{\BE \widehat{L}_{n,M}(\beta)}{n^2 r_n^{\beta+d}} &= \lim_{n\to\infty} \frac{n(n-1)}{2n^2 r_n^{\beta+d}} \int_{W^2} {\bf 1}\{\|x-y\|\leq r_n\} \, \|x-y\|^{\beta} \, f_M(x) \, f_M(y) \, \dint (x,y) \\
& =\frac{d\kappa_d}{2(\beta+d)} \int_W f_M(x)^2 \, \dint x.
\end{align*}
By the Cauchy-Schwarz inequality, we have
$$
\int_W 1 \, \dint x = 1 =\int_W f(x) \, \dint x < \sqrt{\int_W f(x)^2 \, \dint x} \sqrt{\int_W 1 \, \dint x} = \sqrt{\int_W f(x)^2 \, \dint x}
$$
since $f\not\equiv {\bf 1}_W$. Together with the monotone convergence theorem this implies that we can choose $M\geq 1$ such that
$$
\lim_{n\to\infty} \frac{\BE \widehat{L}_{n,M}(\beta)}{n^2 r_n^{\beta+d}}\geq \frac{d\kappa_d}{2(\beta+d)} (1+\varepsilon)
$$
for some $\varepsilon\in (0,\infty)$. Since, by Theorem \ref{thm:ELn},
$$
\lim_{n\to\infty}  \frac{m_{n,j}(\beta)}{n^2 r_n^{\beta+d}} = \frac{d\kappa_d}{2(\beta+d)}
$$
for $j\in\{e,a\}$, this shows that $S_{2,e,n}$ and $S_{2,a,n}$ behave at least as $\frac{\sqrt{d\kappa_d(2\beta+d)}}{\sqrt{2}(\beta+d)}\varepsilon n r_n^{d/2}$ as $n\to\infty$. From the Chebyshev inequality and Lemma \ref{lem:VarianceLnM} it follows that
\begin{align*}
 \lim_{n\to\infty} \P\bigg(|S_{1,n}|\geq \frac{\sqrt{d\kappa_d(2\beta+d)}}{2\sqrt{2}(\beta+d)} \varepsilon n r_n^{d/2} \bigg) 
 & \leq  \lim_{n\to\infty} \frac{\Var \widehat{L}_{n,M}(\beta)}{\frac{d\kappa_d}{2(2\beta+d)} \frac{d\kappa_d (2\beta+d)}{8(\beta+d)^2} n^2 r_n^{2\beta+d} \varepsilon^2 n^2 r_n^d}\\
& = \frac{16(\beta+d)^2}{(d\kappa_d)^2\varepsilon^2} \lim_{n\to\infty} \frac{\sigma^{(1)}_{\beta,f_M} n^2 r_n^{2\beta+d} + \sigma^{(2)}_{\beta,f_M} n^3 r_n^{2\beta+2d} }{n^4 r_n^{2\beta+2d}}\\
& = \frac{16(\beta+d)^2}{(d\kappa_d)^2\varepsilon^2} \lim_{n\to\infty} \frac{\sigma^{(1)}_{\beta,f_M} }{ n^2 r_n^{d}} + \frac{\sigma^{(2)}_{\beta,f_M}}{n} = 0,
\end{align*}
which implies \eqref{eqn:ConvergenceToInfinity} for $j\in\{a,e\}$.
 \hfill$\square$

\subsection{Proofs of Theorem \ref{thm:CloseAlternatives} and Theorem \ref{thm:CloseAlternativesBoundary}}\label{app:proof:4.1u4.2}
We prepare the proofs of Theorem \ref{thm:CloseAlternatives} and Theorem \ref{thm:CloseAlternativesBoundary} with several lemmas. By $\widetilde{L}_n(\beta)$ we denote the statistic $L_n(\beta)$ with respect to i.i.d.\ points $\widetilde{X}_1,\ldots,\widetilde{X}_n$ distributed according to the density $1+a_n g$, while $L_n(\beta)$ is with respect to $n$ i.i.d.\ points uniformly distributed in $W$.

\begin{lemma}\label{lem:MomentBounds}
Assume that $W$ and $g$ satisfy \eqref{eqn:AssumptionWh} and let $n\geq 2$. Then, for any $\beta>-d$,
\begin{equation} \label{eqn:ApproximationExpectation}
\begin{split}
&
 \left| \BE \widetilde{L}_n(\beta) - \BE L_n(\beta) - \frac{n (n-1) a_n^2}{2} \int_{W^2} \mathbf{1}\{\|x-y\|\leq r_n\} \|x-y\|^\beta g(x) g(y) \, \dint(x,y) \right| \\
&
 \leq \frac{d\kappa_d C_{W,g}}{\beta+d} n^2 r_n^{\beta+d+1} a_n.
\end{split}
\end{equation}
Moreover, for any $\beta>-d/2$,
\begin{equation}\label{eqn:BoundVariances}
\begin{split}
 \big| \Var \widetilde{L}_n(\beta) - \Var L_n(\beta) \big|
 \leq C \big( n^2 r_n^{2\beta+d} a_n (a_n+r_n) + n^3 r_n^{2\beta+2d} a_n (a_n+r_n + a_n^2 + a_n^3 + a_n^2 r_n) \big)
\end{split}
\end{equation}
with some constant $C\in(0,\infty)$ depending on $\beta$, $d$, $C_{W,g}$ and $g$.
\end{lemma}

\noindent {\it Proof:}
It follows from \eqref{eqn:ELnexact} in Theorem \ref{thm:ELn} that
\begin{equation}\label{eqn:DifferenceExpectations}
\begin{split}
\BE \widetilde{L}_n(\beta) - \BE L_n(\beta) & = \frac{n(n-1)a_n^2}{2} \int_{W^2} {\bf 1}\{\|x-y\|\leq  r_n\} \|x-y\|^{\beta}  g(x)g(y)\, \dint(x,y) \\
& \quad + n(n-1)a_n  \int_{W^2} {\bf 1}\{\|x-y\|\leq  r_n\} \|x-y\|^{\beta} g(x) \, \dint(x,y).
\end{split}
\end{equation}
We have
\begin{align*}
& \int_{W^2} {\bf 1}\{\|x-y\|\leq r_n\} \|x-y\|^{\beta} g(x) \, \dint(x,y)\\
 & = \frac{d \kappa_d r_n^{\beta+d}}{\beta+d} \int_W g(x) \, \dint x \\
& \quad
 + \int_W {\bf 1}\{d(x,\partial W)\leq r_n\} \bigg(  \int_W {\bf 1}\{\|x-y\|\leq r_n\} \|x-y\|^2 \, \dint y - \frac{d \kappa_d r_n^{\beta+d}}{\beta+d}\bigg) g(x) \, \dint x.
\end{align*}
Here, the first term is zero since $\int_W g(x)  \, \dint x=0$. By \eqref{eqn:AssumptionWh}, the absolute value of the second term can be bounded by
$$
\frac{d \kappa_d }{\beta+d} r_n^{\beta+d}  \int_W {\bf 1}\{d(x,\partial W)\leq r_n\} |g(x)| \, \dint x  \leq \frac{d \kappa_d C_{W,g}}{\beta+d} r_n^{\beta+d+1},
$$
which proves \eqref{eqn:ApproximationExpectation}.

From Theorem \ref{thm:Variance}(a) we can deduce
\begin{align*}
\Var \widetilde{L}_n(\beta) - \Var L_n(\beta)
& = \BE \widetilde{L}_n(2\beta) - \BE L_n(2\beta)\\
& \quad + n(n-1)(n-2) \int_{W^3}  {\bf 1}\{\|x-y_1\|,\|x-y_2\|\leq r_n\} \|x-y_1\|^\beta \|x-y_2\|^\beta\\
& \hskip 4cm \times \big(a_n(2g(y_1)+g(x))+a_n^2(g(y_1)g(y_2)+2g(y_1)g(x))\\
& \hskip 4.5cm+a_n^3 g(y_1)g(y_2)g(x) \big) \, \dint(y_1,y_2,x)\\
& \quad -\frac{4n-6}{n(n-1)} (\BE \widetilde{L}_n(\beta)- \BE L_n(\beta)) (\BE \widetilde{L}_n(\beta)+ \BE L_n(\beta))\\
& =: \bar{R}_{1,n} + \bar{R}_{2,n} - \bar{R}_{3,n}.
\end{align*}
It follows from
\begin{align*}
 \frac{n (n-1) a_n^2}{2} \bigg| \int_{W^2} \mathbf{1}\{\|x-y\|\leq r_n\} \|x-y\|^{2\beta} g(x) g(y) \, \dint(x,y) \bigg|
 \leq \frac{d\kappa_d}{2(2\beta+d)} \int_W g(x)^2 \, \dint x \ n^2 r_n^{2\beta+d} a_n^2
\end{align*}
and \eqref{eqn:ApproximationExpectation} that
$$
|\bar{R}_{1,n}| \leq  \frac{d \kappa_d}{2(2\beta+d)} \int_W g(x)^2 \, \dint x \ n^2 r_n^{2\beta+d} a_n^2 + \frac{d\kappa_d C_{W,g}}{2\beta+d} n^2 r_n^{2\beta+d+1} a_n.
$$
From
\begin{align*}
\BE \widetilde{L}_n(\beta) + \BE L_n(\beta) & \leq \frac{d\kappa_d}{2(\beta+d)} \bigg( 1 + \int_W (1+a_n g(x))^2 \, \dint x \bigg) n^2r_n^{\beta+d} \\
& = \frac{d\kappa_d}{2(\beta+d)} \bigg( 2 + a_n^2 \int_W g(x)^2 \, \dint x \bigg) n^2r_n^{\beta+d},
\end{align*}
\begin{align*}
 \frac{n (n-1) a_n^2}{2} \bigg| \int_{W^2} \mathbf{1}\{\|x-y\|\leq r_n\} \|x-y\|^{\beta} g(x) g(y) \, \dint(x,y) \bigg| 
 \leq \frac{d\kappa_d}{2(\beta+d)} \int_W g(x)^2 \, \dint x \ n^2 r_n^{\beta+d} a_n^2
\end{align*}
and  \eqref{eqn:ApproximationExpectation} we conclude
\begin{align*}
|\bar{R}_{3,n}| & \leq C_3 \frac{1}{n} (n^2r_n^{\beta+d}a_n^2 + n^2r_n^{\beta+d+1}a_n ) (1+a_n^2)  n^2 r_n^{\beta+d}\\
& \leq C_3 n^3 r_n^{2\beta+2d} a_n (a_n + r_n + a_n^3+a_n^2 r_n)
\end{align*}
with some constant $C_3\in(0,\infty)$ depending on $\beta$, $d$, $C_{W,g}$ and $g$.

By similar arguments as for the second term in \eqref{eqn:DifferenceExpectations}, one obtains
\begin{align*}
& n^3 \bigg| \int_{W^3}  {\bf 1}\{\|x-y_1\|,\|x-y_2\|\leq r_n\} \|x-y_1\|^\beta \|x-y_2\|^\beta a_n\big(2g(y_1)+g(x)\big) \, \dint(y_1,y_2,x) \bigg| \\
& \leq \frac{6d^2\kappa_d^2 C_{W,g}}{(\beta+d)^2} n^3 r_n^{2\beta+2d+1} a_n.
\end{align*}
Moreover, one can show the inequality
\begin{align*}
& n^3 \bigg| \int_{W^3}  {\bf 1}\{\|x-y_1\|,\|x-y_2\|\leq r_n\} \|x-y_1\|^\beta \|x-y_2\|^\beta\\
&\hskip 2cm \times \big(a_n^2(g(y_1)g(y_2)+2g(y_1)g(x))+a_n^3 g(y_1)g(y_2)g(x) \big) \, \dint(y_1,y_2,x) \bigg|\\
& \leq \frac{3d^2\kappa_d^2 }{(\beta+d)^2}  \int_W g(x)^2 \, \dint x \ n^3 r_n^{2\beta+2d} a_n^2 + \frac{d^2\kappa_d^2 }{(\beta+d)^2} \int_W |g(x)|^3 \, \dint x \ n^3 r_n^{2\beta+2d} a_n^3.
\end{align*}
Summarising, it follows that
$$
|\bar{R}_{2,n}| \leq C_2 n^3 r_n^{2\beta+2d} a_n (r_n+a_n+a_n^2)
$$
with some constant $C_2\in(0,\infty)$ depending on $\beta$, $d$, $C_{W,g}$ and $g$. Combining the estimates for $\bar{R}_{1,n}$, $\bar{R}_{2,n}$ and $\bar{R}_{3,n}$ completes the proof of \eqref{eqn:BoundVariances}. \hfill$\square$

\begin{lemma}\label{lem:CLTcontiguous}
Let $\beta>-d/2$ and assume that the observation window $W$ satisfies \eqref{eqn:AdditionalAssumptionW}, that $n^2r_n^d\to \infty$ and \linebreak$\max\{n r_n^{d+1},nr_n^d a_n^3\}\to 0$ as $n\to\infty$ and that
\begin{equation}\label{eqn:AssumptionVarianceLtilde}
\lim_{n\to\infty} \frac{\Var \widetilde{L}_n(\beta) -\Var L_n(\beta)}{n^2 r_n^{2\beta+d}}=0.
\end{equation}
Then,
$$
\frac{\widetilde{L}_n(\beta)- \BE \widetilde{L}_n(\beta)}{\sqrt{\frac{d\kappa_d}{2(2\beta+d)}}n r_n^{\beta+d/2}} \overset{\mathcal{D}}{\longrightarrow} N(0,1) \quad \text{as} \quad n\to\infty.
$$
\end{lemma}

We prepare the proof of Lemma \ref{lem:CLTcontiguous} with the following inequality.

\begin{lemma}\label{lem:IntegralInequality}
For $p,q>0$, $v\in L^{p+q}(W)$ and $a>0$,
$$
\int_W \max\{v(x)-a,0\}^p \, \dint x \leq \frac{1}{a^q} \int_W |v(x)|^{p+q} \, \dint x.
$$
\end{lemma}

\noindent {\it Proof:} We have that
$$
\int_W |v(x)|^{p+q} \, \dint x \geq \int_W  \mathbf{1}\{v(x)\geq a\} (v(x)-a)^p a^q \, \dint x = a^q \int_W \max\{v(x)-a,0\}^p \, \dint x,
$$
which is the desired inequality.
\hfill$\square$

\noindent {\it Proof of Lemma \ref{lem:CLTcontiguous}:}
In the following we consider the framework from the Lemmas \ref{lem:VarianceLnaM}, \ref{lem:CLTLnaM} and \ref{lem:ApproximationM} with $f \equiv \mathbf{1}_W$ and $f_n:=\mathbf{1}_W+a_n g$, $n\in\N$. Then, $\widetilde{L}_n(\beta)$ has the same distribution as $\widehat{L}_n(\beta)$. For the latter we will prove convergence to $N(0,1)$ after an appropriate rescaling.

It follows from \eqref{eqn:AssumptionVarianceLtilde} and Theorem \ref{thm:Variance}(c) that
\begin{equation}\label{eqn:VarHatLn}
\lim_{n\to\infty} \frac{\Var \widehat{L}_n(\beta)}{\sigma^2_{\beta,f,n}} = \lim_{n\to\infty} \frac{\Var \widetilde{L}_n(\beta) - \Var L_n(\beta)}{\sigma^2_{\beta,f,n}} + \lim_{n\to\infty} \frac{\Var L_n(\beta)}{\sigma^2_{\beta,f,n}}=1.
\end{equation}
For the rest of this proof we choose $M=2$. Lemma \ref{lem:ApproximationM} yields
\begin{equation}\label{eqn:R1Square}
\begin{split}
& \BE \bigg(\frac{\widehat{L}_{n}(\beta)-\BE \widehat{L}_{n}(\beta)}{\sigma_{\beta,f,n}} - \frac{\widehat{L}_{n,M}(\beta)-\BE \widehat{L}_{n,M}(\beta)}{\sigma_{\beta,f,n}}\bigg)^2\\
& \leq \frac{d\kappa_d}{2\beta+d} \frac{n^2r_n^{2\beta+d}}{\sigma_{\beta,f,n}^2} \int_{W}  \overline{f}_{n,M}(x)^2+M \overline{f}_{n,M}(x) \, \dint x\\
& \quad + \frac{18d^2\kappa_d^2}{(\beta+d)^2} \frac{n^3r_n^{2\beta+2d}}{\sigma_{\beta,f,n}^2} \int_W M^2 \overline{f}_{n,M}(x) + M \overline{f}_{n,M}(x)^2 + \overline{f}_{n,M}(x)^3 \, \dint x.
\end{split}
\end{equation}
It follows from Lemma \ref{lem:IntegralInequality} (with $p=1$, $q=2$ and $p=2$, $q=1$, respectively) that
$$
\int_W \overline{f}_{n,M}(x) \, \dint x  = a_n \int_W \max\{g(x)-1/a_n,0\} \, \dint x  \leq a_n^3 \int_W |g(x)|^3 \, \dint x
$$
and
$$
\int_W \overline{f}_{n,M}(x)^2 \, \dint x = a_n^2 \int_W \max\{g(x)-1/a_n,0\}^2 \, \dint x \leq a_n^3 \int_W |g(x)|^3 \, \dint x.
$$
Moreover, we have
$$
\int_W \overline{f}_{n,M}(x)^3 \, \dint x = a_n^3 \int_W \max\{g(x)-1/a_n,0\}^3 \, \dint x \leq a_n^3 \int_W |g(x)|^3 \, \dint x.
$$
Since $\sigma^2_{\beta,f,n}=\sigma_{\beta,f}^{(1)} n^2 r_n^{2\beta+d}$, the right-hand side of \eqref{eqn:R1Square} is at most of order
$$
\frac{n^2r_n^{2\beta+d}}{\sigma_{\beta,f,n}^2} a_n^3 + \frac{n^3r_n^{2\beta+2d}}{\sigma_{\beta,f,n}^2} a_n^3 = \frac{1 + n r_n^d}{\sigma_{\beta,f}^{(1)}}  a_n^3,
$$
which vanishes as $n\to\infty$. This means that
\begin{equation}\label{eqn:L2ConvergenceHatLnMHatLn}
\lim_{n\to\infty}\BE \bigg(\frac{\widehat{L}_{n}(\beta)-\BE \widehat{L}_{n}(\beta)}{\sigma_{\beta,f,n}} - \frac{\widehat{L}_{n,M}(\beta)-\BE \widehat{L}_{n,M}(\beta)}{\sigma_{\beta,f,n}}\bigg)^2=0.
\end{equation}
Together with \eqref{eqn:VarHatLn} we see that
\begin{equation}\label{eqn:VarHatLnM}
\lim_{n\to\infty} \frac{\Var \widehat{L}_{n,M}(\beta)}{\sigma^2_{\beta,f,n}}=1.
\end{equation}
It follows from Lemma \ref{lem:CLTLnaM}, where the variance condition is satisfied because of \eqref{eqn:VarHatLnM} and $\sigma^2_{\beta,f,n}=\sigma^2_{\beta,f_M,n}$, that
$$
\frac{\widehat{L}_{n,M}(\beta)- \BE \widehat{L}_{n,M}(\beta)}{\sigma_{\beta,f,n}} \overset{\mathcal{D}}{\longrightarrow} N(0,1) \quad \text{as} \quad n\to\infty.
$$
Because of the $L^2$-convergence in \eqref{eqn:L2ConvergenceHatLnMHatLn} this yields
$$
\frac{\widehat{L}_{n}(\beta)- \BE \widehat{L}_{n}(\beta)}{\sigma_{\beta,f,n}} \overset{\mathcal{D}}{\longrightarrow} N(0,1) \quad \text{as} \quad n\to\infty,
$$
which completes the proof.
\hfill$\square$

\noindent {\it Proof of Theorem \ref{thm:CloseAlternatives}:}
By Lemma \ref{lem:MomentBounds} we have that
\begin{equation}\label{eqn:decompositionRT}
\frac{\BE \widetilde{L}_n(\beta) - \BE L_n(\beta)}{\sqrt{\frac{d\kappa_d}{2(2\beta+d)}} n r_n^{\beta+d/2}} = T_n +R_n
\end{equation}
with
$$
T_n := \frac{\sqrt{2\beta+d}}{\sqrt{2d\kappa_d}} \frac{(n-1)a_n^2}{r_n^{\beta+d/2}} \int_{W^2} \mathbf{1}\{\|x-y\|\leq r_n\} \|x-y\|^\beta g(x) g(y) \, \dint(x,y)
$$
and a remainder term $R_n$ satisfying
\begin{equation}\label{eqn:BoundRn}
|R_n|\leq \frac{C_{W,g}\sqrt{2d\kappa_d(2\beta+d) }}{\beta+d} n r_n^{d/2+1} a_n.
\end{equation}
As in the proof of Theorem \ref{thm:ELn} one can show that
\begin{equation}\label{eqn:LimitTn}
\lim_{n\to\infty} \frac{T_n}{n r_n^{d/2} a_n^2} = \frac{\sqrt{d\kappa_d(2\beta+d)}}{\sqrt{2}(\beta+d)} \int_W g(x)^2 \, \dint x.
\end{equation}
For $\gamma=0$ one obtains $\lim_{n\to\infty} T_n=0$ and $\lim_{n\to\infty} R_n=0$. The latter follows from the assumption $\min\{n r_n^{d/2+1}a_n,r_n/a_n \} \to 0$ as $n\to\infty$, whence, by \eqref{eqn:BoundRn}, $R_n$ vanishes directly or is of a lower order than $T_n$ and, thus, also vanishes.

For $\gamma>0$ or $n r_n^{d/2} a_n^2\to\infty$ as $n\to\infty$, we have that $\lim_{n\to\infty} r_n/a_n=0$. Indeed, if there was a subsequence $(n_m)$ such that $r_{n_m}/a_{n_m}\geq c$ for some $c>0$, we would have $n_m r_{n_m}^{d/2+1} a_{n_m}\geq c n_m r_{n_m}^{d/2} a_{n_m}^2$. Then $\min\{n_m r_{n_m}^{d/2+1}a_{n_m},r_{n_m}/a_{n_m} \}$ would not converge to $0$ as $m\to\infty$, which is a contradiction. Because of \eqref{eqn:BoundRn} and \eqref{eqn:LimitTn} it follows from $\lim_{n\to\infty}r_n/a_n = 0$ that
$\lim_{n\to\infty} R_n/T_n=0$, whence $T_n$ is the leading summand in \eqref{eqn:decompositionRT}.

Assume that $n r_n^{d/2}a_n^2\to \gamma\in[0,\infty)$ as $n\to\infty$. By \eqref{eqn:BoundVariances}, we have
\begin{align*}
 \lim_{n\to\infty}\frac{|\Var \widetilde{L}_n(\beta) - \Var L_n(\beta)|}{n^2 r_n^{2\beta+d}} 
 \leq  C \lim_{n\to\infty} a_n (a_n+r_n) + n r_n^{d} a_n (a_n+r_n + a_n^2 + a_n^3 + a_n^2 r_n) =0,
\end{align*}
where we also used that $a_n,r_n, n r_n^{d+1}\to0$ as $n\to\infty$. Now Lemma \ref{lem:CLTcontiguous} implies
$$
\frac{\widetilde{L}_n(\beta)- \BE \widetilde{L}_n(\beta)}{\sqrt{\frac{d\kappa_d}{2(2\beta+d)}}n r_n^{\beta+d/2}} \overset{\mathcal{D}}{\longrightarrow} N(0,1) \quad \text{as} \quad n\to\infty.
$$
This together with \eqref{eqn:decompositionRT} and the above analysis of the asymptotic behaviour of $T_n$ and $R_n$ yields
$$
\frac{\widetilde{L}_n(\beta)- \BE L_n(\beta)}{\sqrt{\frac{d\kappa_d}{2(2\beta+d)}}n r_n^{\beta+d/2}} \overset{\mathcal{D}}{\longrightarrow} N\left(\frac{\sqrt{d\kappa_d(2\beta+d)}}{\sqrt{2}(\beta+d)}\int_W g(x)^2 \, \dint x \ \gamma,1\right) \quad \text{as} \quad n\to\infty.
$$
Now (a) follows from the continuous mapping theorem.

Next we show part (b). It follows from \eqref{eqn:BoundVariances} that
$$
\frac{\Var \widetilde{L}_n(\beta)}{(n^2 r_n^{\beta+d}a_n^2)^2} \leq C \frac{a_n (a_n+r_n) + nr_n^d a_n (a_n+r_n+a_n^2 + a_n^3 + a_n^2 r_n)}{(n r_n^{d/2}a_n^2)^2} + \frac{\Var L_n(\beta)}{(n^2 r_n^{\beta+d}a_n^2)^2}.
$$
The first term on the right-hand side vanishes as $n\to\infty$ since $a_n,r_n,n r_n^{d+1}\to 0$ and $n r_n^{d/2}a_n^2\to\infty$ as $n\to\infty$. Because $\Var L_n(\beta)$ behaves as $n^2 r_n^{2\beta+d}$ by Theorem \ref{thm:Variance}(c), the second term is of order $1/(n r_n^{d/2} a_n^2)^2$ and converges to zero as $n\to\infty$. We thus have
$$
\lim_{n\to\infty} \frac{\Var \widetilde{L}_n(\beta)}{(n^2 r_n^{\beta+d}a_n^2)^2} \ =0 \quad \text{and} \quad \frac{\widetilde{L}_n(\beta) - \BE \widetilde{L}_n(\beta)}{n^2 r_n^{\beta+d}a_n^2}  \overset{\mathbb{P}}{\longrightarrow} 0 \quad \text{as} \quad n\to\infty.
$$
Together with the fact that $T_n$ is the dominating term in \eqref{eqn:decompositionRT} and \eqref{eqn:LimitTn}, this means that
$$
\frac{\widetilde{L}_n(\beta) - \BE L_n(\beta)}{\sqrt{\frac{d\kappa_d}{2(2\beta+d)}}n^2 r_n^{\beta+d}a_n^2}  \overset{\mathbb{P}}{\longrightarrow} \frac{\sqrt{d\kappa_d(2\beta+d)}}{\sqrt{2}(\beta+d)} \int_W g(x)^2 \, \dint x \quad \text{as} \quad n\to\infty.
$$
Because of $n r_n^{d/2} a_n^2\to \infty$ as $n\to\infty$ this implies
$$
\frac{\widetilde{L}_n(\beta) - \BE L_n(\beta)}{\sqrt{\frac{d\kappa_d}{2(2\beta+d)}}n r_n^{\beta+d/2}}  \overset{\mathbb{P}}{\longrightarrow} \infty \quad \text{as} \quad n\to\infty,
$$
which proves part (b).

Part (c) follows from \eqref{eqn:DiffExactAsymptotic} in the proof of Corollary \ref{cor:CLTuniform}. \hfill$\square$

\noindent {\it Proof of Theorem \ref{thm:CloseAlternativesBoundary}:} Without loss of generality we can assume that $r_n<d(\operatorname{supp} g, \partial W)$ for each $n$. Consequently, the assumption \eqref{eqn:AssumptionWh} is satisfied with $C_{W,g}=0$ for $r=r_n$. Now the proof of Theorem \ref{thm:CloseAlternatives} works without the additional assumption that $\min\{n r_n^{d/2+1} a_n,r_n/a_n\}\to 0$ as $n\to\infty$ because $R_n=0$.
\hfill $\square$

\subsection{Proof of Lemma \ref{lem:expect}}\label{app:proof:5.1}
\noindent {\it Proof:} Let $d \in \{2,3\}$, $W=[0,1]^d$ and $r_n\leq 1$. We apply Corollary \ref{cor:Mean} to obtain
\begin{align*}
\BE L_n(\beta) &= \frac{n(n-1)}{2}\int_{W^2}{\bf 1}\{\Vert x-y \Vert \le r_n\}\Vert x-y \Vert^\beta \, \dint (x,y)\\
&= \frac{n(n-1)}{2}\int_{B^d(0,r_n)}\Vert y\Vert^\beta \int_{\R^d}{\bf 1}\{x\in W, \, x-y\in W\} \, \dint x\, \dint y\\
&= \frac{n(n-1)}{2}\int_{B^d(0,r_n)}\Vert y\Vert^\beta \Vol(W\cap (W+y)) \, \dint y \\
&= \frac{n(n-1)}{2}\int_{B^d(0,r_n)}\Vert y\Vert^\beta \prod_{j=1}^d \left( 1-|y_j| \right) \, \dint y,
\end{align*}
with $y=(y_1, \dotsc, y_d) \in \R^d$. The formulae in (a) and (b) follow now from a longer calculation with polar coordinates. \hfill$\square$

\section{A consequence of Lebesgue's differentiation theorem}

\begin{lemma}\label{lem:Differentiation}
Let $g: \R^d\to\R$ be a measurable function with $\|g\|_\infty:=\sup_{y\in\R^d}|g(y)|<\infty$ and let $\beta>-d$. Then, for almost all $x\in\R^d$,
$$
\lim_{r\to 0} \frac{1}{r^{\beta+d}} \int_{B^d(x,r)} \|x-y\|^\beta \, g(y) \, \dint y = \frac{d\kappa_d}{\beta+d} g(x).
$$
\end{lemma}

\noindent {\it Proof:}
We choose $p>1$ subject to $p\beta>-d$. Then for any $x\in\R^d$ and $r>0$,
\begin{align*}
& \bigg| \frac{1}{r^{\beta+d}} \int_{B^d(x,r)} \|x-y\|^\beta g(y) \, \dint y - \frac{d\kappa_d}{\beta+d} g(x) \bigg| \\
& \leq \frac{1}{r^d} \int_{B^d(x,r)} \frac{\|y-x\|^\beta}{r^\beta} \, |g(y)-g(x)| \, \dint y \\
& \leq \bigg( \frac{1}{r^d} \int_{B^d(x,r)} \frac{\|y-x\|^{p\beta}}{r^{p\beta}} \, \dint y  \bigg)^{1/p} \bigg( \frac{1}{r^d} \int_{B^d(x,r)} |g(y)-g(x)|^{p/(p-1)} \, \dint y  \bigg)^{(p-1)/p} \\
& = \bigg( \frac{d\kappa_d}{p\beta+d}  \bigg)^{1/p} \bigg( \frac{1}{r^d} \int_{B^d(x,r)} |g(y)-g(x)|^{p/(p-1)} \, \dint y  \bigg)^{(p-1)/p},
\end{align*}
where we have used the H\"older inequality in the second last step. By Lebesgue's differentiation theorem (see, for example, \cite[Theorem 8.8]{Rudin1974}), we have
$$
\lim_{r\to\infty} \frac{1}{r^d} \int_{B^d(x,r)} |g(y)-g(x)| \, \dint y =0
$$
for almost all $x\in\R^d$. Since $|g(y)-g(x)|^{p/(p-1)}\leq \big( 2\|g\|_\infty \big)^{1/(p-1)} |g(y)-g(x)|$, we have
$$
\lim_{r\to\infty} \frac{1}{r^d} \int_{B^d(x,r)} |g(y)-g(x)|^{p/(p-1)} \, \dint y = 0
$$
for almost all $x\in\R^d$. Together with the above inequalities this proves the statement.
\hfill$\square$

\section{A central limit theorem for a triangular scheme of $U$-statistics}

In the following we provide a central limit theorem for second-order $U$-statistics of a triangular scheme of random vectors, which is a slight generalisation of \cite[Theorem 2.1]{19}.

For each $n\in\N$ let $Y_1^{(n)},\ldots,Y_n^{(n)}$ be i.i.d.\ random vectors in $\R^d$, whose distribution may depend on $n$. We use the shorthand notation $\mathscr{Y}_n=\{Y_1^{(n)},\ldots,Y_n^{(n)}\}$, $n\in\N$, in the sequel. For $n\in\N$ let $h_n: \R^d\times\R^d\to\R$ be a bounded, symmetric and measurable function and let
$$
S_n:=\frac{1}{2} \sum_{(y_1,y_2)\in \mathscr{Y}^2_{n,\neq}} h_n(y_1,y_2).
$$
The random variables $S_n$, $n\in\N$, are so-called second order $U$-statistics. The following theorem provides a sufficient criterion for the convergence of $(S_n)$, after rescaling, to a standard Gaussian random variable.

\begin{theorem}\label{thm:CLTUstatistics}
Let $S_n$, $n\in\N$, be as above. Assume that $\Var S_n>0$ for all $n\in\N$ and let $\sigma_n>0$, $n\in\N$, be such that
$\lim_{n\to\infty} \Var S_n/\sigma_n^2 =1$.
If
$$
\lim_{n\to\infty} \frac{1}{\sigma_n}\sup_{y_1,y_2\in\R^d} |h_n(y_1,y_2)|=0 \quad \text{and} \quad
\lim_{n\to\infty} \frac{n}{\sigma_n} \sup_{y\in\R^d} \BE |h_n(y,Y_1^{(n)})|=0,
$$
then
$$
\frac{S_n-\BE S_n}{\sigma_{n}} \overset{\mathcal{D}}{\longrightarrow} N\left(0,1\right) \quad \text{as} \quad n\to\infty.
$$
\end{theorem}

\noindent {\it Proof:}
In the special case that $(Y_i^{(n)})_{1\leq i\leq n<\infty}$ are identically distributed, this is a slightly re-written version of \cite[Theorem 2.1]{19}. Otherwise, there are measurable maps $T_n: [0,1]\to\R^d$, $n\in\N$, such that $Y_i^{(n)}$, $i\in\{1,\ldots,n\}$, has the same distribution as $T_n(U)$, where $U$ is a uniformly distributed random variable on $[0,1]$ (see, for example, the proof of Theorem 29.6 in \cite{Billingsley1979}). For $n\in\N$ define $\tilde{h}_n: [0,1]^2\ni(u_1,u_2)\mapsto h_n(T_n(u_1),T_n(u_2))$ and let $\mathcal{U}_n:=\{U_1,\ldots,U_n\}$, where $U_1,\ldots,U_n$ are independent and uniformly distributed on $[0,1]$. Then, $S_n$ has the same distribution as
$$
\widetilde{S}_n:=\frac{1}{2} \sum_{(u_1,u_2)\in\mathcal{U}^2_{n\,\neq}} \tilde{h}_n(u_1,u_2).
$$
Since the assumptions of the theorem are satisfied for the $U$-statistics $(S_n)$, they must also hold for the $U$-statistics $(\widetilde{S}_n)$. As the underlying random variables of $(\widetilde{S}_n)$ are identically distributed, we are in the previously discussed special case for which the central limit theorem holds. This completes the proof.
\hfill $\square$


\begin{table}[h]
\centering
\begin{tabular}{c|c|cccc|cccc}
Alt. & $n$ & ${BR}_n^2(0.1)$ & ${BR}_n^2(0.25)$ & ${BR}_n^2(0.5)$ & ${BR}_n^2(1)$ & $NN_{n,1}^{(0.5)}$ & $NN_{n,15}^{(0.5)}$ & DB & MS \\\hline
\multirow{4}{*}{CON}
& 50  & \textbf{74} & 40 & 33 & 6  & 16 & 66 & 31 & 6\\
& 100 & \textbf{96} & 66 & 56 & 9  & 19 & 90 & 58 & 14 \\
& 200 & \textbf{*}  & 91 & 83 & 14 & 25 & 98 & 89 & 25 \\
& 500 & \textbf{*}  & \textbf{*}  & 99 & 36 & 41  & \textbf{*}  &  \textbf{*} & 41\\\hline
\multirow{4}{*}{CLU}
& 50  & \textbf{80} & 34 & 31 & 42 & 78 & 67 & 28 & 36\\
& 100 & 73 & 30 & 27 & 41 & 74 & \textbf{82} & 28 & 48 \\
& 200 & 61 & 26 & 24 & 41  & 58 & \textbf{90} & 28 & 52 \\
& 500 & 45 & 23 & 22 & 41 & 32 & \textbf{96} & 29 & 47\\\hline
\multirow{4}{*}{$H_0$}
& 50  & 5 & 5 & 5 & 5 & 5 & 5 & 4 & 2\\
& 100 & 5 & 5 & 5 & 5 & 5 & 5 & 5 & 3\\
& 200 & 5 & 5 & 5 & 5 & 4 & 5 & 5 & 4\\
& 500 & 5 & 5 & 5 & 5 & 4 & 5 & 5 & 5\\\hline
\end{tabular}
\caption{Empirical rejection rates of the different competitors ($d=2$)}\label{tab:Comp.d.2}
\end{table}

\begin{table}[h]
\centering
\small
\begin{tabular}{c|c|c|ccccccccccccc}
Alt. & $\beta$ & $n \backslash k$ & 1 & 2 & 3 & 4 & 5 & 6 & 7 & 8 & 9 & 10 & 15 & 20 & 25 \\
\hline
\multirow{4}{*}{CON} & \multirow{12}{*}{$-0.5$} & 50 & 39 & 54 & 61 & 66 & 69 & 71 & \textbf{72} & 72 & 72 & 72 & 67 & 60 & 51 \\
 & & 100 & 59 & 77 & 85 & 90 & 92 & 94 & 95 & 95 & \textbf{96} & 96 & 96 & 96 & 95 \\
 & & 200 & 82 & 95 & 98 & 99 & 99 & \textbf{*} & * & * & * & * & * & * & * \\
 & & 500 & 99 & \textbf{*} & * & * & * & * & * & * & * & * & * & * & * \\
\cline{3-16}
\multirow{4}{*}{CLU} & & 50 & 95 & \textbf{97} & 97 & 96 & 94 & 91 & 88 & 86 & 83 & 80 & 68 & 59 & 53 \\
 & & 100 & 91 & 96 & \textbf{97} & 97 & 97 & 96 & 95 & 94 & 93 & 92 & 82 & 73 & 65 \\
 & & 200 & 81 & 92 & 95 & \textbf{96} & 96 & 96 & 96 & 96 & 96 & 95 & 91 & 86 & 79 \\
 & & 500 & 59 & 77 & 85 & 89 & 91 & 92 & 93 & \textbf{94} & 94 & 94 & 94 & 92 & 90 \\
\cline{3-16}
\multirow{4}{*}{$H_0$} & & 50 & 5 & 5 & 5 & 5 & 5 & 5 & 5 & 5 & 5 & 5 & 5 & 5 & 5 \\
 & & 100 & 5 & 5 & 5 & 5 & 5 & 5 & 5 & 5 & 5 & 5 & 5 & 6 & 5 \\
 & & 200 & 5 & 5 & 5 & 5 & 5 & 5 & 5 & 5 & 5 & 5 & 5 & 5 & 5 \\
 & & 500 & 5 & 5 & 5 & 5 & 5 & 5 & 5 & 5 & 5 & 5 & 5 & 5 & 5 \\
 \hline
\multirow{4}{*}{CON} & \multirow{12}{*}{0} & 50 & 41 & 57 & 64 & 68 & 70 & \textbf{71} & 71 & 71 & 71 & 69 & 59 & 45 & 34 \\
 & & 100 & 64 & 80 & 88 & 91 & 93 & 94 & 95 & \textbf{96} & 96 & 96 & 96 & 94 & 92 \\
 & & 200 & 85 & 96 & 99 & \textbf{*} & * & * & * & * & * & * & * & * & * \\
 & & 500 & \textbf{*} & * & * & * & * & * & * & * & * & * & * & * & * \\
 \cline{3-16}
\multirow{4}{*}{CLU} & & 50 & 95 & \textbf{96} & 95 & 91 & 87 & 81 & 75 & 68 & 64 & 59 & 43 & 35 & 31 \\
 & & 100 & 92 & \textbf{96} & 96 & 96 & 95 & 93 & 92 & 89 & 86 & 82 & 64 & 51 & 43 \\
 & & 200 & 83 & 92 & 94 & \textbf{96} & 95 & 95 & 95 & 94 & 93 & 92 & 84 & 74 & 63 \\
 & & 500 & 61 & 80 & 86 & 89 & 91 & 92 & \textbf{93} & 93 & 93 & 93 & 91 & 88 & 84 \\
 \cline{3-16}
\multirow{4}{*}{$H_0$} & & 50 & 4 & 6 & 5 & 5 & 5 & 5 & 5 & 5 & 5 & 5 & 5 & 5 & 5 \\
 & & 100 & 5 & 5 & 6 & 5 & 5 & 5 & 5 & 5 & 5 & 5 & 5 & 5 & 5 \\
 & & 200 & 5 & 5 & 5 & 5 & 5 & 5 & 5 & 5 & 5 & 5 & 5 & 5 & 5 \\
 & & 500 & 5 & 5 & 5 & 5 & 5 & 4 & 5 & 5 & 5 & 5 & 5 & 5 & 5 \\
\hline
\multirow{4}{*}{CON} & \multirow{12}{*}{1} & 50 & 40 & 53 & 59 & \textbf{63} & 64 & 64 & 64 & 63 & 61 & 58 & 40 & 25 & 16 \\
 & & 100 & 60 & 77 & 85 & 89 & 91 & 93 & \textbf{94} & 94 & 94 & 94 & 93 & 89 & 82 \\
 & & 200 & 83 & 96 & 98 & 99 & 99 & \textbf{*} & * & * & * & * & * & * & * \\
 & & 500 & \textbf{*} & * & * & * & * & * & * & * & * & * & * & * & * \\
\cline{3-16}
\multirow{4}{*}{CLU} & & 50 & \textbf{91} & 91 & 86 & 78 & 67 & 56 & 47 & 41 & 37 & 34 & 29 & 28 & 29 \\
 & & 100 & 88 & \textbf{92} & 92 & 91 & 88 & 85 & 80 & 74 & 68 & 63 & 40 & 32 & 29 \\
 & & 200 & 79 & 88 & \textbf{91} & 91 & 91 & 90 & 89 & 88 & 86 & 83 & 69 & 54 & 42 \\
 & & 500 & 56 & 75 & 81 & 85 & 87 & 88 & \textbf{89} & 89 & 89 & 89 & 85 & 79 & 73 \\
 \cline{3-16}
\multirow{4}{*}{$H_0$} & & 50 & 5 & 5 & 5 & 5 & 5 & 5 & 5 & 5 & 5 & 5 & 5 & 5 & 5 \\
 & & 100 & 5 & 5 & 5 & 5 & 5 & 5 & 6 & 5 & 5 & 5 & 5 & 5 & 5 \\
 & & 200 & 5 & 5 & 5 & 5 & 5 & 5 & 5 & 5 & 5 & 5 & 5 & 5 & 5 \\
 & & 500 & 5 & 5 & 5 & 5 & 5 & 4 & 5 & 5 & 5 & 5 & 5 & 4 & 5 \\
\hline
\multirow{4}{*}{CON} & \multirow{12}{*}{5} & 50 & 29 & 38 & 43 & \textbf{45} & 44 & 43 & 41 & 38 & 35 & 31 & 15 & 8 & 8 \\
 & & 100 & 44 & 61 & 71 & 77 & 80 & 81 & 83 & 83 & \textbf{84} & 83 & 77 & 64 & 45 \\
 & & 200 & 65 & 86 & 93 & 96 & 97 & 98 & \textbf{99} & 99 & 99 & 99 & 99 & 99 & 99 \\
 & & 500 & 97 & \textbf{*} & * & * & * & * & * & * & * & * & * & * & * \\
\cline{3-16}
\multirow{4}{*}{CLU} & & 50 & \textbf{72} & 68 & 56 & 42 & 30 & 24 & 24 & 24 & 25 & 27 & 29 & 30 & 30 \\
 & & 100 & \textbf{69} & 75 & 73 & 69 & 62 & 55 & 48 & 41 & 35 & 30 & 25 & 26 & 27 \\
 & & 200 & 57 & 69 & 73 & \textbf{74} & 72 & 70 & 68 & 65 & 61 & 57 & 39 & 29 & 26 \\
 & & 500 & 36 & 53 & 61 & 65 & 69 & 70 & \textbf{71} & 71 & 71 & 71 & 65 & 57 & 50 \\
 \cline{3-16}
\multirow{4}{*}{$H_0$} & & 50 & 5 & 5 & 5 & 5 & 5 & 5 & 5 & 5 & 5 & 5 & 5 & 5 & 5 \\
 & & 100 & 5 & 5 & 5 & 5 & 5 & 5 & 5 & 5 & 5 & 5 & 5 & 5 & 5 \\
 & & 200 & 5 & 5 & 5 & 5 & 5 & 5 & 5 & 5 & 5 & 5 & 5 & 5 & 5 \\
 & & 500 & 5 & 5 & 5 & 5 & 5 & 5 & 5 & 5 & 5 & 5 & 5 & 5 & 5 \\
\hline
 \end{tabular}
\caption{Empirical rejection rates for $T_e$ in case $d=2$ and $r_n =\Big(\frac{k}{n\kappa_d}\Big)^{\frac1d}$, see \eqref{eq:rn.1}}\label{tab:Te.d.2}
\end{table}

\begin{table}[h]
\centering
\small
\begin{tabular}{c|c|c|ccccccccccccc}
Alt. & $\beta$ & $n \backslash k$ & 1 & 2 & 3 & 4 & 5 & 6 & 7 & 8 & 9 & 10 & 15 & 20 & 25 \\
\hline
\multirow{4}{*}{CON} & \multirow{12}{*}{$-0.5$} & 50 & 16 & 18 & 22 & 25 & 27 & 30 & 31 & 33 & 34 & 35 & \textbf{36} & 35 & 32 \\
 & & 100 & 17 & 22 & 27 & 32 & 36 & 40 & 43 & 46 & 48 & 50 & 57 & 60 & \textbf{61} \\
 & & 200 & 18 & 26 & 34 & 41 & 45 & 51 & 56 & 60 & 63 & 66 & 77 & 82 & \textbf{85} \\
 & & 500 & 24 & 36 & 47 & 56 & 64 & 71 & 76 & 80 & 83 & 86 & 94 & 97 & \textbf{99} \\
\cline{3-16}
\multirow{4}{*}{CLU} & & 50 & 65 & 79 & 85 & 89 & 91 & 92 & \textbf{93} & 93 & 93 & 93 & 92 & 88 & 83 \\
 & & 100 & 41 & 55 & 66 & 72 & 77 & 81 & 83 & 85 & 87 & 87 & 89 & \textbf{90} & 88 \\
 & & 200 & 22 & 31 & 38 & 45 & 50 & 55 & 59 & 62 & 65 & 68 & 75 & 78 & \textbf{80} \\
 & & 500 & 11 & 14 & 17 & 20 & 22 & 25 & 26 & 29 & 30 & 32 & 40 & 46 & \textbf{50} \\
  \cline{3-16}
\multirow{4}{*}{$H_0$} & & 50 & 5 & 5 & 5 & 5 & 5 & 5 & 5 & 5 & 4 & 5 & 5 & 5& 5 \\
 & & 100 & 5 & 5 & 5 & 5 & 5 & 5 & 5 & 5 & 5 & 5 & 5 & 5 & 5 \\
 & & 200 & 5 & 5 & 5 & 5 & 5 & 5 & 5 & 5 & 5 & 5 & 5 & 5 & 5 \\
 & & 500 & 5 & 5 & 5 & 5 & 5 & 5 & 5 & 5 & 5 & 5 & 5 & 5 & 5 \\
 \hline
\multirow{4}{*}{CON} & \multirow{12}{*}{0} & 50 & 11 & 15 & 23 & 29 & 27 & 32 & 34 & 32 & 34 & \textbf{36} & 36 & 32 & 27 \\
 & & 100 & 15 & 20 & 28 & 33 & 37 & 41 & 43 & 46 & 47 & 49 & 57 & 59 & \textbf{60} \\
 & & 200 & 16 & 27 & 33 & 43 & 50 & 56 & 61 & 64 & 68 & 71 & 78 & 82 & \textbf{86} \\
 & & 500 & 27 & 39 & 51 & 60 & 67 & 75 & 79 & 83 & 85 & 88 & 95 & 98 & \textbf{99} \\
 \cline{3-16}
\multirow{4}{*}{CLU} & & 50 & 59 & 76 & 86 & 91 & 90 & 92 & \textbf{93} & 92 & 92 & 92 & 89 & 81 & 71 \\
 & & 100 & 38 & 53 & 66 & 74 & 78 & 81 & 83 & 84 & 85 & 85 & \textbf{88} & 87 & 84 \\
 & & 200 & 20 & 31 & 37 & 47 & 54 & 59 & 63 & 66 & 68 & 70 & 75 & 77 & \textbf{78} \\
 & & 500 & 12 & 14 & 18 & 21 & 23 & 27 & 28 & 31 & 31 & 34 & 41 & 46 & \textbf{51} \\
  \cline{3-16}
\multirow{4}{*}{$H_0$} & & 50 & 3 & 3 & 4 & 5 & 3 & 4 & 5 & 5 & 6 & 7 & 4 & 4 & 4 \\
 & & 100 & 4 & 3 & 5 & 5 & 6 & 5 & 4 & 5 & 5 & 5 & 5 & 6 & 5 \\
 & & 200 & 5 & 6 & 4 & 5 & 4 & 5 & 5 & 5 & 5 & 5 & 5 & 5 & 5 \\
 & & 500 & 7 & 6 & 4 & 6 & 5 & 6 & 5 & 6 & 5 & 5 & 5 & 5 & 5 \\
\hline
\multirow{4}{*}{CON} & \multirow{12}{*}{1} & 50 & 17 & 19 & 21 & 24 & 27 & 29 & \textbf{30} & 30 & 30 & 30 & 29 & 24 & 20 \\
 & & 100 & 17 & 22 & 27 & 32 & 35 & 39 & 41 & 44 & 46 & 48 & 52 & \textbf{53} & 53 \\
 & & 200 & 18 & 26 & 33 & 40 & 46 & 51 & 55 & 59 & 62 & 65 & 75 & 79 & \textbf{83} \\
 & & 500 & 24 & 36 & 47 & 55 & 63 & 70 & 75 & 79 & 82 & 85 & 93 & 97 & \textbf{98} \\
\cline{3-16}
\multirow{4}{*}{CLU} & & 50 & 63 & 75 & 81 & 84 & 85 & \textbf{86} & 86 & 86 & 85 & 84 & 76 & 63 & 48 \\
 & & 100 & 39 & 54 & 63 & 69 & 73 & 76 & 78 & 80 & \textbf{81} & 81 & 81 & 79 & 74 \\
 & & 200 & 21 & 30 & 37 & 43 & 48 & 52 & 56 & 59 & 61 & 63 & 69 & \textbf{70} & 70 \\
 & & 500 & 11 & 13 & 16 & 19 & 22 & 23 & 25 & 27 & 28 & 30 & 37 & 41 & \textbf{43} \\
  \cline{3-16}
\multirow{4}{*}{$H_0$} & & 50 & 8 & 5 & 5 & 5 & 5 & 5 & 5 & 5 & 5 & 5 & 5 & 5 & 5 \\
 & & 100 & 6 & 5 & 5 & 5 & 5 & 5 & 5 & 5 & 5 & 5 & 5 & 5 & 5 \\
 & & 200 & 5 & 5 & 5 & 5 & 5 & 5 & 5 & 5 & 5 & 5 & 5 & 5 & 5 \\
 & & 500 & 5 & 5 & 5 & 5 & 5 & 5 & 5 & 5 & 5 & 5 & 5 & 5 & 5 \\
\hline
\multirow{4}{*}{CON} & \multirow{12}{*}{5} & 50 & 13 & 15 & 17 & 18 & 20 & \textbf{22} & 21 & 22 & 21 & 21 & 19 & 14 & 10 \\
 & & 100 & 15 & 17 & 20 & 22 & 25 & 27 & 29 & 31 & 33 & 34 & 37 & \textbf{38} & 37 \\
 & & 200 & 16 & 19 & 23 & 27 & 32 & 36 & 39 & 41 & 44 & 47 & 56 & 61 & \textbf{65} \\
 & & 500 & 18 & 24 & 31 & 38 & 43 & 49 & 54 & 58 & 63 & 65 & 78 & 86 & \textbf{90} \\
\cline{3-16}
\multirow{4}{*}{CLU} & & 50 & 45 & 57 & 61 & 64 & \textbf{65} & 65 & 65 & 64 & 62 & 60 & 46 & 30 & 18 \\
 & & 100 & 30 & 38 & 44 & 48 & 52 & 55 & 57 & 58 & \textbf{60} & 60 & 58 & 54 & 48 \\
 & & 200 & 17 & 22 & 25 & 29 & 32 & 35 & 37 & 39 & 41 & 43 & \textbf{47} & 47 & 47 \\
 & & 500 & 10 & 11 & 12 & 14 & 15 & 16 & 17 & 18 & 19 & 19 & 23 & \textbf{25} & 25 \\
  \cline{3-16}
\multirow{4}{*}{$H_0$} & & 50 & 5 & 5 & 5 & 5 & 5 & 5 & 5 & 5 & 5 & 5 & 5 & 5 & 5 \\
 & & 100 & 6 & 5 & 5 & 6 & 5 & 5 & 5 & 5 & 5 & 5 & 5 & 5 & 5 \\
 & & 200 & 5 & 5 & 5 & 5 & 5 & 5 & 5 & 5 & 5 & 5 & 5 & 5 & 5 \\
 & & 500 & 5 & 5 & 5 & 5 & 5 & 5 & 5 & 5 & 5 & 5 & 5 & 5 & 5 \\
\hline
 \end{tabular}
\caption{Empirical rejection rates for $T_a$ in case $d=2$ and $r_n = \Big(\frac{k}{n^{\frac32}\kappa_d}\Big)^{\frac1d}$, see \eqref{eq:rn.2}}\label{tab:Ta.d.2}
\end{table}


\begin{table}[h]
\centering
\small
\begin{tabular}{c|c|c|ccccccccccccc|c}
Alt. & $\beta$ & $n \backslash k$ & 1 & 2 & 3 & 4 & 5 & 6 & 7 & 8 & 9 & 10 & 15 & 20 & 25 & DB \\
\hline
\multirow{4}{*}{CON} & \multirow{12}{*}{$-0.5$} & 50 & 92 & 95 & \textbf{96} & 96 & 96 & 96 & 95 & 95 & 94 & 93 & 88 & 82 & 75 & 59 \\
 & & 100 & \textbf{*} & * & * & * & * & * & * & * & * & * & * & * & * & 89\\
 & & 200 & \textbf{*} & * & * & * & * & * & * & * & * & * & * & * & * & *\\
 & & 500 & \textbf{*} & * & * & * & * & * & * & * & * & * & * & * & * & * \\
\cline{3-17}
\multirow{4}{*}{CLU} & & 50 & \textbf{*} & 99 & 98 & 96 & 94 & 92 & 89 & 86 & 83 & 80 & 67 & 58 & 52 & 22\\
 & & 100 & \textbf{*} & * & * & * & 99 & 98 & 97 & 96 & 94 & 92 & 81 & 71 & 62 & 22\\
 & & 200 & \textbf{*} & * & * & * & * & * & * & 99 & 99 & 99 & 93 & 85 & 77 & 22\\
 & & 500 &\textbf{*} & * & * & * & * & * & * & * & * & * & * & 98 & 94 & 23 \\
\cline{3-17}
\multirow{4}{*}{$H_0$} & & 50 & 5 & 5 & 5 & 5 & 5 & 5 & 5 & 5 & 5 & 5 & 5 & 5 & 5 & 5 \\
 & & 100 & 5 & 5 & 5 & 5 & 5 & 5 & 5 & 5 & 5 & 5 & 5 & 5 & 5 & 4 \\
 & & 200 & 5 & 5 & 5 & 5 & 5 & 5 & 5 & 5 & 5 & 5 & 5 & 5 & 5 & 5 \\
 & & 500 & 5 & 5 & 5 & 5 & 5 & 5 & 5 & 5 & 5 & 5 & 5 & 5 & 5 & 5 \\
 \cline{1-17}
\multirow{4}{*}{CON} & \multirow{12}{*}{0} & 50 & 92 & 95 & \textbf{96} & 95 & 95 & 95 & 94 & 92 & 91 & 90 & 80 & 66 & 54 \\
 & & 100 &\textbf{*} & * & * & * & * & * & * & * & * & * & * & * & 99 \\
 & & 200 &\textbf{*} & * & * & * & * & * & * & * & * & * & * & * & * \\
 & & 500 &\textbf{*} & * & * & * & * & * & * & * & * & * & * & * & * \\
 \cline{3-16}
\multirow{4}{*}{CLU} & & 50 & \textbf{99} & 97 & 92 & 84 & 77 & 70 & 64 & 58 & 54 & 50 & 38 & 32 & 30 \\
 & & 100 & \textbf{*} & * & 99 & 98 & 95 & 90 & 85 & 80 & 75 & 70 & 53 & 43 & 36 \\
 & & 200 & \textbf{*} & * & * & * & * & 99 & 98 & 96 & 94 & 91 & 75 & 61 & 51 \\
 & & 500 &\textbf{*} & * & * & * & * & * & * & * & * & * & 97 & 89 & 79 \\
\cline{3-16}
\multirow{4}{*}{$H_0$} & & 50 & 5 & 5 & 5 & 4 & 5 & 5 & 5 & 5 & 5 & 5 & 5 & 5 & 5 \\
 & & 100 & 5 & 5 & 5 & 5 & 5 & 5 & 5 & 5 & 5 & 5 & 5 & 5 & 5 \\
 & & 200 & 5 & 5 & 5 & 4 & 5 & 5 & 5 & 5 & 5 & 5 & 5 & 5 & 5 \\
 & & 500 & 5 & 5 & 5 & 5 & 6 & 5 & 5 & 5 & 5 & 5 & 5 & 5 & 5 \\
\cline{1-16}
\multirow{4}{*}{CON} & \multirow{12}{*}{1} & 50 & 91 & \textbf{94} & 94 & 93 & 92 & 90 & 88 & 84 & 81 & 78 & 54 & 35 & 23 \\
 & & 100 & 99 & \textbf{*} & * & * & * & * & * & * & * & * & 99 & 97 & 92 \\
 & & 200 &\textbf{*} & * & * & * & * & * & * & * & * & * & * & * & * \\
 & & 500 &\textbf{*} & * & * & * & * & * & * & * & * & * & * & * & * \\
\cline{3-16}
\multirow{4}{*}{CLU} & & 50 & \textbf{98} & 83 & 60 & 47 & 39 & 34 & 31 & 29 & 29 & 28 & 26 & 26 & 25 \\
 & & 100 & \textbf{*} & 99 & 94 & 80 & 68 & 57 &  49 & 43 & 39 & 36 & 28 & 26 & 25 \\
 & & 200 & \textbf{*} & * & * & 99 & 97 & 91 & 83 & 75 & 68 & 61 & 41 & 32 & 29 \\
 & & 500 & \textbf{*} & * & * & * & * & * & * & 99 & 99 & 98 & 82 & 62 & 48 \\
\cline{3-16}
\multirow{4}{*}{$H_0$} & & 50 & 5 & 5 & 5 & 5 & 5 & 5 & 5 & 5 & 5 & 5 & 5 & 5 & 5 \\
 & & 100 & 5 & 5 & 5 & 5 & 5 & 5 & 5 & 5 & 5 & 5 & 5 & 5 & 5 \\
 & & 200 & 5 & 5 & 5 & 5 & 5 & 5 & 5 & 5 & 5 & 5 & 5 & 5 & 5 \\
 & & 500 & 5 & 5 & 5 & 5 & 6 & 5 & 5 & 5 & 5 & 5 & 5 & 5 & 5 \\
\cline{1-16}
\multirow{4}{*}{CON} & \multirow{12}{*}{5} & 50 & 82 & \textbf{84} & 82 & 78 & 72 & 65 & 57 & 49 & 43 & 35 & 14 & 8 & 8 \\
 & & 100 & 98 & 99 & 99 & 99 & \textbf{*} & 99 & 99 & 98 & 98 & 97 & 86 & 63 & 40 \\
 & & 200 &\textbf{*} & * & * & * & * & * & * & * & * & * & * & * & * \\
 & & 500 &\textbf{*} & * & * & * & * & * & * & * & * & * & * & * & * \\
\cline{3-16}
\multirow{4}{*}{CLU} & & 50 & \textbf{75} & 24 & 23 & 24 & 25 & 25 & 25 & 26 & 26 & 27 & 27 & 27 & 27 \\
 & & 100 & \textbf{98} & 75 & 39 & 25 & 24 & 24 & 24 & 23 & 24 & 24 & 25 & 25 & 25 \\
 & & 200 & \textbf{*} & 99 & 91 & 73 & 53 & 38 & 30 & 27 & 24 & 24 & 24 & 24 & 25 \\
 & & 500 & \textbf{*} & * & * & 99 & 98 & 96 & 91 & 84 & 75 & 66 & 33 & 25 & 23 \\
\cline{3-16}
\multirow{4}{*}{$H_0$} & & 50 & 5 & 5 & 5 & 5 & 5 & 5 & 5 & 5 & 5 & 5 & 5 & 5 & 5 \\
 & & 100 & 5 & 5 & 5 & 5 & 5 & 5 & 5 & 5 & 5 & 5 & 5 & 5 & 5 \\
 & & 200 & 5 & 5 & 5 & 5 & 5 & 5 & 5 & 5 & 5 & 5 & 5 & 5 & 5 \\
 & & 500 & 5 & 5 & 5 & 5 & 5 & 5 & 5 & 5 & 5 & 5 & 5 & 5 & 5 \\
\cline{1-16}
 \end{tabular}
\caption{Empirical rejection rates for $T_e$ in case $d=3$ and $r_n =\Big(\frac{k}{n\kappa_d}\Big)^{\frac1d}$, see \eqref{eq:rn.1}}\label{tab:Te.d.3}
\end{table}

\begin{table}[h]
\centering
\small
\begin{tabular}{c|c|c|ccccccccccccc|c}
Alt. & $\beta$ & $n \backslash k$ & 1 & 2 & 3 & 4 & 5 & 6 & 7 & 8 & 9 & 10 & 15 & 20 & 25 & DB\\
\hline
\multirow{4}{*}{CON} & \multirow{12}{*}{$-0.5$} & 50 & 58 & 69 & 74 & 77 & \textbf{78} & 78 & 78 & 78 & 77 & 76 & 69 & 61 & 51 & 59\\
 & & 100 & 74 & 86 & 91 & 94 & 95 & 96 & 96 & \textbf{97} & 97 & 97 & 97 & 96 & 95 & 89 \\
 & & 200 & 86 & 96 & 99 & 99 & \textbf{*} & * & * & * & * & * & * & * & * & *\\
 & & 500 & 97 & \textbf{*} & * & * & * & * & * & * & * & * & * & * & * & *\\
\cline{3-17}
\multirow{4}{*}{CLU} & & 50 & 99 & \textbf{*} & * & 99 & 99 & 99 & 99 & 98 & 98 & 97 & 89 & 76 & 58 & 22\\
 & & 100 & 99 & \textbf{*} & * & * & * & * & * & * & * & * & 99 & 98 & 95 & 22\\
 & & 200 & 97 & \textbf{*} & * & * & * & * & * & * & * & * & * & * & * & 22\\
 & & 500 & 74 & 91 & 96 & 98 & 99 & 99 & 99 & \textbf{*} & * & * & * & * & * & 23\\
 \cline{3-17}
 \multirow{4}{*}{$H_0$} & & 50 & 3 & 5 & 5 & 5 & 5 & 5 & 5 & 5 & 5 & 5 & 5 & 5 & 5 & 5 \\
 & & 100 & 5 & 5 & 5 & 5 & 5 & 5 & 5 & 5 & 5 & 5 & 5 & 5 & 5 & 4 \\
 & & 200 & 5 & 5 & 5 & 5 & 5 & 5 & 5 & 5 & 5 & 5 & 5 & 5 & 5 & 5 \\
 & & 500 & 5 & 5 & 5 & 5 & 5 & 5 & 5 & 5 & 5 & 5 & 5 & 5 & 5 & 5 \\
 \cline{1-17}
\multirow{4}{*}{CON} & \multirow{12}{*}{0} & 50 & 62 & 72 & 75 & 76 & 76 & \textbf{78} & 77 & 75 & 74 & 72 & 66 & 54 & 42 \\
 & & 100 & 70 & 88 & 91 & 94 & 95 & 96 & 96 & \textbf{97} & 97 & 96 & 96 & 95 & 94 \\
 & & 200 & 85 & 96 & 99 & \textbf{*} & * & * & * & * & * & * & * & * & * \\
 & & 500 & 98 & \textbf{*} & * & * & * & * & * & * & * & * & * & * & * \\
 \cline{3-16}
\multirow{4}{*}{CLU} & & 50 & 99 & \textbf{*} & 99 & 99 & 99 & 99 & 98 & 97 & 95 & 93 & 73 & 45 & 26 \\
 & & 100 & 99 & \textbf{*} & * & * & * & * & * & * & * & * & 99 & 94 & 84 \\
 & & 200 & 96 & \textbf{*} & * & * & * & * & * & * & * & * & * & * & 99 \\
 & & 500 & 76 & 91 & 96 & 98 & 99 & 99 & 99 & \textbf{*} & * & * & * & * & * \\
  \cline{3-16}
 \multirow{4}{*}{$H_0$} & & 50 & 4 & 4 & 3 & 6 & 3 & 8 & 4 & 6 & 3 & 3 & 7 & 6 & 6 \\
 & & 100 & 3 & 8 & 3 & 7 & 5 & 7 & 6 & 7 & 5 & 6 & 5 & 4 & 6 \\
 & & 200 & 5 & 6 & 4 & 5 & 6 & 3 & 6 & 4 & 6 & 5 & 4 & 4 & 5 \\
 & & 500 & 7 & 6 & 5 & 5 & 7 & 6 & 5 & 6 & 5 & 5 & 6 & 5 & 4 \\
 \cline{1-16}
\multirow{4}{*}{CON} & \multirow{12}{*}{1} & 50 & 57 & 67 & 71 & \textbf{73} & 73 & 73 & 72 & 70 & 69 & 67 & 54 & 40 & 26 \\
 & & 100 & 72 & 85 & 90 & 93 & 94 & 95 & 95 & \textbf{96} & 96 & 95 & 95 & 93 & 91 \\
 & & 200 & 85 & 96 & 98 & 99 & \textbf{*} & * & * & * & * & * & * & * & * \\
 & & 500 & 97 & \textbf{*} & * & * & * & * & * & * & * & * & * & * & * \\
\cline{3-16}
\multirow{4}{*}{CLU} & & 50 & \textbf{99} & 99 & 99 & 98 & 97 & 96 & 93 & 88 & 81 & 72 & 29 & 12 & 6 \\
 & & 100 & 99 & \textbf{*} & * & * & * & * & * & 99 & 99 & 99 & 94 & 75 & 48 \\
 & & 200 & 96 & 99 & \textbf{*} & * & * & * & * & * & * & * & * & 99 & 96 \\
 & & 500 & 72 & 89 & 94 & 97 & 98 & 99 & 99 & 99 & 99 & 99 & \textbf{*} & 99 & 99 \\
 \cline{3-16}
 \multirow{4}{*}{$H_0$} & & 50 & 7 & 5 & 5 & 5 & 5 &5 & 5 & 5 & 5 & 5 & 5 & 5 & 5 \\
 & & 100 & 5 & 5 & 5 & 5 & 5 & 5 & 5 & 5 & 5 & 5 & 5 & 5 & 5 \\
 & & 200 & 5 & 5 & 5 & 5 & 5 & 5 & 5 & 5 & 5 & 5 & 5 & 5 & 5 \\
 & & 500 & 5 & 5 & 5 & 5 & 5 & 5 & 5 & 5 & 4 & 5 & 5 & 5 & 5 \\
 \cline{1-16}
\multirow{4}{*}{CON} & \multirow{12}{*}{5} & 50 & 50 & 55 & 59 & \textbf{60} & 60 & 59 & 57 & 54 & 52 & 48 & 30 & 16 & 7 \\
 & & 100 & 61 & 75 & 81 & 85 & 87 & 88 & \textbf{89} & 89 & 89 & 89 & 87 & 82 & 76 \\
 & & 200 & 74 & 90 & 95 & 97 & 98 & 99 & 99 & 99 & 99 & \textbf{*} & * & * & * \\
 & & 500 & 90 & 99 & \textbf{*} & * & * & * & * & * & * & * & * & * & * \\
\cline{3-16}
\multirow{4}{*}{CLU} & & 50 & \textbf{96} & 96 & 93 & 89 & 79 & 66 & 49 & 34 & 22 & 15 & 9 & 14 & 18 \\
 & & 100 & 95 & 97 & \textbf{98} & 98 & 97 & 96 & 95 & 93 & 90 & 85 & 50 & 18 & 5 \\
 & & 200 & 86 & 95 & 97 & \textbf{98} & 98 & 98 & 98 & 98 & 98 & 97 & 93 & 83 & 66 \\
 & & 500 & 55 & 74 & 82 & 87 & 90 & 92 & 93 & 94 & 94 & \textbf{95} & 94 & 93 & 89 \\
 \cline{3-16}
 \multirow{4}{*}{$H_0$} & & 50 & 5 & 5 & 5 & 5 & 5 & 5 & 5 & 5 & 5 & 5 & 5 & 5 & 5 \\
 & & 100 & 5 & 5 & 5 & 5 & 5 & 5 & 5 & 5 & 5 & 5 & 5 & 5 & 5 \\
 & & 200 & 5 & 5 & 5 & 5 & 5 & 5 & 5 & 5 & 5 & 5 & 5 & 5 & 5 \\
 & & 500 & 5 & 5 & 5 & 5 & 5 & 5 & 5 & 5 & 5 & 5 & 5 & 5 & 5 \\
 \cline{1-16}
 \end{tabular}
\caption{Empirical rejection rates for $T_a$ in case $d=3$ and $r_n = \Big(\frac{k}{n^{\frac32}\kappa_d}\Big)^{\frac1d}$, see \eqref{eq:rn.2}}\label{tab:Ta.d.3}
\end{table}

\begin{table}[h]
\renewcommand{\arraystretch}{.99}
\centering
\small
\begin{tabular}{c|c|c|ccccccccccccc}
$d$ & $\beta$ & $n \backslash k$ & 1 & 2 & 3 & 4 & 5 & 6 & 7 & 8 & 9 & 10 & 15 & 20 & 25 \\
\hline
\multirow{16}{*}{$2$} & \multirow{4}{*}{$-0.5$} & 50 & \textbf{26}&23&21&20&19&19&19&18&18&18&18&18&18 \\
& & 100 &\textbf{51}&45&42&40&38&37&35&34&33&32&30&28&28 \\
& & 200 &\textbf{86}&83&80&77&75&73&72&70&68&67&61&56&52 \\
& & 500 &\textbf{*}&*&*&*&*&*&*&*&*&99&99&98&97 \\
\cline{2-16}
& \multirow{4}{*}{$0$} & 50 & \textbf{15}&13&12&11&11&11&11&11&11&11&12&12&13 \\
& & 100 & \textbf{40}&29&26&22&19&18&17&17&16&15&15&15&16\\
& & 200 & \textbf{81}&70&62&57&50&46&43&40&37&36&29&26&23\\
& & 500 & \textbf{*}&*&*&99&99&98&97&96&95&94&86&79&72\\
\cline{2-16}
& \multirow{4}{*}{$1$} & 50 & 8&7&8&8&8&8&9&9&9&10&11&11&\textbf{12}\\
& & 100 & \textbf{19}&12&11&11&10&10&10&10&10&10&11&12&12\\
& & 200 & \textbf{64}&37&25&20&18&16&16&15&15&15&14&14&13\\
& & 500 & \textbf{*}&99&96&89&80&71&63&56&51&46&33&27&25\\
\cline{2-16}
& \multirow{4}{*}{$5$} & 50 & 6&7&7&8&7&9&8&9&9&\textbf{10}&10&10&10\\
& & 100 & 8&8&8&9&8&9&9&9&9&9&11&12&\textbf{13}\\
& & 200 & \textbf{20}&10&9&10&11&10&12&12&12&11&12&13&13\\
& & 500 & \textbf{97}&71&36&22&18&16&16&17&16&17&17&17&17\\
\hline
\hline
\multirow{16}{*}{$3$} & \multirow{4}{*}{$-0.5$} & 50 & \textbf{34}&30&28&26&25&24&24&24&24&23&23&24&25 \\
& & 100 &\textbf{67}&60&55&52&50&48&46&44&44&43&40&38&36 \\
& & 200 &\textbf{95}&92&30&88&86&84&83&81&80&79&71&67&64\\
& & 500 &\textbf{*}&*&*&*&*&*&*&*&*&*&*&99&99 \\
\cline{2-16}
& \multirow{4}{*}{$0$} & 50 &16&14&14&12&14&14&14&14&15&15&16&17&\textbf{18} \\
& & 100 &\textbf{43}&31&27&26&24&23&22&23&22&22&22&22&23\\
& & 200 &\textbf{85}&74&65&59&55&51&48&46&44&41&37&35&33\\
& & 500 &\textbf{*}&*&*&99&99&98&97&96&95&94&88&82&77\\
\cline{2-16}
& \multirow{4}{*}{$1$} & 50 &8&9&10&9&10&11&12&12&13&13&14&15&\textbf{17}\\
& & 100 &13&12&12&12&12&13&13&14&14&15&16&17&\textbf{19}\\
& & 200 &\textbf{34}&23&20&18&18&17&18&18&18&18&20&21&22\\
& & 500 &\textbf{98}&85&69&59&52&47&43&41&39&37&34&33&32\\
\cline{2-16}
& \multirow{4}{*}{$5$} & 50 &7&8&9&9&10&11&12&12&12&12&13&16&\textbf{18}\\
& & 100 &9&10&10&11&11&12&13&13&14&14&16&17&\textbf{19}\\
& & 200 &11&13&13&13&14&14&15&15&15&16&18&20&\textbf{22}\\
& & 500 &18&19&19&20&22&21&22&21&22&22&24&25&\textbf{27}\\
\hline
\end{tabular}
\caption{Empirical rejection rates for $T_e$ under the SPS-alternative with $r_n =\Big(\frac{k}{n\kappa_d}\Big)^{\frac1d}$, see \eqref{eq:rn.1}}\label{tab:SPS.Te}
\end{table}


\begin{table}[h]
\renewcommand{\arraystretch}{.99}
\centering
\small
\begin{tabular}{c|c|c|ccccccccccccc}
$d$ & $\beta$ & $n \backslash k$ & 1 & 2 & 3 & 4 & 5 & 6 & 7 & 8 & 9 & 10 & 15 & 20 & 25 \\
\hline
\multirow{16}{*}{$2$} & \multirow{4}{*}{$-0.5$} & 50 & \textbf{34}&30&28&26&24&22&21&20&19&18&14&11&8\\
& & 100 &56&\textbf{58}&58&56&54&53&51&49&47&46&39&33&28\\
& & 200 &75&83&86&87&\textbf{88}&88&88&88&87&87&83&80&76\\
& & 500 &93&98&99&99&99&\textbf{*}&*&*&*&*&*&*&* \\
\cline{2-16}
& \multirow{4}{*}{$0$} & 50 & \textbf{26}&21&20&19&14&13&13&11&10&10&6&5&4 \\
& & 100 & 53&\textbf{54}&53&50&47&43&39&35&32&30&22&16&11\\
& & 200 & 74&83&85&\textbf{87}&87&87&87&86&85&83&76&67&59\\
& & 500 & 95&98&99&99&\textbf{*}&*&*&*&*&*&*&*&*\\
\cline{2-16}
& \multirow{4}{*}{$1$} & 50 & \textbf{27}&16&11&9&8&7&7&6&6&6&6&5&5\\
& & 100 & \textbf{52}&49&44&37&31&25&21&18&15&13&8&6&5\\
& & 200 & 74&81&\textbf{82}&82&82&80&78&76&73&69&52&36&24\\
& & 500 & 93&98&99&99&99&99&\textbf{*}&*&*&*&*&*&99\\
\cline{2-16}
& \multirow{4}{*}{$5$} & 50 & \textbf{13}&8&6&6&6&6&7&6&6&6&6&6&6\\
& & 100 & \textbf{39}&30&21&14&11&9&7&7&7&7&7&8&7\\
& & 200 & 62&\textbf{67}&67&64&60&55&48&41&36&31&13&7&7\\
& & 500 & 84&93&96&97&\textbf{98}&98&98&98&98&98&98&96&93\\
\hline
\hline
\multirow{16}{*}{$3$} & \multirow{4}{*}{$-0.5$} & 50 & \textbf{46}&39&34&29&25&23&20&18&16&14&9&5&4 \\
& & 100 &\textbf{82}&75&71&66&62&59&55&52&49&45&33&23&16 \\
& & 200 &\textbf{99}&98&97&97&96&95&94&93&92&91&84&76&69\\
& & 500 &\textbf{*}&*&*&*&*&*&*&*&*&*&*&*&* \\
\cline{2-16}
& \multirow{4}{*}{$0$} & 50 &\textbf{36}&24&17&14&9&11&7&7&5&4&4&3&3 \\
& & 100 &\textbf{71}&64&52&47&38&34&27&24&20&17&8&4&3\\
& & 200 &\textbf{98}&97&95&93&90&87&84&80&78&74&55&39&27\\
& & 500 &\textbf{*}&*&*&*&*&*&*&*&*&*&*&*&99\\
\cline{2-16}
& \multirow{4}{*}{$1$} & 50 &\textbf{14}&8&6&6&5&5&5&5&5&5&5&5&4\\
& & 100 &\textbf{51}&29&18&13&10&7&6&6&5&5&4&4&4\\
& & 200 &\textbf{96}&89&78&67&56&46&37&31&24&20&7&3&3\\
& & 500 &\textbf{*}&*&*&*&*&*&*&*&*&*&97&88&73\\
\cline{2-16}
& \multirow{4}{*}{$5$} & 50 &\textbf{6}&6&6&6&6&6&6&6&6&6&5&5&5\\
& & 100 &\textbf{9}&7&8&8&8&7&8&7&7&7&7&6&6\\
& & 200 &\textbf{67}&23&10&9&8&9&8&9&9&9&9&8&8\\
& & 500 &\textbf{*}&*&99&97&91&79&63&47&33&23&9&8&9\\
\hline
\end{tabular}
\caption{Empirical rejection rates for $T_a$ under the SPS-alternative with $r_n = \Big(\frac{k}{n^{\frac32}\kappa_d}\Big)^{\frac1d}$, see \eqref{eq:rn.2}}\label{tab:SPS.Ta}
\end{table}

\newpage

\begin{sidewaystable}[h]
\renewcommand{\arraystretch}{.99}
\begin{tabular}{c|c|ccccccccccccc}
$\beta$ & $n \backslash k$ & 1 & 2 & 3 & 4 & 5 & 6 & 7 & 8 & 9 & 10 & 15 & 20 & 25\\
\hline
\multirow{4}{*}{$-0.5$} & 50 & 3.5153 & 3.5611 & 3.647 & 3.7462 & 3.8737 & 4.0189 & 4.1848 & 4.329 & 4.507 & 4.669 & 5.6829 & 6.5529 & 7.326\\
& 100 & 3.684 & 3.71 & 3.8253 & 3.9398 & 4.081 & 4.2427 & 4.3712 & 4.5251 & 4.7462 & 4.9095 & 5.9538 & 7.0291 & 8.1751\\
& 200 & 3.7486 & 3.7786 & 3.8392 & 3.9334 & 4.0411 & 4.1677 & 4.3276 & 4.4797 & 4.6413 & 4.8084 & 5.7686 & 6.7913 & 7.9071\\
& 500 & 3.824 & 3.8526 & 3.9011 & 3.9687 & 4.0495 & 4.1732 & 4.2981 & 4.3739 & 4.5089 & 4.638 & 5.347 & 6.1883 & 7.1429\\
\hline
\multirow{4}{*}{0} & 50 & 3.3375 & 3.5945 & 3.8075 & 3.951 & 4.0252 & 4.3326 & 4.4754 & 4.5129 & 4.9256 & 5.1809 & 6.3216 & 7.3246 & 8.203\\
& 100 & 3.8339 & 3.7618 & 3.7106 & 4.0433 & 4.193 & 4.3351 & 4.6114 & 4.8324 & 5.0195 & 5.2954 & 6.641 & 8.0996 & 9.3564\\
& 200 & 3.6694 & 3.7618 & 3.9031 & 3.9109 & 4.1479 & 4.2694 & 4.5334 & 4.7475 & 4.9236 & 5.1882 & 6.347 & 7.6525 & 9.1298\\
& 500 & 3.7991 & 3.8698 & 3.8935 & 4.0124 & 4.1025 & 4.3069 & 4.4478 & 4.5799 & 4.6927 & 4.8885 & 5.8603 & 6.9065 & 8.198\\
\hline
\multirow{4}{*}{1} & 50 & 3.4684 & 3.5603 & 3.6979 & 3.8496 & 4.0063 & 4.2318 & 4.4481 & 4.6313 & 4.8555 & 5.0735 & 6.2101 & 7.1905 & 7.7949\\
& 100 & 3.6467 & 3.7174 & 3.8564 & 4.0146 & 4.1474 & 4.3562 & 4.5674 & 4.8293 & 5.0666 & 5.3019 & 6.7023 & 7.997 & 9.3607\\
& 200 & 3.7272 & 3.7692 & 3.8592 & 3.9879 & 4.188 & 4.3389 & 4.5259 & 4.7382 & 4.9674 & 5.2089 & 6.4348 & 7.7738 & 9.2293\\
& 500 & 3.7983 & 3.8088 & 3.921 & 4.0282 & 4.1471 & 4.3119 & 4.4313 & 4.5962 & 4.7598 & 4.9432 & 5.9055 & 7.0258 & 8.3193\\
\hline
\multirow{4}{*}{5} & 50 & 3.4067 & 3.396 & 3.5065 & 3.5612 & 3.6555 & 3.7396 & 3.8386 & 3.9407 & 4.0372 & 4.1538 & 4.7196 & 5.0936 & 5.2229\\
& 100 & 3.5517 & 3.5891 & 3.6528 & 3.7325 & 3.8056 & 3.9412 & 4.058 & 4.2367 & 4.3675 & 4.4953 & 5.2748 & 6.0161 & 6.712\\
& 200 & 3.6891 & 3.665 & 3.7209 & 3.7833 & 3.9296 & 4.0139 & 4.1062 & 4.2567 & 4.3563 & 4.5128 & 5.2359 & 6.0929 & 6.9251\\
& 500 & 3.7482 & 3.7563 & 3.8228 & 3.8894 & 3.8786 & 4.0446 & 4.1364 & 4.2289 & 4.3353 & 4.4072 & 5.0608 & 5.7461 & 6.596\\
\hline
\end{tabular}
\caption{Empirical critical values for $T_e$ in case $d=2$ and $r_n = \Big(\frac{k}{n\kappa_d}\Big)^{\frac1d}$, see \eqref{eq:rn.1}}\label{tab:EC2}
\end{sidewaystable}

\newpage

\begin{sidewaystable}[h]
\renewcommand{\arraystretch}{.99}
\begin{tabular}{c|c|ccccccccccccc}
$\beta$ & $n \backslash k$ & 1 & 2 & 3 & 4 & 5 & 6 & 7 & 8 & 9 & 10 & 15 & 20 & 25\\
\hline
\multirow{4}{*}{$-0.5$} & 50 & 3.059 & 3.4495 & 3.4895 & 3.5476 & 3.6013 & 3.6136 & 3.671 & 3.7402 & 3.8304 & 3.8987 & 4.4778 & 5.2045 & 6.1158\\
& 100 & 3.3636 & 3.6174 & 3.6356 & 3.6774 & 3.6665 & 3.6884 & 3.724 & 3.773 & 3.8006 & 3.8531 & 4.1379 & 4.5676 & 5.1179\\
& 200 & 3.6098 & 3.6475 & 3.6837 & 3.7082 & 3.728 & 3.7441 & 3.7669 & 3.7857 & 3.7967 & 3.8159 & 4.0084 & 4.2671 & 4.5564\\
& 500 & 3.6948 & 3.7277 & 3.7479 & 3.7491 & 3.7536 & 3.7463 & 3.7814 & 3.7909 & 3.837 & 3.8487 & 3.8895 & 4.0106 & 4.122\\
\hline
\multirow{4}{*}{0} & 50 & 3.5348 & 3.6357 & 3.8551 & 3.6055 & 3.9197 & 3.6414 & 3.8381 & 4.0619 & 3.9305 & 3.8376 & 4.8105 & 5.8258 & 6.8581\\
& 100 & 3.2805 & 3.721 & 3.4082 & 3.872 & 3.8025 & 3.8163 & 3.8778 & 3.844 & 4.0201 & 4.205 & 4.4408 & 4.84 & 5.7245\\
& 200 & 3.5862 & 3.5359 & 3.728 & 3.6373 & 3.5344 & 3.8531 & 3.8196 & 3.8274 & 3.8627 & 3.7838 & 3.976 & 4.521 & 4.7926\\
& 500 & 3.3917 & 3.8812 & 3.9251 & 3.568 & 3.913 & 3.7914 & 3.7388 & 3.7294 & 3.8093 & 3.7881 & 3.9171 & 4.0679 & 4.1481\\
\hline
\multirow{4}{*}{1} & 50 & 3.0182 & 3.4862 & 3.5104 & 3.5512 & 3.6001 & 3.6774 & 3.7808 & 3.9017 & 4.0562 & 4.1931 & 5.0547 & 6.1696 & 7.4639\\
& 100 & 3.5107 & 3.5938 & 3.6393 & 3.6181 & 3.6692 & 3.7113 & 3.7541 & 3.8079 & 3.8998 & 3.968 & 4.4981 & 5.1274 & 5.9715\\
& 200 & 3.6403 & 3.7021 & 3.7011 & 3.7726 & 3.7307 & 3.7443 & 3.8137 & 3.835 & 3.8984 & 3.8807 & 4.1639 & 4.5612 & 4.9969\\
& 500 & 3.682 & 3.7521 & 3.7587 & 3.781 & 3.7749 & 3.7912 & 3.8095 & 3.8242 & 3.8517 & 3.8635 & 3.931 & 4.1156 & 4.3126\\
\hline
\multirow{4}{*}{5} & 50 & 3.2332 & 3.167 & 3.331 & 3.4236 & 3.4419 & 3.5016 & 3.5726 & 3.6696 & 3.8171 & 3.9049 & 4.5556 & 5.4372 & 6.4344\\
& 100 & 3.154 & 3.4374 & 3.5352 & 3.5248 & 3.6007 & 3.6177 & 3.6573 & 3.6886 & 3.7474 & 3.8571 & 4.2155 & 4.6772 & 5.3163\\
& 200 & 3.3493 & 3.5723 & 3.6399 & 3.6894 & 3.6745 & 3.6561 & 3.7296 & 3.7904 & 3.7717 & 3.7667 & 4.0233 & 4.2971 & 4.6231\\
& 500 & 3.5735 & 3.6706 & 3.7452 & 3.7067 & 3.7459 & 3.7472 & 3.7381 & 3.7855 & 3.7764 & 3.813 & 3.8827 & 3.9822 & 4.1625\\
\hline
\end{tabular}
\caption{Empirical critical values for $T_a$ in case $d=2$ and $r_n = \Big(\frac{k}{n^{\frac32}\kappa_d}\Big)^{\frac1d}$, see \eqref{eq:rn.2}}\label{tab:EC1}
\end{sidewaystable}

\newpage

\begin{sidewaystable}[h]
\renewcommand{\arraystretch}{.99}
\begin{tabular}{c|c|ccccccccccccc}
$\beta$ & $n \backslash k$ & 1 & 2 & 3 & 4 & 5 & 6 & 7 & 8 & 9 & 10 & 15 & 20 & 25\\
\hline
\multirow{4}{*}{$-0.5$} & 50 & 3.291 & 3.4152 & 3.6171 & 3.8073 & 3.9982 & 4.179 & 4.3982 & 4.5908 & 4.7811 & 4.9581 & 5.9169 & 6.687 & 7.2964\\
& 100 & 3.5113 & 3.6402 & 3.8774 & 4.0916 & 4.3269 & 4.5864 & 4.8626 & 5.1119 & 5.3263 & 5.6027 & 6.8917 & 8.1342 & 9.2808\\
& 200 & 3.6032 & 3.7398 & 3.9406 & 4.1741 & 4.414 & 4.715 & 4.9566 & 5.2263 & 5.4921 & 5.7925 & 7.3647 & 8.778 & 10.19\\
& 500 & 3.7113 & 3.8272 & 4.0474 & 4.2397 & 4.4686 & 4.7166 & 5.0007 & 5.2743 & 5.5108 & 5.7829 & 7.2489 & 8.8147 & 10.3807\\
\hline
\multirow{4}{*}{0} & 50 & 3.1545 & 3.3134 & 3.4272 & 3.9481 & 4.0208 & 4.3308 & 4.4546 & 4.5199 & 4.7144 & 4.9395 & 5.8986 & 6.6639 & 7.391\\
& 100 & 3.586 & 3.7246 & 3.9569 & 4.1844 & 4.2736 & 4.6567 & 4.9779 & 5.1929 & 5.5419 & 5.8016 & 7.1195 & 8.4343 & 9.559\\
& 200 & 3.5915 & 3.7587 & 3.9103 & 4.2356 & 4.439 & 4.7096 & 5.1048 & 5.3459 & 5.6768 & 5.9974 & 7.6054 & 9.1797 & 10.6337\\
& 500 & 3.6335 & 3.7929 & 4.0581 & 4.2893 & 4.5389 & 4.7968 & 5.1122 & 5.3755 & 5.6666 & 5.944 & 7.5246 & 9.2364 & 10.8763\\
\hline
\multirow{4}{*}{1} & 50 & 3.1978 & 3.2985 & 3.4761 & 3.7027 & 3.8968 & 4.0746 & 4.282 & 4.4805 & 4.6267 & 4.797 & 5.7106 & 6.2775 & 6.6809\\
& 100 & 3.4313 & 3.5603 & 3.7796 & 4.0027 & 4.2631 & 4.5141 & 4.7796 & 5.0626 & 5.314 & 5.5933 & 6.8797 & 8.0612 & 9.1192\\
& 200 & 3.5185 & 3.6853 & 3.8858 & 4.1833 & 4.4413 & 4.7079 & 4.9798 & 5.2699 & 5.5453 & 5.8156 & 7.4537 & 8.887 & 10.2664\\
& 500 & 3.6457 & 3.7906 & 4.0316 & 4.2263 & 4.4733 & 4.721 & 5.0232 & 5.335 & 5.5669 & 5.9183 & 7.4302 & 9.0472 & 10.723\\
\hline
\multirow{4}{*}{5} & 50 & 3.0234 & 3.0631 & 3.1429 & 3.2483 & 3.3453 & 3.4346 & 3.5348 & 3.6421 & 3.7032 & 3.7805 & 4.2232 & 4.4042 & 4.4337\\
& 100 & 3.2868 & 3.3212 & 3.4543 & 3.6047 & 3.7154 & 3.8567 & 4.0369 & 4.1956 & 4.3376 & 4.5078 & 5.277 & 5.9176 & 6.4958\\
& 200 & 3.3949 & 3.4875 & 3.6135 & 3.7985 & 3.9709 & 4.13 & 4.2746 & 4.4308 & 4.6418 & 4.8225 & 5.8187 & 6.7304 & 7.6005\\
& 500 & 3.5393 & 3.6181 & 3.7635 & 3.8794 & 4.0342 & 4.2151 & 4.3972 & 4.6339 & 4.7776 & 4.9825 & 5.9785 & 7.0771 & 8.2448\\
\hline
\end{tabular}
\caption{Empirical critical values for $T_e$ in case $d=3$ and $r_n = \Big(\frac{k}{n\kappa_d}\Big)^{\frac1d}$, see \eqref{eq:rn.1}}\label{tab:EC4}
\end{sidewaystable}
\newpage
\begin{sidewaystable}[h]
\renewcommand{\arraystretch}{.99}
\begin{tabular}{c|c|ccccccccccccc}
$\beta$ & $n \backslash k$ & 1 & 2 & 3 & 4 & 5 & 6 & 7 & 8 & 9 & 10 & 15 & 20 & 25\\
\hline
\multirow{4}{*}{$-0.5$} & 50 & 3.2597 & 3.4875 & 3.7513 & 4.0104 & 4.3327 & 4.6857 & 5.0597 & 5.4407 & 5.8605 & 6.3121 & 8.7125 & 11.3174 & 14.1969\\
& 100 & 3.2199 & 3.5981 & 3.7745 & 4.0337 & 4.2606 & 4.5426 & 4.8688 & 5.2008 & 5.5566 & 5.9118 & 7.8884 & 10.063 & 12.3647\\
& 200 & 3.6179 & 3.6867 & 3.7927 & 3.9873 & 4.1672 & 4.3922 & 4.6249 & 4.8978 & 5.1767 & 5.4394 & 7.0396 & 8.7485 & 10.5914\\
& 500 & 3.6869 & 3.7585 & 3.8528 & 3.9549 & 4.1134 & 4.2601 & 4.4272 & 4.6144 & 4.8305 & 5.0434 & 6.1219 & 7.3977 & 8.7644\\
\hline
\multirow{4}{*}{0} & 50 & 3.3955 & 3.4367 & 3.8551 & 4.3677 & 4.918 & 4.5172 & 5.1173 & 5.7192 & 6.3221 & 6.9259 & 9.1035 & 12.138 & 15.1724\\
& 100 & 3.2805 & 3.481 & 4.1082 & 3.872 & 4.6225 & 4.563 & 5.3235 & 5.329 & 6.0867 & 6.125 & 8.5008 & 10.89 & 13.2845\\
& 200 & 3.5862 & 3.5359 & 3.9098 & 4.1676 & 4.195 & 4.7623 & 4.6985 & 5.282 & 5.2746 & 5.8607 & 7.6772 & 9.5313 & 11.4006\\
& 500 & 3.3917 & 3.8812 & 3.9251 & 4.1553 & 3.913 & 4.2818 & 4.6642 & 4.5909 & 4.9962 & 5.3339 & 6.2479 & 7.9489 & 9.2872\\
\hline
\multirow{4}{*}{1} & 50 & 3.1833 & 3.4018 & 3.7233 & 4.1124 & 4.4785 & 4.855 & 5.31 & 5.7959 & 6.2698 & 6.7732 & 9.603 & 12.6136 & 15.9305\\
& 100 & 3.2867 & 3.5733 & 3.8195 & 4.1181 & 4.3869 & 4.7105 & 5.0822 & 5.4852 & 5.8748 & 6.327 & 8.5753 & 11.1363 & 13.8202\\
& 200 & 3.4985 & 3.6399 & 3.8368 & 4.0294 & 4.2826 & 4.5365 & 4.8227 & 5.1292 & 5.4417 & 5.7801 & 7.6549 & 9.6059 & 11.761\\
& 500 & 3.6849 & 3.7417 & 3.8808 & 3.9683 & 4.1648 & 4.3516 & 4.5512 & 4.7898 & 5.0224 & 5.2462 & 6.5491 & 8.0747 & 9.6426\\
\hline
\multirow{4}{*}{5} & 50 & 2.4588 & 3.2078 & 3.4256 & 3.7305 & 4.0282 & 4.3659 & 4.7462 & 5.1398 & 5.5282 & 5.9654 & 8.2712 & 10.6489 & 13.3085\\
& 100 & 2.9862 & 3.3926 & 3.6124 & 3.8096 & 4.0335 & 4.3441 & 4.6368 & 4.9501 & 5.2748 & 5.6764 & 7.5187 & 9.5585 & 11.7511\\
& 200 & 3.3645 & 3.5133 & 3.6722 & 3.8172 & 4.0221 & 4.255 & 4.4754 & 4.7343 & 4.9896 & 5.2368 & 6.8046 & 8.4616 & 10.1664\\
& 500 & 3.5962 & 3.6743 & 3.7579 & 3.859 & 4.0147 & 4.1344 & 4.3064 & 4.4773 & 4.7088 & 4.8687 & 5.9855 & 7.2618 & 8.4907\\
\hline
\end{tabular}
\caption{Empirical critical values for $T_a$ in case $d=3$ and $r_n = \Big(\frac{k}{n^{\frac32}\kappa_d}\Big)^{\frac1d}$, see \eqref{eq:rn.2}}\label{tab:EC3}
\end{sidewaystable}

\end{appendices}

\bibliographystyle{abbrv}

\end{document}